\pdfoutput=1
\documentclass[11pt]{article}
\usepackage[left=1in,right=1in,top=1in,bottom=1in]{geometry}
\usepackage{times}
\usepackage{expl3}
\usepackage{cite}
\usepackage[table]{xcolor}
\usepackage{multirow}
\usepackage{stackengine} 
\usepackage{hhline}
\usepackage{lipsum}
\usepackage{titlesec}
\usepackage{wrapfig}
\usepackage{enumerate}
\usepackage{epsfig}
\usepackage{amsmath}
\usepackage{tabularx}
\usepackage{array}
\usepackage{booktabs}
\usepackage{enumitem}
\usepackage{bbm}
\usepackage{calc}
\usepackage{graphicx}
\usepackage{amsmath}
\usepackage[title]{appendix}
\usepackage{amssymb}
\usepackage{epstopdf}
\usepackage{boldline}
\usepackage{arydshln}
\usepackage{calligra}
\usepackage{bm}
\usepackage{url}
\usepackage{blindtext}
\usepackage{accents}

\newcommand{\define}{\stackrel{\mbox{\tiny def}}{=}}

\newtheorem{theorem}{Theorem}
\newtheorem{proposition}{Proposition}

\newtheorem{lemma}{Lemma}

\newtheorem{example}{Example}

\usepackage{mathtools}
\usepackage{epstopdf}
\usepackage{balance}
\usepackage{thmtools}
\usepackage{thm-restate}
\usepackage{hyperref}
\usepackage{cleveref}
\usepackage[mathscr]{euscript}

\usepackage[ruled,vlined]{algorithm2e}
\include{pythonlisting}
\newcommand{\stirlingii}{\genfrac{\{}{\}}{0pt}{}}
\newcommand{\ostar}{\mathbin{\mathpalette\make@circled\star}}

\makeatletter
\newcommand{\removelatexerror}{\let\@latex@error\@gobble}
\makeatother
\makeatletter
\newcommand*{\rom}[1]{\expandafter\@slowromancap\romannumeral #1@}
\makeatother

\ExplSyntaxOn
\newcommand\latinabbrev[1]{
  \peek_meaning:NTF . {% Same as \@ifnextchar
    #1\@}%
  { \peek_catcode:NTF a {% Check whether next char has same catcode as \'a, i.e., is a letter
      #1.\@ }%
    {#1.\@}}}
\ExplSyntaxOff

\titleclass{\subsubsubsection}{straight}[\subsubsection]

\begin{document}
\vspace{1cm}
\title{Operator Characterization via Projectors and Nilpotents}
\vspace{1.8cm}
\author{Shih-Yu~Chang
% <-this % stops a space
\thanks{Shih-Yu Chang is with the Department of Applied Data Science,
San Jose State University, San Jose, CA, U. S. A. (e-mail: {\tt
shihyu.chang@sjsu.edu}). 
           }}

\maketitle

\begin{abstract}
This paper presents significant contributions to the study of operators with countable, continuous, and hybrid spectra, with applications across both finite-dimensional and infinite-dimensional contexts, particularly in non-Hermitian (non-self-adjoint) systems. For finite-dimensional operators, a new concept of analogous matrices is introduced, where two matrices are deemed analogous if they share the same projector and nilpotent structures. This perspective highlights structural equivalences beyond simple spectral similarities. To aid in classifying such matrices, a graph-based representation of these projectors and nilpotent elements is developed, offering both visual and computational insights. Additionally, the paper calculates the number of distinct families of analogous matrices based on matrix size, providing a fundamental tool for matrix classification. The paper also extends the spectral mapping theorem to multivariate functions of matrices, incorporating both Hermitian and non-Hermitian matrices. With functions assumed to be holomorphic, this generalized theorem expands the applicability of spectral theory to a wider class of operators. The framework established for finite-dimensional matrices is further extended to infinite-dimensional settings, including countable spectrum operators, enhancing understanding of operator behavior in broader contexts. For continuous spectrum operators, the paper advances von Neumann's spectral theorem for self-adjoint operators to include a larger class of spectral operators, accommodating both self-adjoint and non-self-adjoint cases. This expansion provides a unified spectral framework that generalizes spectral decomposition, allowing the spectral mapping theorem to be applied in diverse cases. Furthermore, analogous to finite dimensions, the concept of analogous operators is introduced for continuous spectrum operators, enhancing operator classification. For operators with hybrid spectrum—comprising both discrete and continuous elements—analogous properties and spectral mapping are also explored. 
\end{abstract}


\section{Introduction}\label{sec: Introduction}

Spectral analysis of operators is crucial in mathematics because it provides deep insights into the structure and behavior of linear operators, particularly through the study of their spectra (the set of eigenvalues and corresponding eigenspaces)~\cite{fucik2006spectral,aiello2022spectral}. This analysis generalizes concepts like diagonalization and eigenvalues from finite-dimensional matrices to infinite-dimensional spaces, enabling us to solve various functional and differential equations. In quantum mechanics, for instance, operators represent physical observables, and their spectra correspond to possible measurement outcomes. Spectral theory helps predict system behavior, solve linear PDEs, and provides foundational tools in areas like harmonic analysis and probability theory~\cite{chang2020spectral,gudder1965spectral,karlovich2018pseudodifferential}.

In science and technology, spectral analysis plays a key role in diverse fields such as signal processing, image compression, and machine learning. In signal processing, for example, the Fourier transform is an application of spectral theory, decomposing signals into their frequency components. In engineering, eigenvalue problems appear in vibration analysis and stability studies, while in machine learning, techniques like Principal Component Analysis (PCA) rely on spectral methods to reduce data dimensions, identify patterns, and improve computational efficiency~\cite{guo2009principal}. Spectral methods are also essential in network analysis, where the spectrum of the graph Laplacian reveals information about the structure and dynamics of complex networks~\cite{chung1996lectures}.

In an excellent physics review~\cite{ashida2020non}, the authors provides a comprehensive exploration of the foundations and applications of \emph{non-Hermitian physics} in both classical and quantum systems. It begins by introducing key concepts from non-Hermitian linear algebra, such as Jordan normal form, biorthogonality, exceptional points, pseudo-Hermiticity, and parity-time (PT) symmetry. These mathematical tools are fundamental for understanding phenomena unique to non-Hermitian systems. The review then discusses various classical systems—such as photonics, mechanics, and acoustics—that simulate non-Hermitian wave physics, highlighting phenomena like unidirectional invisibility, enhanced sensitivity, and coherent perfect absorption. In the quantum realm, the review explains how non-Hermitian operators describe open quantum systems using approaches like the Feshbach projection and quantum trajectory. These frameworks are applied to a wide range of physical systems, from atomic and molecular physics to nuclear physics, where phenomena such as quantum resonances, superradiance, and the quantum Zeno effect are also discussed. The review also covers the emerging field of band topology in non-Hermitian systems, offering a complete classification and examples. Additionally, topics like nonreciprocal transport, speed limits, and nonunitary quantum walks are explored, demonstrating the broad relevance of non-Hermitian physics across multiple scientific domains. All these fascinating studies in physics related to non-Hermitian systems are grounded in either non-Hermitian matrices (finite dimensions) or non-self-adjoint operators (infinite dimensions).

This paper introduces several important contributions to the study of countable spectrum (include both finite-dimensional and infinite-dimesional operators), continous spectrum, and hybrid spectrum operators, with a focus on developing new frameworks for analyzing matrices and operators, especially in the context of projectors and nilpotents. For finite-dimensional matrices, the paper proposes a novel concept of analogous matrices, where two matrices are considered analogous if they share the same projectors and nilpotent structures. This concept provides a new lens through which to view matrix equivalence, emphasizing their internal structures rather than just their spectra. Building on this, the paper presents a graph representation of these shared projectors and nilpotent elements, enabling a visual and computational way to classify matrices within the same family of analogous matrices. Another major contribution is the enumeration of the total number of different families of analogous matrices, given the matrix size, which provides a foundational result for understanding matrix classification.

Furthermore, the paper extends the spectral mapping theorem to functions with multivariate matrix inputs, allowing for both Hermitian and non-Hermitian matrices as arguments. The functions are assumed to be holomorphic and represented by infinite series, extending classical results into a more general, multivariate setting. This result is particularly valuable for applications where matrix inputs are not limited to Hermitian matrices, broadening the applicability of spectral theory in physics and other domains. Analogous properties and spctral mapping theorem for finite-dimensional operators are also extended to countable infinite-dimensional operators.  

For continous spectrum operators, the paper first extends the von Neumann spectral theorem for self-adjoint operators by deriving the \emph{spectral mapping theorem} for a broader class of spectral operators, which includes both self-adjoint and non-self-adjoint operators. This extension generalizes the spectral decomposition concept, allowing the spectral mapping theorem to apply to operators beyond the self-adjoint case, accommodating a wider variety of spectral types while still preserving many of the core spectral properties established in the self-adjoint framework. The result is a unified spectral framework that encompasses the von Neumann spectral theorem as a specific case within this broader setting, thereby expanding the applicability of spectral theory to include non-self-adjoint operators as well. Additionally, analogous to the finite-dimensional case, the paper introduces the concept of analogous operators for continous spectral operators and studies their properties, expanding the understanding of operator equivalence beyond finite-dimensional settings. Lastly, the paper presents the multivariate spectral mapping theorem for continous spectral operators, where the inputs are several continous spectrum operators of holomorphic functions, and the results hold for both self-adjoint and non-self-adjoint operators. Finally, we consider hybrid spectrum operators, where spectrum is composed by discrete set and continous set. Analogous properties and spctral mapping theorem for hybrid operators operators are also explored.

A key remark made in the paper is that while the resolvent approach is the traditional method for the spectral mapping theorem (via functional calculus), our proposed method based on projectors and nilpotents (Jordan Decomposition) offers significant advantages~\cite{taylor1970joint,colombo2008new,aiello2022spectral}. The resolvent approach falls short in capturing the nilpotent structure of operators, especially in identifying generalized eigenspaces and Jordan blocks, which are essential for understanding nilpotent behavior. This limitation becomes especially pronounced when the eigenvalue is zero, making it challenging to detect or analyze nilpotency. In contrast, our projector and nilpotent method directly addresses these gaps, providing a more complete framework for spectral analysis of countable spectrum (include both finite-dimensional and infinite-dimesional operators), continous spectrum, and hybrid spectrum operators. Furthermore, our approach effectively accommodates the spectral analysis of noncommutative operators, expanding the scope of tools available for complex operator behavior.

In Section~\ref{sec: Finite-Dimensional Operators}, we begin by reviewing the Jordan decomposition for square matrices, followed by a discussion of analogous operator properties, an exploration of finite-dimensional operator type enumeration, and the establishment of the spectral mapping theorem for multivariate square matrices. Section~\ref{sec: Countable Infinite-Dimensional Operators} extends the concepts of analogous operators and the spectral mapping theorem to countable infinite-dimensional operators. In Section~\ref{sec: Continous Spectrum Operators}, we begin with the study of analogous continuous spectrum operators and establish the multivariate spectral mapping theorem for continuous spectrum operators. Finally, Section~\ref{sec: Hybrid Spectrum Operators} presents analogous properties and the spectral mapping theorem for hybrid spectrum operators.

\section{Finite-Dimensional Operators}\label{sec: Finite-Dimensional Operators}

The section is organized as follows:
\begin{itemize}
\item Section~\ref{sec: Jordan Decomposition Review} reviews the Jordan decomposition for square matrices.
\item Section~\ref{sec: Analogous Matrices and Their Properties} introduces the concept of ``Analogous Operators" between square matrices and discusses their properties.
\item Section~\ref{sec: Enumeration of Finite-Dimensional Operator Types} explores the enumeration of analogous operator families for matrices with given dimension.
\item Section~\ref{sec: Spectral Mapping Theorem for Multivariate Finite-Dimensional Operators}establishes a spectral mapping theorem for multivariate square matrices.
\end{itemize}

\subsection{Jordan Decomposition Review}\label{sec: Jordan Decomposition Review}

The Jordan decomposition, also known as the \emph{Jordan canonical form} or \emph{Jordan normal form}, is a mathematical concept used in linear algebra to simplify the representation of matrices. It expresses a matrix in a nearly diagonal form, making it easier to analyze and understand the matrix’s properties, particularly for matrices that are not diagonalizable.

We review the following important concepts. A Jordan Block with dimension $m \times m$ and the eigenvalue $\lambda$, denoted by $J_m(\lambda)$, can be expressed by  
\begin{eqnarray}\label{eq: Jordan block def}
\bm{J}_m(\lambda)= \begin{pmatrix}
   \lambda & 1 & 0 & \cdots & 0 \\
   0 & \lambda & 1 & \cdots & 0 \\
   \vdots & \vdots & \ddots & \ddots & \vdots \\
   0 & 0 & \cdots & \lambda & 1 \\
   0 & 0 & \cdots & 0 & \lambda
   \end{pmatrix}_{m \times m}
\end{eqnarray}
where $\lambda$ is a scalar (an eigenvalue of the original matrix), and the off-diagonal elements are either 0 or 1. 

The algebraic multiplicity of an eigenvalue $\lambda_k$ of a matrix $\bm{X}$ refers to how many times $\lambda_k$ appears as a root of the \emph{characteristic polynomial} of $\bm{X}$. If $\lambda_k$ is an eigenvalue, the algebraic multiplicity, denoted by $\alpha_k^{(\mathrm{A})}$, is the power (or count) of $\lambda_k$ in the factorization of the characteristic polynomial, i.e., $\mbox{det}(\bm{X} - \lambda_k \bm{I})$$=$$0$, where $\bm{I}$ is the identity matrix of the same size as $\bm{X}$. On the other hand, the geometric multiplicity of an eigenvalue $\lambda_k$ of a matrix $\bm{X}$, denoted by $\alpha_k^{(\mathrm{G})}$, is the dimension of the eigenspace associated with $\lambda_k$. If the eigenspace of $\lambda_k$ is defined as: $\mbox{Null}(\bm{X} - \lambda_k \bm{I})$, where $\mbox{Null}(\bm{X} - \lambda_k \bm{I})$ refers to the null space (or kernel) of the matrix  $(\bm{X} - \lambda_k \bm{I})$. The dimension of this null space gives the geometric multiplicity of $\lambda_k$.

Jordan decomposition theorem says that we can decompose a square matrix $\bm{X} \in \mathbb{C}^{m \times m}$ as~\cite{gohberg1996simple}:  
\begin{eqnarray}\label{eq: Jordan decompostion}
\bm{X}&=& \bm{U}\left(\bigoplus\limits_{k=1}^{K}\bigoplus\limits_{i=1}^{\alpha_k^{(\mathrm{G})}}\bm{J}_{m_{k,i}}(\lambda_k)\right)\bm{U}^{-1}.
\end{eqnarray}
where $\bm{U} \in \mathbb{C}^{m \times m}$ is an invertible matrix, and $\alpha_k^{(\mathrm{G})}$ is the geometry multiplicity with respect to the $k$-th eigenvalue $\lambda_k$. We have the following relationships about $\alpha_k^{(\mathrm{A})}$ and $\alpha_k^{(\mathrm{G})}$:
\begin{eqnarray}
\sum\limits_{k=1}^K \alpha_k^{(\mathrm{A})}&=&m,
\end{eqnarray}
and 
\begin{eqnarray}
\sum\limits_{i=1}^{\alpha_k^{(\mathrm{G})}}m_{k,i}&=&\alpha_k^{(\mathrm{A})}.
\end{eqnarray}

We have the following lemma about Jordan decomposition.
\begin{lemma}\label{lma: projectors and nilpotents of a matrix}
Given a matrix $\bm{X} \in \mathbb{C}^{m \times m}$, we can express it by
\begin{eqnarray}\label{eq1: lma: projectors and nilpotents of a matrix}
\bm{X}&=&\sum\limits_{k=1}^K\sum\limits_{i=1}^{\alpha_k^{\mathrm{G}}} \lambda_k \bm{P}_{k,i}+
\sum\limits_{k=1}^K\sum\limits_{i=1}^{\alpha_k^{\mathrm{G}}} \bm{N}_{k,i},
\end{eqnarray}
where $\bm{P}_{k,i}\in \mathbb{C}^{m \times m}$ is the projector matrices with respect to the eigenvalue $\lambda_k$ and $i$-th dimensoin of the null space $\mbox{Null}(\bm{X} - \lambda_k \bm{I})$; and $\bm{N}_{k,i}\in \mathbb{C}^{m \times m}$ is the nilpotent matrices with respect to the eigenvalue $\lambda_k$ and $i$-th dimensoin of the null space $\mbox{Null}(\bm{X} - \lambda_k \bm{I})$.
\end{lemma}
\textbf{Proof:}
Since each Jordan Block with dimension $m_{k,i}\times m_{k,i}$ and the eigenvalue $\lambda_k$, denoted by $\bm{J}_{m_{k,i}}(\lambda_k)$, can be expressed by  
\begin{eqnarray}\label{eq1.1: lma: projectors and nilpotents of a matrix}
\bm{J}_{m_{k,i}}(\lambda_k)&=&\begin{pmatrix}
   \lambda_k & 1 & 0 & \cdots & 0 \\
   0 & \lambda_k & 1 & \cdots & 0 \\
   \vdots & \vdots & \ddots & \ddots & \vdots \\
   0 & 0 & \cdots & \lambda_k & 1 \\
   0 & 0 & \cdots & 0 & \lambda_k
   \end{pmatrix}_{m_{k,i} \times m_{k,i}} \nonumber \\
&=& \lambda_k \bm{I}_{m_{k,i}}+ \acute{\bm{J}}_{m_{k,i}},
\end{eqnarray}
where $\bm{I}_{m_{k,i}}$ is an identity matrix with the dimension $m_{k,i} \times m_{k,i}$, and $\acute{\bm{J}}_{m_{k,i}}$ is a off-diagonal matrix with the dimension $m_{k,i} \times m_{k,i}$ such that 
\begin{eqnarray}\label{eq2: lma: projectors and nilpotents of a matrix}
\acute{\bm{J}}_{m_{k,i}}&=&\begin{pmatrix}
   0 & 1 & 0 & \cdots & 0 \\
   0 & 0 & 1 & \cdots & 0 \\
   \vdots & \vdots & \ddots & \ddots & \vdots \\
   0 & 0 & \cdots & 0 & 1 \\
   0 & 0 & \cdots & 0 & 0
   \end{pmatrix}_{m_{k,i} \times m_{k,i}}. 
\end{eqnarray}

We construct the identity matrix $\bm{I}_{k,i}$ with the dimension $m \times m$ as
\begin{eqnarray}\label{eq3: lma: projectors and nilpotents of a matrix}
\bm{I}_{k,i}&=&\bigoplus\limits_{k'=1}^{K}\bigoplus\limits_{i'=1}^{\alpha_{k'}^{(\mathrm{G})}}\delta(k,k')\delta(i',i)\bm{I}_{m_{k',i'}},
\end{eqnarray}
where $\delta(a,b)$$=$$0$ if $a$$\neq$$b$, and $\delta(a,b)$$=$$1$ if $a$$=$$b$. Similarly, we also construct the off-diagonal matrix $\acute{\bm{J}}$ with the dimension $m \times m$ as
\begin{eqnarray}\label{eq4: lma: projectors and nilpotents of a matrix}
\acute{\bm{J}}_{k,i}&=&\bigoplus\limits_{k'=1}^{K}\bigoplus\limits_{i'=1}^{\alpha_{k'}^{(\mathrm{G})}}\delta(k,k')\delta(i',i)\acute{\bm{J}}_{m_{k',i'}}.
\end{eqnarray}
Let $\bm{P}_{k,i}$ be constructed by 
\begin{eqnarray}\label{eq5: lma: projectors and nilpotents of a matrix}
\bm{P}_{k,i}&=&\bm{V}\bm{I}_{k,i}\bm{V}^{-1},
\end{eqnarray}
and 
\begin{eqnarray}\label{eq6: lma: projectors and nilpotents of a matrix}
\bm{N}_{k,i}&=&\bm{V}\acute{\bm{J}}_{k,i}\bm{V}^{-1}.
\end{eqnarray}
From Eq.~\eqref{eq: Jordan decompostion}, we have Eq.~\eqref{eq1: lma: projectors and nilpotents of a matrix}.

Finally, we have to show that matrices $\bm{P}_{k,i}$ are projector matrices and matrices $\bm{N}_{k,i}$ are nilpotent matrices. Because we have
\begin{eqnarray}\label{eq7: lma: projectors and nilpotents of a matrix}
\bm{P}_{k,i}\bm{P}_{k',i'}=\bm{P}_{k,i}\delta(k,k')\delta(i,i'),
\end{eqnarray}
and
\begin{eqnarray}\label{eq8: lma: projectors and nilpotents of a matrix}
\sum\limits_{k=1}^K\sum\limits_{i=1}^{\alpha_k^{\mathrm{G}}}\bm{P}_{k,i}=\bm{I}.
\end{eqnarray}
Therefore,  matrices $\bm{P}_{k,i}$ are projectors as they are orthogonal and complete. On the other hand, 
\begin{eqnarray}\label{eq9: lma: projectors and nilpotents of a matrix}
\bm{N}_{k,i}^{\ell}\neq\bm{0}, \mbox{~and~~} \bm{N}_{k,i}^{m_{k,i}} =\bm{0}, 
\end{eqnarray}
where $\ell \in \mathbb{N}$ and $\ell < m_{k,i}$. This shows that matrices $\bm{N}_{k,i}$ are nilpotent matrices.
$\hfill\Box$

\subsection{Analogous Matrices and Their Properties}\label{sec: Analogous Matrices and Their Properties}

We will introduce a new concept, named as \emph{Analogous Operators}, between two square matrices and its properties in Section~\ref{sec: Analogous Operators}. In Section~\ref{sec: Graph Representation of Projectors and Nilpotents in Finite-Dimensional Operators}, a graph representation method is proposed to demonstrate a family of matrices which are analogous.

\subsubsection{Analogous Operators}\label{sec: Analogous Operators}

In linear algebra, when two matrices share the same Jordan form but differ in eigenvalues, it reveals an intriguing mathematical relationship. Both matrices possess the same \emph{generalized eigenvector structure}—with identical Jordan block arrangements, dimensions, and chain lengths. This means that the underlying geometric relationships in both matrices remain unchanged, even though their eigenvalues differ. The geometric multiplicity and algebraic multiplicity remain identical, preserving the fundamental structure of the transformation. However, the matrices differ in how the eigenvector chains are scaled, as eigenvalues represent the stretching or compression along corresponding eigenspaces.

This insight motivates us to explore new notions of \emph{analogous matrices} that go beyond simple similarity. If two matrices can have the same underlying structure yet exhibit different behaviors in terms of scaling, we might design new frameworks to classify matrices based on their geometric and algebraic traits, but with flexibility in eigenvalue-based transformations. Such a notion could help us understand systems that share core dynamics but differ in specific physical or theoretical properties.

Given two matrices $\bm{X}$ and $\bm{Y}$ with the same size $m \times m$ and the same number of distinct eigenvalues $\lambda_k$ for $k=1,2,\ldots,K$, we say that the matrix $\bm{X}$ is \emph{analogous} to the matrix $\bm{Y}$ with respect to the ratios $[c_1, c_2, \ldots, c_K] \in \mathbb{C}^K$ and $c_i \neq 0$ for $k=1,2,\ldots,K$, denoted by $\bm{X} \propto^{\bm{U},\bm{V}}_{[c_1, c_2, \ldots, c_K]} \bm{Y}$, if these two matrices $\bm{X}$ and $\bm{Y}$ can be expressed by
\begin{eqnarray}\label{eq1: Analogous Operators}
\bm{X}&=& \bm{U}\left(\bigoplus\limits_{k=1}^{K}\bigoplus\limits_{i=1}^{\alpha_k^{(\mathrm{G})}}\bm{J}_{m_{k,i}}(\lambda_k)\right)\bm{U}^{-1}, \nonumber \\
\bm{Y}&=& \bm{V}\left(\bigoplus\limits_{k=1}^{K}\bigoplus\limits_{i=1}^{\alpha_k^{(\mathrm{G})}}\bm{J}_{m_{k,i}}(c_k \lambda_k)\right)\bm{V}^{-1}.
\end{eqnarray}

We will have the following properties about two analogous matrices $\bm{X}$ and $\bm{Y}$ given by the relationship $\bm{X} \propto^{\bm{U},\bm{V}}_{[c_1, c_2, \ldots, c_K]} \bm{Y}$. 

\begin{proposition}\label{prop0: Analogous Operators}
Two analogous matrices $\bm{X}$ and $\bm{Y}$ given by the relationship $\bm{X} \propto^{\bm{U},\bm{V}}_{[c_1, c_2, \ldots, c_K]} \bm{Y}$, then, these two matrices $\bm{X}$ and $\bm{Y}$ have the same rank.
\end{proposition}
\textbf{Proof:}
From Eq.~\eqref{eq1: Analogous Operators}, we know that the number of zero-valued eigenvalues of the matrix $\bm{X}$ should be equal to the number of zero-valued eigenvalues of the matrix $\bm{Y}$.
$\hfill\Box$

The next proposition is to show that the similar relationship between two matrices is the special case of the proposed analogous relationship. 

\begin{proposition}\label{prop1: Analogous Operators}
Two analogous matrices $\bm{X}$ and $\bm{Y}$ given by the relationship $\bm{X} \propto^{\bm{U},\bm{V}}_{[1, 1, \ldots, 1]} \bm{Y}$, then, these two matrices $\bm{X}$ and $\bm{Y}$ are similar.
\end{proposition}
\textbf{Proof:}
Since we have 
\begin{eqnarray}\label{eq1: prop1: Analogous Operators}
\bm{X}&=& \bm{U}\left(\bigoplus\limits_{k=1}^{K}\bigoplus\limits_{i=1}^{\alpha_k^{(\mathrm{G})}}\bm{J}_{m_{k,i}}(\lambda_k)\right)\bm{U}^{-1}, \nonumber \\
\bm{Y}&=& \bm{V}\left(\bigoplus\limits_{k=1}^{K}\bigoplus\limits_{i=1}^{\alpha_k^{(\mathrm{G})}}\bm{J}_{m_{k,i}}(1 \times \lambda_k)\right)\bm{V}^{-1},
\end{eqnarray}   
therefore, we can obtain the matrix $\bm{Y}$ from the matrix $\bm{X}$ by 
\begin{eqnarray}\label{eq2: prop1: Analogous Operators}
\bm{V}\bm{U}^{-1}\bm{X}\bm{U}\bm{V}^{-1}&=& \bm{Y}.
\end{eqnarray}   
$\hfill\Box$

The next proposition is about the commutative condition of two analogous matrices.

\begin{proposition}\label{prop2: Analogous Operators}
Two analogous matrices $\bm{X}$ and $\bm{Y}$ are given by the relationship $\bm{X} \propto^{\bm{U},\bm{U}}_{[c_1, c_2, \ldots, c_K]} \bm{Y}$, then, these two matrices $\bm{X}$ and $\bm{Y}$ are commute under the conventional matrix multiplication. 
\end{proposition}
\textbf{Proof:}
Becuase two analogous matrices $\bm{X}$ and $\bm{Y}$ are given by the relationship $\bm{X} \propto^{\bm{U},\bm{U}}_{[c_1, c_2, \ldots, c_K]} \bm{Y}$, we have
\begin{eqnarray}\label{eq1: prop2: Analogous Operators}
\bm{X}&=& \bm{U}\left(\bigoplus\limits_{k=1}^{K}\bigoplus\limits_{i=1}^{\alpha_k^{(\mathrm{G})}}\bm{J}_{m_{k,i}}(\lambda_k)\right)\bm{U}^{-1}\nonumber \\
&=&\sum\limits_{k=1}^K\sum\limits_{i=1}^{\alpha_k^{\mathrm{G}}} \lambda_k \bm{P}_{k,i}+
\sum\limits_{k=1}^K\sum\limits_{i=1}^{\alpha_k^{\mathrm{G}}} \bm{N}_{k,i}, 
\end{eqnarray}   
and
\begin{eqnarray}\label{eq2: prop2: Analogous Operators}
\bm{Y}&=& \bm{U}\left(\bigoplus\limits_{k=1}^{K}\bigoplus\limits_{i=1}^{\alpha_k^{(\mathrm{G})}}\bm{J}_{m_{k,i}}(c_k \times \lambda_k)\right)\bm{U}^{-1}\nonumber \\
&=&\sum\limits_{k=1}^K\sum\limits_{i=1}^{\alpha_k^{\mathrm{G}}}c_k \lambda_k \bm{P}_{k,i}+
\sum\limits_{k=1}^K\sum\limits_{i=1}^{\alpha_k^{\mathrm{G}}} \bm{N}_{k,i}.
\end{eqnarray}   

Besides, we also have the following relations among $\bm{P}_{k,i}$ and $\bm{N}_{k,i}$ from Eq.~\eqref{eq5: lma: projectors and nilpotents of a matrix} and Eq.~\eqref{eq6: lma: projectors and nilpotents of a matrix}:
\begin{eqnarray}\label{eq3: prop2: Analogous Operators}
\bm{P}_{k,i}\bm{P}_{k',i'}&=&\bm{P}_{k,i}\delta(k,k')\delta(i,i'), \nonumber \\
\bm{P}_{k',i'}\bm{N}_{k,i}&=&\bm{N}_{k,i}\bm{P}_{k',i'}=\bm{N}_{k,i}\delta(k,k')\delta(i,i'), \nonumber \\
\bm{N}_{k,i}\bm{N}_{k',i'}&=&\bm{N}^2_{k,i}\delta(k,k')\delta(i,i').
\end{eqnarray}   
Then, we have 
\begin{eqnarray}\label{eq4: prop2: Analogous Operators}
\bm{X}\bm{Y}&=& \sum\limits_{k=1}^K\sum\limits_{i=1}^{\alpha_k^{\mathrm{G}}} c_k \lambda^2_k \bm{P}_{k,i}+ \sum\limits_{k=1}^K\sum\limits_{i=1}^{\alpha_k^{\mathrm{G}}} c_k \lambda_k \bm{N}_{k,i} \nonumber \\
&&+ \sum\limits_{k=1}^K\sum\limits_{i=1}^{\alpha_k^{\mathrm{G}}} \lambda_k \bm{N}_{k,i}+
\sum\limits_{k=1}^K\sum\limits_{i=1}^{\alpha_k^{\mathrm{G}}} \bm{N}^2_{k,i} \nonumber \\
&=&  \bm{Y}\bm{X}
\end{eqnarray}   
$\hfill\Box$

The next proposition is to show that the proposed analogous relationship between $\bm{X}$ and $\bm{Y}$ provided by $\bm{X} \propto^{\bm{U},\bm{V}}_{[c_1, c_2, \ldots, c_K]} \bm{Y}$ is an equivalence relation.

\begin{proposition}\label{prop3: Analogous Operators}
The relationship $\bm{X} \propto^{\bm{U},\bm{V}}_{[c_1, c_2, \ldots, c_K]} \bm{Y}$ is an equivalence relation.
\end{proposition}
\textbf{Proof:}
Given the matrices $\bm{X}, \bm{Y}$ and $\bm{Z}$ with the following formats:
\begin{eqnarray}\label{eq1: prop3: Analogous Operators}
\bm{X}&=& \bm{U}\left(\bigoplus\limits_{k=1}^{K}\bigoplus\limits_{i=1}^{\alpha_k^{(\mathrm{G})}}\bm{J}_{m_{k,i}}(\lambda_k)\right)\bm{U}^{-1}, \nonumber \\
\bm{Y}&=& \bm{V}\left(\bigoplus\limits_{k=1}^{K}\bigoplus\limits_{i=1}^{\alpha_k^{(\mathrm{G})}}\bm{J}_{m_{k,i}}(c_k \lambda_k)\right)\bm{V}^{-1},\nonumber \\
\bm{Z}&=& \bm{W}\left(\bigoplus\limits_{k=1}^{K}\bigoplus\limits_{i=1}^{\alpha_k^{(\mathrm{G})}}\bm{J}_{m_{k,i}}(d_k \lambda_k)\right)\bm{W}^{-1},
\end{eqnarray} 
we have 
\begin{eqnarray}\label{eq2s: prop3: Analogous Operators}
\bm{X}&\propto^{\bm{U},\bm{U}}_{[1, 1, \ldots, 1]}&\bm{X}, ~(\mbox{reflexive})\nonumber \\
\bm{X}&\propto^{\bm{U},\bm{V}}_{[c_1, c_2, \ldots, c_K]}&\bm{Y}, \mbox{~and $\bm{Y}\propto^{\bm{V},\bm{U}}_{[1/c_1, 1/c_2, \ldots, 1/c_K]}\bm{X}$}~(\mbox{symmetric})\nonumber \\
\bm{X}&\propto^{\bm{U},\bm{V}}_{[c_1, c_2, \ldots, c_K]}&\bm{Y}, \mbox{~and $\bm{Y}\propto^{\bm{V},\bm{W}}_{[d_1/c_1, d_2/c_2, \ldots, d_K/c_K]}\bm{Z}$}\nonumber \\
&\Longrightarrow& \bm{X}\propto^{\bm{U},\bm{W}}_{[d_1, d_2, \ldots, d_K]}\bm{Z}~(\mbox{transitive})
\end{eqnarray} 
$\hfill\Box$

The next proposition is to show the matrix determinant relationship between two analogous matrices. 
\begin{proposition}\label{prop4: Analogous Operators}
Given two analogous matrices $\bm{X}$ and $\bm{Y}$ provided by the relationship $\bm{X} \propto^{\bm{U},\bm{V}}_{[c_1, c_2, \ldots, c_K]} \bm{Y}$ with
\begin{eqnarray}\label{eq1: prop4: Analogous Operators}
\bm{X}&=& \bm{U}\left(\bigoplus\limits_{k=1}^{K}\bigoplus\limits_{i=1}^{\alpha_k^{(\mathrm{G})}}\bm{J}_{m_{k,i}}(\lambda_k)\right)\bm{U}^{-1}, \nonumber \\
\bm{Y}&=& \bm{V}\left(\bigoplus\limits_{k=1}^{K}\bigoplus\limits_{i=1}^{\alpha_k^{(\mathrm{G})}}\bm{J}_{m_{k,i}}(c_k \lambda_k)\right)\bm{V}^{-1};
\end{eqnarray}   
we have 
\begin{eqnarray}\label{eq2: prop4: Analogous Operators}
\det(\bm{X})&=&\frac{\det(\bm{Y})}{\left(\prod\limits_{k=1}^K c_k^{\alpha_k^{(\mathrm{A})}}\right)}. 
\end{eqnarray}
\end{proposition}
\textbf{Proof:}
Because we have
\begin{eqnarray}\label{eq3: prop4: Analogous Operators}
\det(\bm{X})&=&\prod\limits_{k=1}^K \lambda_k^{\alpha_k^{(\mathrm{A})}},
\end{eqnarray}
and 
\begin{eqnarray}\label{eq4: prop4: Analogous Operators}
\det(\bm{Y})&=&\prod\limits_{k=1}^K\left(c_k\lambda_k\right)^{\alpha_k^{(\mathrm{A})}},
\end{eqnarray}
this proposition is proved by comparing Eq.~\eqref{eq3: prop4: Analogous Operators} and Eq.~\eqref{eq4: prop4: Analogous Operators}.
$\hfill\Box$

The next proposition is to show the matrix trace relationship for two analogous matrices.

\begin{proposition}\label{prop5: Analogous Operators}
Given two analogous matrices $\bm{X}$ and $\bm{Y}$ provided by the relationship $\bm{X} \propto^{\bm{U},\bm{V}}_{[c, c, \ldots, c]} \bm{Y}$ with
\begin{eqnarray}\label{eq1: prop5: Analogous Operators}
\bm{X}&=& \bm{U}\left(\bigoplus\limits_{k=1}^{K}\bigoplus\limits_{i=1}^{\alpha_k^{(\mathrm{G})}}\bm{J}_{m_{k,i}}(\lambda_k)\right)\bm{U}^{-1}, \nonumber \\
\bm{Y}&=& \bm{V}\left(\bigoplus\limits_{k=1}^{K}\bigoplus\limits_{i=1}^{\alpha_k^{(\mathrm{G})}}\bm{J}_{m_{k,i}}(c \lambda_k)\right)\bm{V}^{-1};
\end{eqnarray}   
we have 
\begin{eqnarray}\label{eq2: prop5: Analogous Operators}
\mathrm{Trace}(\bm{X})&=&c \mathrm{Trace}(\bm{Y})
\end{eqnarray}
\end{proposition}
\textbf{Proof:}
Because we have
\begin{eqnarray}\label{eq3: prop5: Analogous Operators}
\mathrm{Trace}(\bm{X})&=&\sum\limits_{k=1}^K \alpha_k^{(\mathrm{A})}\lambda_k,
\end{eqnarray}
and 
\begin{eqnarray}\label{eq4: prop5: Analogous Operators}
\mathrm{Trace}(\bm{Y})&=&\sum\limits_{k=1}^K \alpha_k^{(\mathrm{A})}c\lambda_k,
\end{eqnarray}
this proposition is proved by comparing Eq.~\eqref{eq3: prop5: Analogous Operators} and Eq.~\eqref{eq4: prop5: Analogous Operators}.
$\hfill\Box$

The next proposition is to show the relationship of characteristic polynomials between two analogous matrices. 
\begin{proposition}\label{prop6: Analogous Operators}
Given two analogous matrices $\bm{X}$ and $\bm{Y}$ provided by the relationship $\bm{X} \propto^{\bm{U},\bm{V}}_{[c, c, \ldots, c]} \bm{Y}$ with
\begin{eqnarray}\label{eq1: prop6: Analogous Operators}
\bm{X}&=& \bm{U}\left(\bigoplus\limits_{k=1}^{K}\bigoplus\limits_{i=1}^{\alpha_k^{(\mathrm{G})}}\bm{J}_{m_{k,i}}(\lambda_k)\right)\bm{U}^{-1}, \nonumber \\
\bm{Y}&=& \bm{V}\left(\bigoplus\limits_{k=1}^{K}\bigoplus\limits_{i=1}^{\alpha_k^{(\mathrm{G})}}\bm{J}_{m_{k,i}}(c \lambda_k)\right)\bm{V}^{-1};
\end{eqnarray}   
we have 
\begin{eqnarray}\label{eq2: prop6: Analogous Operators}
\mathrm{CP}_{\bm{X}}(x)&=&\mathrm{CP}_{\bm{Y}}\left(\frac{x}{c}\right),
\end{eqnarray}
where $\mathrm{CP}_{\bm{X}}$ is the  characteristic polynomial of the matrix $\bm{X}$ and $\mathrm{CP}_{\bm{Y}}$ is the  characteristic polynomial of the matrix $\bm{Y}$.
\end{proposition}
\textbf{Proof:}
Since two polynomials $\mathrm{CP}_{\bm{X}}$ and $\mathrm{CP}_{\bm{Y}}$ have roots that share the same scalar ratio $c$, it suggests a proportional relationship $c$ between the roots of the polynomials. 
$\hfill\Box$

\subsubsection{Graph Representation of Projectors and Nilpotents in Finite-Dimensional Operators}\label{sec: Graph Representation of Projectors and Nilpotents in Finite-Dimensional Operators}

To represent the algebraic structure of projectors and nilpotents for matrices come from the same analogous family, we can formulate a graph where the nodes correspond to the elements $\bm{P}_{k,i}$, $\bm{N}_{k,i}$, and $\bm{0}$, while the edges capture the interactions among $\bm{P}_{k,i}$ and $\bm{N}_{k,i}$. The relations described in Eq.~\eqref{eq3: prop2: Analogous Operators} provide the foundation for defining both the nodes and the edges. Such graph, named as \emph{Analogous Structure Graph} (ASG), will be used to represent the whole family of analogous matrices since they share the same elements of $\bm{P}_{k,i}$ and $\bm{N}_{k,i}$.

We begin with graph node representation. Each element $\bm{P}_{k,i}$, $\bm{N}_{k,i}$, and $\bm{0}$ corresponds to a node in the graph:
\begin{itemize}
\item $\bm{P}_{k,i}$: Represents a projection operator that behaves according to the Kronecker delta conditions. It will have a self-loop as $\bm{P}_{k,i}^2 $$=$$ \bm{P}_{k,i}$.
\item $\bm{N}_{k,i}$: Represents a nilpotent matrix operator with specific degrees of nilpotency, i.e., $\bm{N}^\ell_{k,i} \neq \bm{0}$ for $\ell < m_{k,i}$ and $\bm{N}^\ell_{k,i} $$=$$ \bm{0}$ for $\ell \geq m_{k,i}$. Thus, $\bm{N}_{k,i}$ is linked to $\bm{0}$ as it becomes the zero operator after reaching the maximum nilpotency degree.
\item $\bm{0}$: The zero matrix, which can be interpreted as an absorbing node for any operator becoming zero.
\end{itemize}

Edges between the nodes will reflect the multiplication rules described in Eq.~\eqref{eq3: prop2: Analogous Operators}:
\begin{itemize}
\item $\bm{P}_{k,i}\bm{P}_{k',i'}$$=$$ \bm{P}_{k,i}\delta(k,k')\delta(i,i')$: This relation indicates that $\bm{P}_{k,i}$ only interacts with itself, forming self-loops on nodes corresponding to $\bm{P}_{k,i}$, with no direct edges between different $\bm{P}_{k,i}$ and $\bm{P}_{k',i'}$ unless $k = k'$ and $i = i'$.
\item $\bm{P}_{k',i'}\bm{N}_{k,i} $$=$$ \bm{N}_{k,i}\bm{P}_{k',i'} $$=$$ \bm{N}_{k,i}\delta(k,k')\delta(i,i')$: These interactions define edges between $\bm{N}_{k,i}$ and $\bm{P}_{k,i}$, forming directional edges where $\bm{N}_{k,i}$ maps to $\bm{P}_{k,i}$, provided the Kronecker delta condition is satisfied.
\item $\bm{N}_{k,i}\bm{N}_{k',i'} $$=$$ \bm{N}^2_{k,i}\delta(k,k')\delta(i,i')$: This equation illustrates that $\bm{N}_{k,i}$ interacts with itself with the nilpotency condition indicating that higher powers of $\bm{N}_{k,i}$ eventually map to $\bm{0}$.
\end{itemize}

\begin{figure}[htbp!]
	\centerline{\includegraphics[width=0.9\columnwidth,draft=false]
		{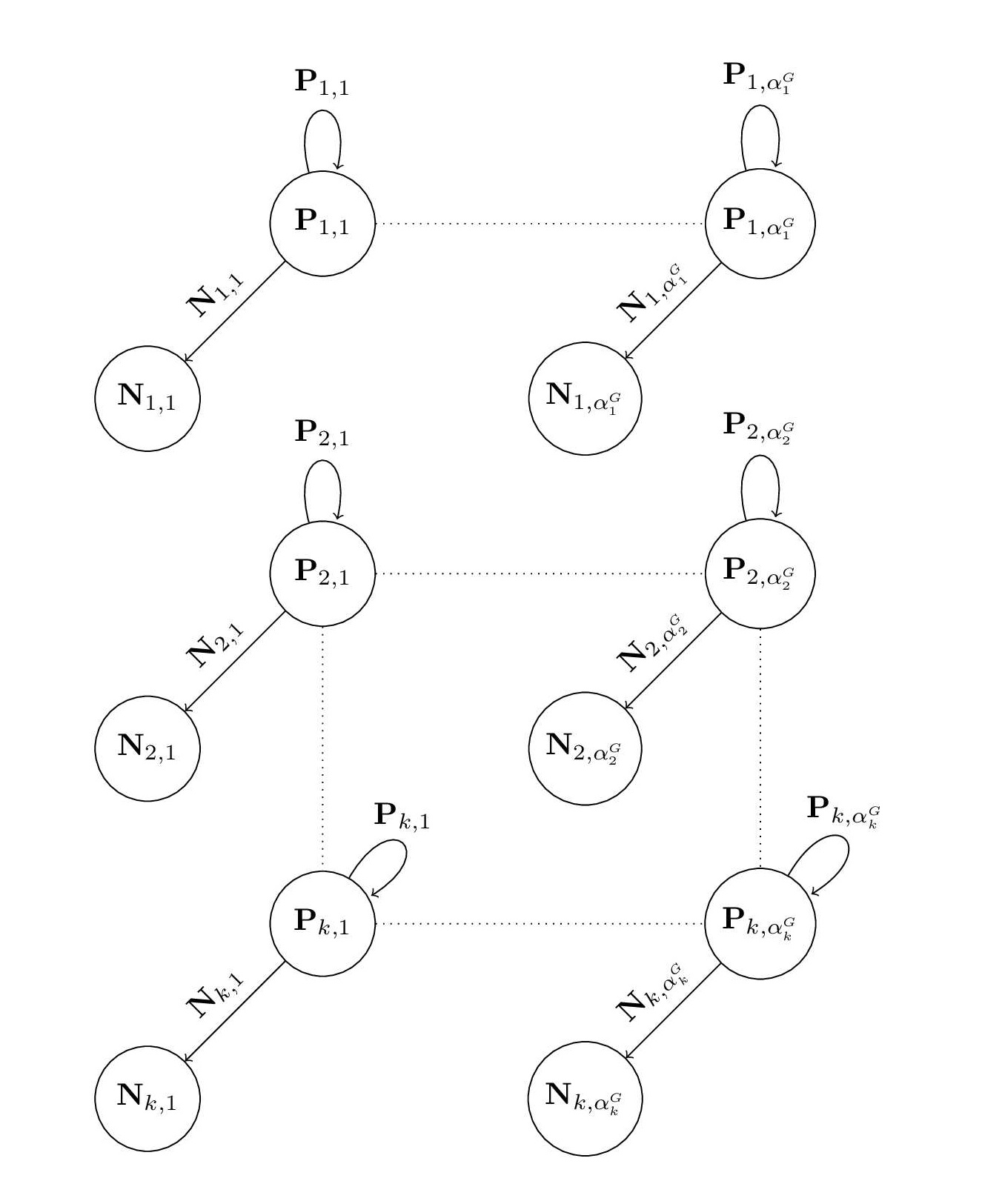}}
	\caption{The graph representation of projectors and their transitions.}\label{fig:Projector}
\end{figure}

\begin{figure}[htbp!]
	\centerline{\includegraphics[width=0.9\columnwidth,draft=false]
		{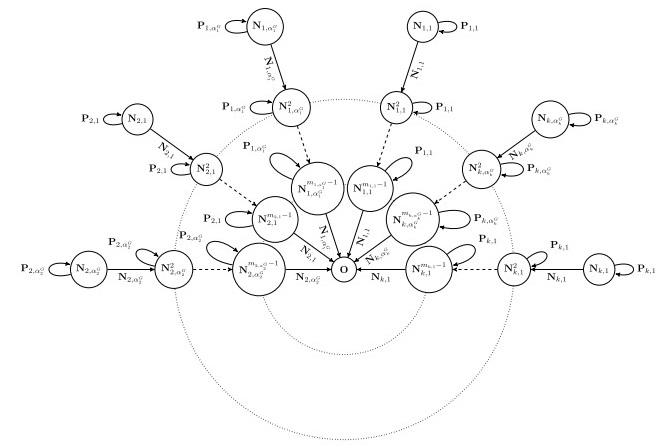}}
	\caption{The graph representation of nilpotants and their transitions.}\label{fig:Nilpotent}
\end{figure}

The resulting graph will be a directed graph with nodes $\bm{P}_{k,i}$, $\bm{N}_{k,i}$, and $\bm{0}$, where:
\begin{itemize}
\item $\bm{P}_{k,i}$ has self-loops due to their self-multiplication properties.
\item $\bm{N}_{k,i}$ will be mapped to $\bm{0}$ finally once it reaches its maximum degree of nilpotency, which is $m_{k,i}$.
\item No direct edges between $\bm{P}_{k,i}$ and $\bm{P}_{k',i'}$ or between different $\bm{N}_{k,i}$ unless the Kronecker delta conditions $k = k'$ and $i = i'$ are satisfied.
\end{itemize}
Figures~\ref{fig:Projector} and~\ref{fig:Nilpotent} demonstrate the relationships among the projection and nilpotent matrices in a structured manner, allowing for a visual and analytical exploration of their algebraic interactions.

Example~\ref{exp: 4 by 4 matrix with graph} below is provided to show an \emph{Analogous Structure Graph} (ASG) for an analogous family of matrices with dimension $4 \times 4$.

\begin{example}\label{exp: 4 by 4 matrix with graph}

Suppose we are give the following matrx $\bm{X}$ as
\begin{eqnarray}
\bm{X}&=&\begin{pmatrix}
8.23 & -1.91 & 0.09 & -2.14 \\
3.32 & 2.73 & 0.73 & -1.59 \\
1.43 & -0.23 & 2.27 & -0.66 \\
3.05 & -1.18 & -1.18 & 1.77
\end{pmatrix} \nonumber \\
&=& \begin{pmatrix}
1 & 2 & 3 & 4 \\
0 & 1 & 4 & 3 \\
2 & 0 & 1 & 1 \\
3 & 4 & 1 & 2
\end{pmatrix}\begin{pmatrix}
2 & 0 & 0 & 0 \\
0 & 3 & 0 & 0 \\
0 & 0 & 5 & 1 \\
0 & 0 & 0 & 5
\end{pmatrix}\begin{pmatrix}
-0.11& -0.05 & 0.45 & 0.07  \\
-0.16& 0.14 & -0.36 & 0.30  \\
-0.52& 0.59 & 0.09 & 0.11  \\
 0.75 & -0.5 & 0 & -0.25 
\end{pmatrix},
\end{eqnarray}
then, we have the parameters: \(K = 3\), \(\alpha_1^{\mathrm{G}} = \alpha_2^{\mathrm{G}} = 1\), and \(\alpha_3^{\mathrm{G}} = 2\) for such Jordan form. From Lemma~\ref{lma: projectors and nilpotents of a matrix}, we have
\begin{eqnarray}
\bm{X}&=&\sum\limits_{k=1}^3\sum\limits_{i=1}^{\alpha_k^{\mathrm{G}}} \lambda_k \bm{P}_{k,i} +
\sum\limits_{k=1}^3\sum\limits_{i=1}^{\alpha_k^{\mathrm{G}}} \bm{N}_{k,i}.
\end{eqnarray}

Step 1: Projection Matrices \(\bm{P}_{k,i}\)

Projection matrices \(\bm{P}_{k,i}\) correspond to the Jordan blocks. Since the problem specifies different \(\alpha_k^{\mathrm{G}}\), we construct each block accordingly.

\begin{eqnarray}
\bm{P}_{1,1}&=& \begin{pmatrix}
1 & 2 & 3 & 4 \\
0 & 1 & 4 & 3 \\
2 & 0 & 1 & 1 \\
3 & 4 & 1 & 2
\end{pmatrix}\begin{pmatrix}
1 & 0 & 0 & 0 \\
0 & 0 & 0 & 0 \\
0 & 0 & 0 & 0 \\
0 & 0 & 0 & 0
\end{pmatrix}\begin{pmatrix}
-0.11& -0.05 & 0.45 & 0.07  \\
-0.16& 0.14 & -0.36 & 0.30  \\
-0.52& 0.59 & 0.09 & 0.11  \\
 0.75 & -0.5 & 0 & -0.25 
\end{pmatrix}\nonumber \\
&=&\begin{pmatrix}
-0.11 & -0.04 & 0.45 & 0.07 \\
0 & 0 & 0 & 0 \\
-0.23 & -0.09 & 0.91 & 0.14 \\
-0.34 & -0.14 & 1.36 & 0.20
\end{pmatrix}, \nonumber \\
\bm{P}_{2,1}&=& \begin{pmatrix}
1 & 2 & 3 & 4 \\
0 & 1 & 4 & 3 \\
2 & 0 & 1 & 1 \\
3 & 4 & 1 & 2
\end{pmatrix}\begin{pmatrix}
0 & 0 & 0 & 0 \\
0 & 1 & 0 & 0 \\
0 & 0 & 0 & 0 \\
0 & 0 & 0 & 0
\end{pmatrix}\begin{pmatrix}
-0.11& -0.05 & 0.45 & 0.07  \\
-0.16& 0.14 & -0.36 & 0.30  \\
-0.52& 0.59 & 0.09 & 0.11  \\
 0.75 & -0.5 & 0 & -0.25 
\end{pmatrix}\nonumber \\
&=&\begin{pmatrix}
-0.32 & 0.27 & -0.72 & 0.59 \\
-0.16 & 0.14 & -0.36 & 0.30 \\
0 & 0 & 0 & 0 \\
-0.63 & 0.54 & -1.45 & 1.18
\end{pmatrix}.
\end{eqnarray}

For \(\alpha_3^{\mathrm{G}} = 2\), we need two projection matrices, \(\bm{P}_{3,1}\) and \(\bm{P}_{3,2}\).

\begin{eqnarray}
\bm{P}_{3,1}&=& \begin{pmatrix}
1 & 2 & 3 & 4 \\
0 & 1 & 4 & 3 \\
2 & 0 & 1 & 1 \\
3 & 4 & 1 & 2
\end{pmatrix} \begin{pmatrix}
0 & 0 & 0 & 0 \\
0 & 0 & 0 & 0 \\
0 & 0 & 1 & 0 \\
0 & 0 & 0 & 0
\end{pmatrix}\begin{pmatrix}
-0.11& -0.05 & 0.45 & 0.07  \\
-0.16& 0.14 & -0.36 & 0.30  \\
-0.52& 0.59 & 0.09 & 0.11  \\
 0.75 & -0.5 & 0 & -0.25 
\end{pmatrix}\nonumber \\
&=&\begin{pmatrix}
-1.57 & 1.77 & 0.27 & 0.34 \\
-2.09 & 2.36 & 0.36 & 0.45 \\
-0.52 & 0.59 & 0.09 & 0.11 \\
-0.52 & 0.59 & 0.09 & 0.11
\end{pmatrix},\nonumber \\
\bm{P}_{3,2}&=&\begin{pmatrix}
1 & 2 & 3 & 4 \\
0 & 1 & 4 & 3 \\
2 & 0 & 1 & 1 \\
3 & 4 & 1 & 2
\end{pmatrix}\begin{pmatrix}
0 & 0 & 0 & 0 \\
0 & 0 & 0 & 0 \\
0 & 0 & 0 & 0 \\
0 & 0 & 0 & 1
\end{pmatrix}\begin{pmatrix}
-0.11& -0.05 & 0.45 & 0.07  \\
-0.16& 0.14 & -0.36 & 0.30  \\
-0.52& 0.59 & 0.09 & 0.11  \\
 0.75 & -0.5 & 0 & -0.25 
\end{pmatrix}\nonumber \\
&=&\begin{pmatrix}
 3  &  -2 &   0 & -1 \\
 2.25 & -1.5 &  0 &   -0.75 \\
 0.75 & -0.5 &  0 &   -0.25 \\
 1.5  & -1 &   0  & -0.5
\end{pmatrix}.
\end{eqnarray}

Step 2: Nilpotent Matrices \(\bm{N}_{k,i}\)

The nilpotent matrices \(\bm{N}_{k,i}\) correspond to the generalized eigenvectors in each Jordan block. They satisfy the property that \(\bm{N}_{k,i}^2 = 0\).
\begin{eqnarray}
\bm{N}_{1,1} = \begin{pmatrix}
0 & 0 & 0 & 0 \\
0 & 0 & 0 & 0 \\
0 & 0 & 0 & 0 \\
0 & 0 & 0 & 0
\end{pmatrix}, \quad
\bm{N}_{2,1} = \begin{pmatrix}
0 & 0 & 0 & 0 \\
0 & 0 & 0 & 0 \\
0 & 0 & 0 & 0 \\
0 & 0 & 0 & 0
\end{pmatrix}
\end{eqnarray}

For \(\alpha_3^{\mathrm{G}} = 2\), the nilpotent matrix \(\bm{N}_{3,1}\) can be non-zero.
\begin{eqnarray}
\bm{N}_{3,1} &=& 
\begin{pmatrix}
2.25 & -1.5 &  0 &   -0.75\\
 3 &  -2 &  0 &  -1  \\
 0.75 & -0.5 &  0 &   -0.25 \\
 0.75 & -0.5 & 0 &  -0.25\\
\end{pmatrix},~~
\bm{N}_{3,2} = \begin{pmatrix}
0 & 0 & 0 & 0 \\
0 & 0 & 0 & 0 \\
0 & 0 & 0 & 0 \\
0 & 0 & 0 & 0
\end{pmatrix}
\end{eqnarray}

Step 3: Example Matrix \(\bm{X}\)

Using the decomposition formula, we have

\[
\bm{X} = \lambda_1 \bm{P}_{1,1} + \lambda_2 \bm{P}_{2,1} + \lambda_3 \bm{P}_{3,1} + \lambda_3 \bm{P}_{3,2} + \bm{N}_{1,1} + \bm{N}_{2,1} + \bm{N}_{3,1} + \bm{N}_{3,2}
\].

Then, we have analogous structure graph for the matrix $\bm{X}$ with the following graph data
\begin{itemize}
\item Nodes are $\bm{P}_{1,1}$, $\bm{P}_{2,1}$, $\bm{P}_{3,1}$, $\bm{P}_{3,2}$, $\bm{N}_{1,1}$, $\bm{N}_{2,1}$, $\bm{N}_{3,1}$, $\bm{N}_{3,2}$ and $\bm{0}$. 
\item Edges are $\bm{P}_{1,1}$, $\bm{P}_{2,1}$, $\bm{P}_{3,1}$, $\bm{P}_{3,2}$, $\bm{N}_{1,1}$, $\bm{N}_{2,1}$, $\bm{N}_{3,1}$, and $\bm{N}_{3,2}$. 
\item Nodes transision by edges: 
\begin{enumerate}
\item Node $\bm{P}_{1,1}$: $\bm{P}^2_{1,1}=\bm{P}_{1,1}$, $\bm{P}_{1,1}\bm{N}_{1,1}=\bm{N}_{1,1}\bm{P}_{1,1}=\bm{N}_{1,1}$, 
\item Node $\bm{P}_{2,1}$: $\bm{P}^2_{2,1}=\bm{P}_{2,1}$, $\bm{P}_{2,1}\bm{N}_{2,1}=\bm{N}_{2,1}\bm{P}_{2,1}=\bm{N}_{2,1}$, 
\item Node $\bm{P}_{3,1}$: $\bm{P}^2_{3,1}=\bm{P}_{3,1}$, $\bm{P}_{3,1}\bm{N}_{3,1}=\bm{N}_{3,1}\bm{P}_{3,1}=\bm{N}_{3,1}$, 
\item Node $\bm{P}_{3,2}$: $\bm{P}^2_{3,2}=\bm{P}_{3,2}$, $\bm{P}_{3,2}\bm{N}_{3,2}=\bm{N}_{3,2}\bm{P}_{3,2}=\bm{N}_{3,2}$, 
\item Node $\bm{N}_{1,1}$: $\bm{N}_{1,1}=\bm{0}$, $\bm{P}_{1,1}\bm{N}_{1,1}=\bm{N}_{1,1}\bm{P}_{1,1}=\bm{N}_{1,1}$, 
\item Node $\bm{N}_{2,1}$: $\bm{N}_{2,1}=\bm{0}$, $\bm{P}_{2,1}\bm{N}_{2,1}=\bm{N}_{2,1}\bm{P}_{2,1}=\bm{N}_{2,1}$, 
\item Node $\bm{N}_{3,1}$: $\bm{N}^2_{3,1}=\bm{0}$, $\bm{P}_{3,1}\bm{N}_{3,1}=\bm{N}_{3,1}\bm{P}_{3,1}=\bm{N}_{3,1}$, 
\item Node $\bm{N}_{3,2}$: $\bm{N}_{3,2}=\bm{0}$, $\bm{P}_{3,2}\bm{N}_{3,2}=\bm{N}_{3,2}\bm{P}_{3,2}=\bm{N}_{3,2}$. 
\end{enumerate}
\end{itemize}
\end{example}

\subsection{Enumeration of Finite-Dimensional Operator Types}\label{sec: Enumeration of Finite-Dimensional Operator Types}

In this section, we will evaluate the number of analogous families for a given matrix with the dimension $m$. We begin by defining two counting functoins $\mathrm{P}_k(m)$ and $\mathrm{P}(m)$. The function $\mathrm{P}_k(m)$ will return the number of the following integer solutions of $x_1, x_2, \ldots x_k$:
\begin{eqnarray}\label{eq3: enum}
x_1 + x_2 + \ldots x_k &=&m,
\end{eqnarray}
where $m$ is a nonnegative natural number with $x_1 \geq x_2 \geq \ldots \geq x_k$ and $x_k \geq 1$. The function $\mathrm{P}_k(m)$ will return the number of partitions of $m$ into $k$ parts. For all $m$, we have $\mathrm{P}_1(m)=1$. We also have the following recurrence relation for $\mathrm{P}_k(m)$:
\begin{eqnarray}\label{eq4: enum}
\mathrm{P}_k(m)&=&\mathrm{P}_k(m-k) + \mathrm{P}_{k-1}(m-1).
\end{eqnarray}

On the other hand, we will define the function $\mathrm{P}(m)$ to evaluate the number of the following integer solutions of $x_1, x_2, \ldots x_k$:
\begin{eqnarray}\label{eq5: enum}
x_1 + x_2 + \ldots x_k &=&m,
\end{eqnarray}
where $k=1,2,\ldots,m$ with $x_1 \geq x_2 \geq \ldots \geq x_k$ and $x_k \geq 0$. By comparing the counting relations between Eq.~\eqref{eq3: enum} and Eq.~\eqref{eq5: enum} Then, we have 
\begin{eqnarray}\label{eq6: enum}
\mathrm{P}(m)&=&\sum\limits_{k=1}^m\mathrm{P}_k(m). 
\end{eqnarray}

Recall we have the following relatonships about $\alpha_k^{(\mathrm{A})}$ and $\alpha_k^{(\mathrm{G})}$:
\begin{eqnarray}\label{eq1: enum}
\sum\limits_{k=1}^K \alpha_k^{(\mathrm{A})}&=&m,
\end{eqnarray}
and 
\begin{eqnarray}\label{eq2: enum}
\sum\limits_{i=1}^{\alpha_k^{(\mathrm{G})}}m_{k,i}&=&\alpha_k^{(\mathrm{A})}.
\end{eqnarray}

For a $m \times m$ matrix with $K$ distinct eigenvalues and their algebraic multiplicities $\alpha_k^{(\mathrm{A})}$ for $k=1,2,\ldots,K$, we can count the number of different Jordan bocks by (multiplication principal) via Eq.~\eqref{eq2: enum}
\begin{eqnarray}\label{eq7: enum}
\prod\limits_{k=1}^{K}\mathrm{P}\left(\alpha_k^{(A)}\right).
\end{eqnarray}
Because the value of $K$ can be ranged from $1$ to $m$, the number of analogous families for a given matrix with the dimension $m$ can be expressed by (addition principal) 
\begin{eqnarray}\label{eq8: enum}
\sum\limits_{K=1}^m\left[\sum\limits_{\sum\limits_{k=1}^K \alpha_k^{(\mathrm{A})}=m}\prod\limits_{k=1}^{K}\mathrm{P}\left(\alpha_k^{(A)}\right)\right].
\end{eqnarray}
Note that there are $\mathrm{P}_K(m)$ valid solutions of $\alpha_k^{(\mathrm{A})}$ in $\sum\limits_{k=1}^K \alpha_k^{(\mathrm{A})}=m$.

\subsection{Spectral Mapping Theorem for Multivariate Finite-Dimensional Operators}\label{sec: Spectral Mapping Theorem for Multivariate Finite-Dimensional Operators}

In this section, we will establish spectral mapping theorem for multivariate finite-dimensional operators. We adopt power series method used in~\cite{ashida2020non}, however, they did not provide condition when those functions can be expressed by power series approach. 

\subsubsection{Single Variable}\label{sec: Single Variable}

We consider a single-variable complex function \( f(z) = \sum\limits_{i=0}^{\infty} a_i z^i \) represents a power series. The conditions under which this series is valid, i.e., when it converges and defines a well-behaved function, depend on the radius of convergence of the power series~\cite{pathak2015functions}.

The radius of convergence, $R$ determines the region within which the series converges. It is given by the formula:
\begin{eqnarray}\label{eq: conv R root}
\frac{1}{R} = \limsup_{i \to \infty} \sqrt[i]{|a_i|}
\end{eqnarray}
or equivalently, using the root test or ratio test:
\begin{eqnarray}\label{eq: conv R ratio}
     \frac{1}{R} = \lim_{i \to \infty} \left| \frac{a_{i+1}}{a_i} \right|
\end{eqnarray}
if the limit exists. The function $f(z) = \sum\limits_{i=0}^{\infty} a_i z^i$ represents a power series that converges within a disk of radius $R$, determined by the root test given by Eq.~\eqref{eq: conv R root} or ratio test given by Eq.~\eqref{eq: conv R ratio}. The series defines an analytic function inside the disk $|z| < R$. 

\begin{theorem}\label{thm: Spectral Mapping Theorem for Single Variable}
Given an analytic function $f(z)$ within the domain for $|z| < R$, a matrix $\bm{X}$ with the dimension $m$ and $K$ distinct eigenvalues $\lambda_k$ for $k=1,2,\ldots,K$ such that
\begin{eqnarray}\label{eq1: thm: Spectral Mapping Theorem for Single Variable}
\bm{X}&=&\sum\limits_{k=1}^K\sum\limits_{i=1}^{\alpha_k^{\mathrm{G}}} \lambda_k \bm{P}_{k,i}+
\sum\limits_{k=1}^K\sum\limits_{i=1}^{\alpha_k^{\mathrm{G}}} \bm{N}_{k,i},
\end{eqnarray}
where $\left\vert\lambda_k\right\vert<R$, then, we have
\begin{eqnarray}\label{eq2: thm: Spectral Mapping Theorem for Single Variable}
f(\bm{X})&=&\sum\limits_{k=1}^K \left[\sum\limits_{i=1}^{\alpha_k^{(\mathrm{G})}}f(\lambda_k)\bm{P}_{k,i}+\sum\limits_{i=1}^{\alpha_k^{(\mathrm{G})}}\sum\limits_{q=1}^{m_{k,i}-1}\frac{f^{(q)}(\lambda_k)}{q!}\bm{N}_{k,i}^q\right].
\end{eqnarray}
\end{theorem}
\textbf{Proof:}
Recall we have
\begin{eqnarray}\label{eq2-1: thm: Spectral Mapping Theorem for Single Variable}
\bm{P}_{k,i}\bm{P}_{k',i'}&=&\bm{P}_{k,i}\delta(k,k')\delta(i,i'), \nonumber \\
\bm{P}_{k',i'}\bm{N}_{k,i}&=&\bm{N}_{k,i}\bm{P}_{k',i'}=\bm{N}_{k,i}\delta(k,k')\delta(i,i'), \nonumber \\
\bm{N}_{k,i}\bm{N}_{k',i'}&=&\bm{N}^2_{k,i}\delta(k,k')\delta(i,i').
\end{eqnarray}   

Because $f(z)$ is an analytic function, we have
\begin{eqnarray}\label{eq3: thm: Spectral Mapping Theorem for Single Variable}
f(\bm{X})&=&\sum\limits_{\ell=0}^{\infty}a_\ell \bm{X}^\ell \nonumber \\
&=&\sum\limits_{\ell=0}^{\infty}a_\ell\left(\sum\limits_{k=1}^K\sum\limits_{i=1}^{\alpha_k^{\mathrm{G}}} \lambda_k \bm{P}_{k,i}+
\sum\limits_{k=1}^K\sum\limits_{i=1}^{\alpha_k^{\mathrm{G}}} \bm{N}_{k,i}\right)^\ell \nonumber \\
&=_1&\sum\limits_{\ell=0}^{\infty}\sum\limits_{k=1}^K\sum\limits_{i=1}^{\alpha_k^{\mathrm{G}}} a_\ell\left(\lambda_k \bm{P}_{k,i}+\bm{N}_{k,i}\right)^\ell  \nonumber \\
&=& \sum\limits_{k=1}^K\sum\limits_{i=1}^{\alpha_k^{\mathrm{G}}}\sum\limits_{\ell=0}^{\infty}\sum\limits_{q=0}^{\ell} \frac{a_\ell \ell!}{q! (\ell-q)!}\lambda^{\ell-q}_k\bm{P}^{\ell-q}_{k,i}\bm{N}^q_{k,i}\nonumber \\
&=& \sum\limits_{k=1}^K\sum\limits_{i=1}^{\alpha_k^{\mathrm{G}}}\sum\limits_{\ell=0}^{\infty}\left[a_{\ell}\lambda_k^{\ell}\bm{P}_{k,i}+\sum\limits_{q=1}^{\ell} \frac{a_\ell \ell!}{q! (\ell-q)!}\lambda^{\ell-q}_k\bm{P}^{\ell-q}_{k,i}\bm{N}^q_{k,i}\right]\nonumber \\
&=_2&\sum\limits_{k=1}^K \left[\sum\limits_{i=1}^{\alpha_k^{(\mathrm{G})}}f(\lambda_k)\bm{P}_{k,i}+\sum\limits_{i=1}^{\alpha_k^{(\mathrm{G})}}\sum\limits_{q=1}^{m_{k,i}-1}\frac{f^{(q)}(\lambda_k)}{q!}\bm{N}_{k,i}^q\right]. 
\end{eqnarray}
where we apply Eq.~\eqref{eq2-1: thm: Spectral Mapping Theorem for Single Variable} in $=_1$ and $=_2$, and $f^{(q)}(z) = \sum\limits_{\ell=q}^{\infty}\frac{a_\ell \ell!}{(\ell-q)!}z^{\ell-q}$ in $=_2$.
$\hfill \Box$

\subsubsection{Multiple Variables}\label{sec: Multiple Variables}

Let \( f(z_1, z_2, \dots, z_n) \) be a function of \( n \) complex variables \( z_1, z_2, \dots, z_n \). A multivariate power series expansion around a point (typically the origin) is of the form:
\begin{eqnarray}
f(z_1, z_2, \dots, z_n)&=&\sum_{i_1, i_2, \dots, i_n = 0}^{\infty} a_{i_1, i_2, \dots,i_n} z_1^{i_1} z_2^{i_2} \dots z_n^{i_n},
\end{eqnarray}
where the coefficients \( a_{i_1,i_2, \dots i_n} \) are constants, and the exponents \( i_1, i_2, \dots, i_n \) are non-negative integers.

The power series converges within a certain domain in \( \mathbb{C}^n \). The convergence is generally within a polydisk or domain of convergence, defined by:
\begin{eqnarray}
   D(R_1, R_2, \dots, R_n)&=&\{ (z_1, z_2, \dots, z_n) \in \mathbb{C}^n : |z_1| < R_1, |z_2| < R_2, \dots, |z_n| < R_n \},
\end{eqnarray}
where \( R_1, R_2, \dots, R_n \) are positive real numbers called the radius of convergence in each direction.

The multivariate radius of convergence for each variable can be determined similarly to the single-variable case, using a generalization of the ratio or root test. For instance, for each variable \( z_j \), the series must satisfy:
\begin{eqnarray}
   \limsup_{i_j \to \infty} \sqrt[i_j]{|a_{i_1 i_2 \dots i_n}|} = \frac{1}{R_j}.
\end{eqnarray}

The function \( f(z_1, z_2, \dots, z_n) \) must be analytic (holomorphic) in each variable within the region of convergence. This means that for each fixed set of values of \( z_2, z_3, \dots, z_n \), the function \( f \) is analytic in \( z_1 \), and similarly for all other variables. This property is called separate analyticity.

We first consider spectral mapping theorem for two input matrices in Theorem~\ref{thm: Spectral Mapping Theorem for Two Variables}.

\begin{theorem}\label{thm: Spectral Mapping Theorem for Two Variables}
Given an analytic function $f(z_1,z_2)$ within the domain for $|z_1| < R_1$ and $|z_2| < R_2$, the first matrix $\bm{X}_1$ with the dimension $m$ and $K_1$ distinct eigenvalues $\lambda_{k_1}$ for $k_1=1,2,\ldots,K_1$ such that
\begin{eqnarray}\label{eq1-1: thm: Spectral Mapping Theorem for Two Variables}
\bm{X}_1&=&\sum\limits_{k_1=1}^{K_1}\sum\limits_{i_1=1}^{\alpha_{k_1}^{\mathrm{G}}} \lambda_{k_1} \bm{P}_{k_1,i_1}+
\sum\limits_{k_1=1}^{K_1}\sum\limits_{i_1=1}^{\alpha_{k_1}^{\mathrm{G}}} \bm{N}_{k_1,i_1},
\end{eqnarray}
where $\left\vert\lambda_{k_1}\right\vert<R_1$, and second matrix $\bm{X}_2$ with the dimension $m$ and $K_2$ distinct eigenvalues $\lambda_{k_2}$ for $k_2=1,2,\ldots,K_2$ such that
\begin{eqnarray}\label{eq1-2: thm: Spectral Mapping Theorem for Two Variables}
\bm{X}_2&=&\sum\limits_{k_2=1}^{K_2}\sum\limits_{i_2=1}^{\alpha_{k_2}^{\mathrm{G}}} \lambda_{k_2} \bm{P}_{k_2,i_2}+
\sum\limits_{k_2=1}^{K_2}\sum\limits_{i_2=1}^{\alpha_{k_2}^{\mathrm{G}}} \bm{N}_{k_2,i_2},
\end{eqnarray}
where $\left\vert\lambda_{k_2}\right\vert<R_2$.

Then, we have
\begin{eqnarray}\label{eq2: thm: Spectral Mapping Theorem for Two Variables}
f(\bm{X}_1, \bm{X}_2)&=&\sum\limits_{k_1=1}^{K_1}\sum\limits_{k_2=1}^{K_2}\sum\limits_{i_1=1}^{\alpha_{k_1}^{(\mathrm{G})}}\sum\limits_{i_2=1}^{\alpha_{k_2}^{(\mathrm{G})}}f(\lambda_{k_1}, \lambda_{k_2})\bm{P}_{k_1,i_1}\bm{P}_{k_2,i_2} \nonumber \\
&&+\sum\limits_{k_1=1}^{K_1}\sum\limits_{k_2=1}^{K_2}\sum\limits_{i_1=1}^{\alpha_{k_1}^{(\mathrm{G})}}\sum\limits_{i_2=1}^{\alpha_{k_2}^{(\mathrm{G})}}\sum_{q_2=1}^{m_{k_2,i_2}-1}\frac{f^{(-,q_2)}(\lambda_{k_1},\lambda_{k_2})}{q_2!}\bm{P}_{k_1,i_1}\bm{N}_{k_2,i_2}^{q_2} \nonumber \\
&&+\sum\limits_{k_1=1}^{K_1}\sum\limits_{k_2=1}^{K_2}\sum\limits_{i_1=1}^{\alpha_{k_1}^{(\mathrm{G})}}\sum\limits_{i_2=1}^{\alpha_{k_2}^{(\mathrm{G})}}\sum_{q_1=1}^{m_{k_1,i_1}-1}\frac{f^{(q_1,-)}(\lambda_{k_1},\lambda_{k_2})}{q_1!}\bm{N}_{k_1,i_1}^{q_1}\bm{P}_{k_2,i_2}  \nonumber \\
&&+\sum\limits_{k_1=1}^{K_1}\sum\limits_{k_2=1}^{K_2}\sum\limits_{i_1=1}^{\alpha_{k_1}^{(\mathrm{G})}}\sum\limits_{i_2=1}^{\alpha_{k_2}^{(\mathrm{G})}}\sum_{q_1=1}^{m_{k_1,i_1}-1}\sum_{q_2=1}^{m_{k_2,i_2}-1}\frac{f^{(q_1,q_2)}(\lambda_{k_1},\lambda_{k_2})}{q_1!q_2!}\bm{N}_{k_1,i_1}^{q_1}\bm{N}_{k_2,i_2}^{q_2}
\end{eqnarray}
\end{theorem}
\textbf{Proof:}
Recall we still have
\begin{eqnarray}\label{eq2-1: thm: Spectral Mapping Theorem for Two Variables}
\bm{P}_{k_1,i_1}\bm{P}_{k'_1,i'_1}&=&\bm{P}_{k_1,i_1}\delta(k_1,k'_1)\delta(i_1,i'_1), \nonumber \\
\bm{P}_{k'_1,i'_1}\bm{N}_{k_1,i_1}&=&\bm{N}_{k_1,i_1}\bm{P}_{k'_1,i'_1}=\bm{N}_{k_1,i_1}\delta(k_1,k'_1)\delta(i_1,i'_1), \nonumber \\
\bm{N}_{k_1,i_1}\bm{N}_{k'_1,i'_1}&=&\bm{N}^2_{k_1,i_1}\delta(k_1,k'_1)\delta(i_1,i'_1).
\end{eqnarray}   
\begin{eqnarray}\label{eq2-2: thm: Spectral Mapping Theorem for Two Variables}
\bm{P}_{k_2,i_2}\bm{P}_{k'_2,i'_2}&=&\bm{P}_{k_2,i_2}\delta(k_2,k'_2)\delta(i_2,i'_2), \nonumber \\
\bm{P}_{k'_2,i'_2}\bm{N}_{k_2,i_2}&=&\bm{N}_{k_2,i_2}\bm{P}_{k'_2,i'_2}=\bm{N}_{k_2,i_2}\delta(k_2,k'_2)\delta(i_2,i'_2), \nonumber \\
\bm{N}_{k_2,i_2}\bm{N}_{k'_2,i'_2}&=&\bm{N}^2_{k_2,i_2}\delta(k_2,k'_2)\delta(i_2,i'_2).
\end{eqnarray}   

Because $f(z)$ is an analytic function, we have
\begin{eqnarray}\label{eq3: thm: Spectral Mapping Theorem for Two Variables}
f(\bm{X}_1,\bm{X}_2)&=&\sum\limits_{\ell_1=0,\ell_2=0}^{\infty}a_{\ell_1,\ell_2} \bm{X}_1^{\ell_1}\bm{X}_2^{\ell_2}\nonumber \\
&=&\sum\limits_{\ell_1=0,\ell_2=0}^{\infty}a_{\ell_1,\ell_2}\left(\sum\limits_{k_1=1}^{K_1}\sum\limits_{i_1=1}^{\alpha_{k_1}^{\mathrm{G}}} \lambda_{k_1} \bm{P}_{k_1,i_1}+
\sum\limits_{k_1=1}^{K_1}\sum\limits_{i_1=1}^{\alpha_{k_1}^{\mathrm{G}}} \bm{N}_{k_1,i_1}\right)^{\ell_1}\nonumber \\
&&\times \left(\sum\limits_{k_2=1}^{K_2}\sum\limits_{i_2=1}^{\alpha_{k_2}^{\mathrm{G}}} \lambda_{k_2} \bm{P}_{k_2,i_2}+
\sum\limits_{k_2=1}^{K_2}\sum\limits_{i_2=1}^{\alpha_{k_2}^{\mathrm{G}}} \bm{N}_{k_2,i_2}\right)^{\ell_2}\nonumber \\
&=_1&\sum\limits_{\ell_1=0,\ell_2=0}^{\infty}a_{\ell_1,\ell_2}\left[\sum\limits_{k_1=1}^{K_1}\sum\limits_{i_1=1}^{\alpha_{k_1}^{\mathrm{G}}}\left(\lambda_{k_1} \bm{P}_{k_1,i_1}+\bm{N}_{k_1,i_1}\right)^{\ell_1}\right]\left[\sum\limits_{k_2=1}^{K_2}\sum\limits_{i_2=1}^{\alpha_{k_2}^{\mathrm{G}}}\left(\lambda_{k_2} \bm{P}_{k_2,i_2}+\bm{N}_{k_2,i_2}\right)^{\ell_2}\right]\nonumber \\
&=& \sum\limits_{k_1=1}^{K_1}\sum\limits_{k_2=1}^{K_2}\sum\limits_{i_1=1}^{\alpha_{k_1}^{(\mathrm{G})}}\sum\limits_{i_2=1}^{\alpha_{k_2}^{(\mathrm{G})}}\sum\limits_{\ell_1=0,\ell_2=0}^{\infty}a_{\ell_1,\ell_2}\left[\sum\limits_{q_1=0}^{m_{k_1,i_1}-1}\frac{\ell_1!}{q_1! (\ell_1-q_1)!}\lambda^{\ell_1-q_1}_{k_1}\bm{P}^{\ell_1-q_1}_{k_1,i_1}\bm{N}^{q_1}_{k_1,i_1}\right]\nonumber \\
&&\times \left[\sum\limits_{q_2=0}^{m_{k_2,i_2}-1}\frac{\ell_2!}{q_2! (\ell_2-q_2)!}\lambda^{\ell_2-q_2}_{k_2}\bm{P}^{\ell_2-q_2}_{k_2,i_2}\bm{N}^{q_2}_{k_2,i_2}\right]\nonumber \\
&=& \sum\limits_{k_1=1}^{K_1}\sum\limits_{k_2=1}^{K_2}\sum\limits_{i_1=1}^{\alpha_{k_1}^{(\mathrm{G})}}\sum\limits_{i_2=1}^{\alpha_{k_2}^{(\mathrm{G})}}\sum\limits_{\ell_1=0,\ell_2=0}^{\infty}a_{\ell_1,\ell_2}\left[\lambda_{k_1}^{\ell_1}\bm{P}^{\ell_1}_{k_1,i_1}+\sum\limits_{q_1=1}^{m_{k_1,i_1}-1}\frac{\ell_1!}{q_1! (\ell_1-q_1)!}\lambda^{\ell_1-q_1}_{k_1}\bm{P}^{\ell_1-q_1}_{k_1,i_1}\bm{N}^{q_1}_{k_1,i_1}\right]\nonumber \\
&&\times \left[\lambda_{k_2}^{\ell_2}\bm{P}^{\ell_2}_{k_2,i_2}+\sum\limits_{q_2=1}^{m_{k_2,i_2}-1}\frac{\ell_2!}{q_2! (\ell_2-q_2)!}\lambda^{\ell_2-q_2}_{k_2}\bm{P}^{\ell_2-q_2}_{k_2,i_2}\bm{N}^{q_2}_{k_2,i_2}\right]\nonumber \\
&=_2&\sum\limits_{k_1=1}^{K_1}\sum\limits_{k_2=1}^{K_2}\sum\limits_{i_1=1}^{\alpha_{k_1}^{(\mathrm{G})}}\sum\limits_{i_2=1}^{\alpha_{k_2}^{(\mathrm{G})}}f(\lambda_{k_1}, \lambda_{k_2})\bm{P}_{k_1,i_1}\bm{P}_{k_2,i_2} \nonumber \\
&&+\sum\limits_{k_1=1}^{K_1}\sum\limits_{k_2=1}^{K_2}\sum\limits_{i_1=1}^{\alpha_{k_1}^{(\mathrm{G})}}\sum\limits_{i_2=1}^{\alpha_{k_2}^{(\mathrm{G})}}\sum_{q_2=1}^{m_{k_2,i_2}-1}\frac{f^{(-,q_2)}(\lambda_{k_1},\lambda_{k_2})}{q_2!}\bm{P}_{k_1,i_1}\bm{N}_{k_2,i_2}^{q_2} \nonumber \\
&&+\sum\limits_{k_1=1}^{K_1}\sum\limits_{k_2=1}^{K_2}\sum\limits_{i_1=1}^{\alpha_{k_1}^{(\mathrm{G})}}\sum\limits_{i_2=1}^{\alpha_{k_2}^{(\mathrm{G})}}\sum_{q_1=1}^{m_{k_1,i_1}-1}\frac{f^{(q_1,-)}(\lambda_{k_1},\lambda_{k_2})}{q_1!}\bm{N}_{k_1,i_1}^{q_1}\bm{P}_{k_2,i_2}  \nonumber \\
&&+\sum\limits_{k_1=1}^{K_1}\sum\limits_{k_2=1}^{K_2}\sum\limits_{i_1=1}^{\alpha_{k_1}^{(\mathrm{G})}}\sum\limits_{i_2=1}^{\alpha_{k_2}^{(\mathrm{G})}}\sum_{q_1=1}^{m_{k_1,i_1}-1}\sum_{q_2=1}^{m_{k_2,i_2}-1}\frac{f^{(q_1,q_2)}(\lambda_{k_1},\lambda_{k_2})}{q_1!q_2!}\bm{N}_{k_1,i_1}^{q_1}\bm{N}_{k_2,i_2}^{q_2}. 
\end{eqnarray}
where we apply Eq.~\eqref{eq2-1: thm: Spectral Mapping Theorem for Two Variables} and Eq.~\eqref{eq2-2: thm: Spectral Mapping Theorem for Two Variables} in $=_1$ and $=_2$, and following relations about the partial derivatives of $f(z_1,z_2)$ in $=_2$:
\begin{eqnarray}\label{eq4: thm: Spectral Mapping Theorem for Two Variables}
f^{(-,q_2)}(z_1,z_2)&=&\sum_{\ell_1=0}^{\infty}\sum_{\ell_2=q_2}^{\infty}\frac{a_{\ell_1,\ell_2}\ell_2!}{(\ell_2-q_2)!}z_1^{\ell_1}z_2^{\ell_2-q_2},\nonumber\\
f^{(q_1,-)}(z_1,z_2)&=&\sum_{\ell_1=q_1}^{\infty}\sum_{\ell_2=0}^{\infty}\frac{a_{\ell_1,\ell_2}\ell_1!}{(\ell_1-q_1)!}z_1^{\ell_1-q_1}z_2^{\ell_2},\nonumber\\
f^{(q_1,q_2)}(z_1,z_2)&=&\sum_{\ell_1=q_1}^{\infty}\sum_{\ell_2=q_2}^{\infty}\frac{a_{\ell_1,\ell_2}\ell_1!\ell_2!}{(\ell_1-q_1)!(\ell_2-q_2)!}z_1^{\ell_1-q_1}z_2^{\ell_2-q_2}.
\end{eqnarray}
$\hfill \Box$

Before we present the spectralmapping theorem for any number of input matrices, we have to introduce several special ntations. Given $r$ positive integers $q_1,q_2,\ldots,q_r$, we define $\alpha_{\kappa}(q_1,\ldots,q_r)$ to be the selection of these $r$ arguments $q_1,\ldots,q_r$ to $\kappa$ arguments, i.e., we have
\begin{eqnarray}
\alpha_{\kappa}(q_1,\ldots,q_r)&=&\{q_{\iota_1},q_{\iota_2},\ldots,q_{\iota_\kappa}\}.
\end{eqnarray}
We use $\mbox{Ind}(\alpha_{\kappa}(q_1,\ldots,q_r))$ to obtain indices of those $\kappa$ positive integers $\{q_{\iota_1},q_{\iota_2},\ldots,q_{\iota_\kappa}\}$, i.e., we have
\begin{eqnarray}
\mbox{Ind}(\alpha_{\kappa}(q_1,\ldots,q_r))&=&\{\iota_1,\iota_2,\ldots,\iota_\kappa\}.
\end{eqnarray}
We use $\alpha_{\kappa}(q_1,\ldots,q_r)=1$ to represent $q_{\iota_1}=1,q_{\iota_2}=1,\ldots,q_{\iota_\kappa}=1$. We also use \\
$m_{k_{\mbox{Ind}(\alpha_{\kappa}(q_1,\ldots,q_r))},i_{\mbox{Ind}(\alpha_{\kappa}(q_1,\ldots,q_r))}}-1$ to represent $m_{k_{\iota_1},i_{\iota_1}}-1,m_{k_{\iota_2},i_{\iota_2}}-1,\ldots,m_{k_{\iota_\kappa},i_{\iota_\kappa}}-1$. 

\iffalse
\begin{eqnarray}
\prod\limits_{\substack{\beta =\mbox{Ind}(\alpha_{\kappa}(q_1,\ldots,q_r)), \bm{Y}=\bm{N}\\ \beta \neq \mbox{Ind}(\alpha_{\kappa}(q_1,\ldots,q_r)), \bm{Y}=\bm{P}} 
}^{r} \bm{Y}_{k_\beta,i_\beta} 
\end{eqnarray}
\fi

We are ready to present Theorem~\ref{thm: Spectral Mapping Theorem for r Variables} about spectral mapping theorem for $r$ matrices. 
\begin{theorem}\label{thm: Spectral Mapping Theorem for r Variables}
Given an analytic function $f(z_1,z_2,\ldots,z_r)$ within the domain for $|z_l| < R_l$, and the matrix $\bm{X}_l$ with the dimension $m$ and $K_l$ distinct eigenvalues $\lambda_{k_l}$ for $k_l=1,2,\ldots,K_l$ such that
\begin{eqnarray}\label{eq1-1: thm: Spectral Mapping Theorem for r Variables}
\bm{X}_l&=&\sum\limits_{k_l=1}^{K_l}\sum\limits_{i_l=1}^{\alpha_{k_l}^{\mathrm{G}}} \lambda_{k_l} \bm{P}_{k_l,i_l}+
\sum\limits_{k_l=1}^{K_l}\sum\limits_{i_l=1}^{\alpha_{k_l}^{\mathrm{G}}} \bm{N}_{k_l,i_l},
\end{eqnarray}
where $\left\vert\lambda_{k_l}\right\vert<R_l$ for $l=1,2,\ldots,r$.

Then, we have
\begin{eqnarray}\label{eq2: thm: Spectral Mapping Theorem for kappa Variables}
\lefteqn{f(\bm{X}_1,\ldots,\bm{X}_r)=}\nonumber \\
&& \sum\limits_{k_1=\ldots=k_r=1}^{K_1,\ldots,K_r} \sum\limits_{i_1=\ldots=i_r=1}^{\alpha_{k_1}^{(\mathrm{G})},\ldots,\alpha_{k_r}^{(\mathrm{G})}}
f(\lambda_{k_1},\ldots,\lambda_{k_r})\bm{P}_{k_1,i_1}\ldots\bm{P}_{k_r,i_r}\nonumber \\
&&+\sum\limits_{k_1=\ldots=k_r=1}^{K_1,\ldots,K_r} \sum\limits_{i_1=\ldots=i_r=1}^{\alpha_{k_1}^{(\mathrm{G})},\ldots,\alpha_{k_r}^{(\mathrm{G})}}\sum\limits_{\kappa=1}^{r-1}\sum\limits_{\alpha_\kappa(q_1,\ldots,q_r)}\Bigg(\sum\limits_{\alpha_{\kappa}(q_1,\ldots,q_r)=1}^{m_{k_{\mbox{Ind}(\alpha_{\kappa}(q_1,\ldots,q_r))},i_{\mbox{Ind}(\alpha_{\kappa}(q_1,\ldots,q_r))}}-1}\nonumber \\
&&~~~~~ \frac{f^{\alpha_{\kappa}(q_1,\ldots,q_r)}(\lambda_{k_1},\ldots,\lambda_{k_r})}{q_{\iota_1}!q_{\iota_2}!\ldots q_{\iota_\kappa}!}\times \prod\limits_{\substack{\beta =\mbox{Ind}(\alpha_{\kappa}(q_1,\ldots,q_r)), \bm{Y}=\bm{N}^{q_\beta}_{k_\beta,i_\beta} \\ \beta \neq \mbox{Ind}(\alpha_{\kappa}(q_1,\ldots,q_r)), \bm{Y}=\bm{P}_{k_\beta,i_\beta} }
}^{r} \bm{Y}\Bigg) 
\nonumber \\
&&+\sum\limits_{k_1=\ldots=k_r=1}^{K_1,\ldots,K_r} \sum\limits_{i_1=\ldots=i_r=1}^{\alpha_{k_1}^{(\mathrm{G})},\ldots,\alpha_{k_r}^{(\mathrm{G})}} \sum\limits_{q_1=\ldots=q_r=1}^{m_{k_1,i_1}-1,\ldots,m_{k_r,i_r}-1}
\frac{f^{(q_1,\ldots,q_r)}(\lambda_{k_1},\ldots,\lambda_{k_r})}{q_1!\cdots q_r!}\bm{N}^{q_1}_{k_1,i_1}\ldots\bm{N}^{q_r}_{k_r,i_r}
\end{eqnarray}
where we have
\begin{itemize}
\item $\sum\limits_{\alpha_{\kappa}(q_1,\ldots,q_r)}$ is the summation running over all selection of $\alpha_{\kappa}(q_1,\ldots,q_r)$ given $\kappa$;
\item $f^{\alpha_{\kappa}(q_1,\ldots,q_r)}(\lambda_1,\ldots,\lambda_r)$ represents the partial derivatives with respect to variables with indices $\iota_1,\iota_2,\ldots,\iota_\kappa$ and the orders of derivatives given by $q_{\iota_1},q_{\iota_2},\ldots,q_{\iota_\kappa}$. 
\end{itemize}
\end{theorem}
\textbf{Proof:}
The proof follows directly from Theorem~\ref{thm: Spectral Mapping Theorem for Two Variables}, with the remaining steps consisting primarily of routine notational manipulations.
$\hfill\Box$

\section{Countable Infinite-Dimensional Operators}\label{sec: Countable Infinite-Dimensional Operators}

\subsection{Jordan Decomposition for Countable Infinite-Dimensional Operators}\label{sec: Jordan Decomposition for Countable Infinite-Dimensional Operators}

In this section, we consider a square matrix with countable infinite dimension $\bm{X} \in \mathbb{C}^{\infty \times \infty}$, which can be expressed by  
\begin{eqnarray}\label{eq: Jordan decompostion count-inf}
\bm{X}&=& \bm{U}\left(\bigoplus\limits_{k=1}^{\infty}\bigoplus\limits_{i=1}^{\alpha_k^{(\mathrm{G})}}\bm{J}_{m_{k,i}}(\lambda_k)\right)\bm{U}^{-1}.
\end{eqnarray}
where $\bm{U} \in \mathbb{C}^{\infty \times \infty}$ is an invertible infinite-dimensional matrix, and $\alpha_k^{(\mathrm{G})}$ is the geometry multiplicity with respect to the $k$-th eigenvalue $\lambda_k$. We have the following relationship about $\alpha_k^{(\mathrm{A})}$ and $\alpha_k^{(\mathrm{G})}$:
\begin{eqnarray}
\sum\limits_{i=1}^{\alpha_k^{(\mathrm{G})}}m_{k,i}&=&\alpha_k^{(\mathrm{A})}.
\end{eqnarray}

We have the following lemma about Jordan decomposition for countable infinite-dimensional operators.
\begin{lemma}\label{lma: projectors and nilpotents of a matrix count-inf}
Given a square infinite dimensional matrix $\bm{X} \in \mathbb{C}^{\infty \times \infty}$, we can express it by
\begin{eqnarray}\label{eq1: lma: projectors and nilpotents of a matrix count-inf}
\bm{X}&=&\sum\limits_{k=1}^\infty\sum\limits_{i=1}^{\alpha_k^{\mathrm{G}}} \lambda_k \bm{P}_{k,i}+
\sum\limits_{k=1}^\infty\sum\limits_{i=1}^{\alpha_k^{\mathrm{G}}} \bm{N}_{k,i},
\end{eqnarray}
where $\bm{P}_{k,i}\in \mathbb{C}^{\infty \times \infty}$ is the projector matrices with respect to the eigenvalue $\lambda_k$ and $i$-th dimensoin of the null space $\mbox{Null}(\bm{X} - \lambda_k \bm{I})$; and $\bm{N}_{k,i}\in \mathbb{C}^{\infty \times \infty}$ is the nilpotent matrices with respect to the eigenvalue $\lambda_k$ and $i$-th dimensoin of the null space $\mbox{Null}(\bm{X} - \lambda_k \bm{I})$.
\end{lemma}
\textbf{Proof:}
Since each Jordan Block with dimension $m_{k,i}\times m_{k,i}$ and the eigenvalue $\lambda_k$, denoted by $\bm{J}_{m_{k,i}}(\lambda_k)$, can be expressed by  
\begin{eqnarray}\label{eq1.1: lma: projectors and nilpotents of a matrix count-inf}
\bm{J}_{m_{k,i}}(\lambda_k)&=&\begin{pmatrix}
   \lambda_k & 1 & 0 & \cdots & 0 \\
   0 & \lambda_k & 1 & \cdots & 0 \\
   \vdots & \vdots & \ddots & \ddots & \vdots \\
   0 & 0 & \cdots & \lambda_k & 1 \\
   0 & 0 & \cdots & 0 & \lambda_k
   \end{pmatrix}_{m_{k,i} \times m_{k,i}} \nonumber \\
&=& \lambda_k \bm{I}_{m_{k,i}}+ \acute{\bm{J}}_{m_{k,i}},
\end{eqnarray}
where $\bm{I}_{m_{k,i}}$ is an identity matrix with the dimension $m_{k,i} \times m_{k,i}$, and $\acute{\bm{J}}_{m_{k,i}}$ is a off-diagonal matrix with the dimension $m_{k,i} \times m_{k,i}$ such that 
\begin{eqnarray}\label{eq2: lma: projectors and nilpotents of a matrix count-inf}
\acute{\bm{J}}_{m_{k,i}}&=&\begin{pmatrix}
   0 & 1 & 0 & \cdots & 0 \\
   0 & 0 & 1 & \cdots & 0 \\
   \vdots & \vdots & \ddots & \ddots & \vdots \\
   0 & 0 & \cdots & 0 & 1 \\
   0 & 0 & \cdots & 0 & 0
   \end{pmatrix}_{m_{k,i} \times m_{k,i}}. 
\end{eqnarray}

We construct the identity matrix $\bm{I}_{k,i}$ with the dimension $\infty \times \infty$ as
\begin{eqnarray}\label{eq3: lma: projectors and nilpotents of a matrix count-inf}
\bm{I}_{k,i}&=&\bigoplus\limits_{k'=1}^{\infty}\bigoplus\limits_{i'=1}^{\alpha_{k'}^{(\mathrm{G})}}\delta(k,k')\delta(i',i)\bm{I}_{m_{k',i'}},
\end{eqnarray}
where $\delta(a,b)$$=$$0$ if $a$$\neq$$b$, and $\delta(a,b)$$=$$1$ if $a$$=$$b$. Similarly, we also construct the off-diagonal matrix $\acute{\bm{J}}$ with the dimension $\infty \times \infty$ as
\begin{eqnarray}\label{eq4: lma: projectors and nilpotents of a matrix count-inf}
\acute{\bm{J}}_{k,i}&=&\bigoplus\limits_{k'=1}^{\infty}\bigoplus\limits_{i'=1}^{\alpha_{k'}^{(\mathrm{G})}}\delta(k,k')\delta(i',i)\acute{\bm{J}}_{m_{k',i'}}.
\end{eqnarray}
Let $\bm{P}_{k,i}$ be constructed by 
\begin{eqnarray}\label{eq5: lma: projectors and nilpotents of a matrix count-inf}
\bm{P}_{k,i}&=&\bm{V}\bm{I}_{k,i}\bm{V}^{-1},
\end{eqnarray}
and 
\begin{eqnarray}\label{eq6: lma: projectors and nilpotents of a matrix count-inf}
\bm{N}_{k,i}&=&\bm{V}\acute{\bm{J}}_{k,i}\bm{V}^{-1}.
\end{eqnarray}
From Eq.~\eqref{eq: Jordan decompostion count-inf}, we have Eq.~\eqref{eq1: lma: projectors and nilpotents of a matrix count-inf}.

Finally, we have to show that matrices $\bm{P}_{k,i}$ are projector matrices and matrices $\bm{N}_{k,i}$ are nilpotent matrices. Because we have
\begin{eqnarray}\label{eq7: lma: projectors and nilpotents of a matrix count-inf}
\bm{P}_{k,i}\bm{P}_{k',i'}=\bm{P}_{k,i}\delta(k,k')\delta(i,i'),
\end{eqnarray}
and
\begin{eqnarray}\label{eq8: lma: projectors and nilpotents of a matrix count-inf}
\sum\limits_{k=1}^\infty\sum\limits_{i=1}^{\alpha_k^{\mathrm{G}}}\bm{P}_{k,i}=\bm{I}.
\end{eqnarray}
Therefore,  matrices $\bm{P}_{k,i}$ are projectors as they are orthogonal and complete. On the other hand, 
\begin{eqnarray}\label{eq9: lma: projectors and nilpotents of a matrix count-inf}
\bm{N}_{k,i}^{\ell}\neq\bm{0}, \mbox{~and~~} \bm{N}_{k,i}^{m_{k,i}} =\bm{0}, 
\end{eqnarray}
where $\ell \in \mathbb{N}$ and $\ell < m_{k,i}$. This shows that matrices $\bm{N}_{k,i}$ are nilpotent matrices.
$\hfill\Box$

\subsection{Analogous Infinite-Dimensional Matrices and Their Properties}\label{sec: Analogous Infinite-Dimensional Matrices and Their Properties}

In this secion, we will extend those results from Section~\ref{sec: Analogous Matrices and Their Properties} to infinite-dimensional matrices.

Given two matrices $\bm{X}$ and $\bm{Y}$ with the same size $\infty \times \infty$ and the same number of distinct eigenvalues $\lambda_k$ for $k=1,2,\ldots,\infty$, we say that the matrix $\bm{X}$ is analogous to the matrix $\bm{Y}$ with respect to the ratios $[c_1, c_2, \ldots] \in \mathbb{C}^\infty$ and $c_i \neq 0$ for $k=1,2,\ldots$, denoted by $\bm{X} \propto^{\bm{U},\bm{V}}_{[c_1, c_2, \ldots]} \bm{Y}$, if these two matrices $\bm{X}$ and $\bm{Y}$ can be expressed by
\begin{eqnarray}\label{eq1: Analogous Operators count-inf}
\bm{X}&=& \bm{U}\left(\bigoplus\limits_{k=1}^{\infty}\bigoplus\limits_{i=1}^{\alpha_k^{(\mathrm{G})}}\bm{J}_{m_{k,i}}(\lambda_k)\right)\bm{U}^{-1}, \nonumber \\
\bm{Y}&=& \bm{V}\left(\bigoplus\limits_{k=1}^{\infty}\bigoplus\limits_{i=1}^{\alpha_k^{(\mathrm{G})}}\bm{J}_{m_{k,i}}(c_k \lambda_k)\right)\bm{V}^{-1}.
\end{eqnarray}

We will have the following properties about two analogous matrices $\bm{X}$ and $\bm{Y}$ given by the relationship $\bm{X} \propto^{\bm{U},\bm{V}}_{[c_1, c_2, \ldots]} \bm{Y}$. 

\begin{proposition}\label{prop0: Analogous Operators count-inf}
% same rank
Two analogous matrices $\bm{X}$ and $\bm{Y}$ given by the relationship $\bm{X} \propto^{\bm{U},\bm{V}}_{[c_1, c_2, \ldots]} \bm{Y}$, then, these two matrices $\bm{X}$ and $\bm{Y}$ have the same rank.
\end{proposition}
\textbf{Proof:}
From Eq.~\eqref{eq1: Analogous Operators count-inf}, we know that the number of zero-valued eigenvalues of the matrix $\bm{X}$ should be equal to the number of zero-valued eigenvalues of the matrix $\bm{Y}$.
$\hfill\Box$

The next proposition is to show that the similar relationship between two infinite-dimensional matrices is the special case of the proposed analogous relationship. 

\begin{proposition}\label{prop1: Analogous Operators count-inf}
Two analogous infinite-dimensional matrices $\bm{X}$ and $\bm{Y}$ given by the relationship $\bm{X} \propto^{\bm{U},\bm{V}}_{[1, 1, \ldots]} \bm{Y}$, then, these two infinite-dimensional matrices $\bm{X}$ and $\bm{Y}$ are similar.
\end{proposition}
\textbf{Proof:}
Since we have 
\begin{eqnarray}\label{eq1: prop1: Analogous Operators count-inf}
\bm{X}&=& \bm{U}\left(\bigoplus\limits_{k=1}^{\infty}\bigoplus\limits_{i=1}^{\alpha_k^{(\mathrm{G})}}\bm{J}_{m_{k,i}}(\lambda_k)\right)\bm{U}^{-1}, \nonumber \\
\bm{Y}&=& \bm{V}\left(\bigoplus\limits_{k=1}^{\infty}\bigoplus\limits_{i=1}^{\alpha_k^{(\mathrm{G})}}\bm{J}_{m_{k,i}}(1 \times \lambda_k)\right)\bm{V}^{-1},
\end{eqnarray}   
therefore, we can obtain the matrix $\bm{Y}$ from the matrix $\bm{X}$ by 
\begin{eqnarray}\label{eq2: prop1: Analogous Operators count-inf}
\bm{V}\bm{U}^{-1}\bm{X}\bm{U}\bm{V}^{-1}&=& \bm{Y}.
\end{eqnarray}   
$\hfill\Box$

The next proposition is about the commutative condition of two analogous infinite-dimensional matrices.

\begin{proposition}\label{prop2: Analogous Operators count-inf}
Two analogous  infinite-dimensional matrices $\bm{X}$ and $\bm{Y}$ are given by the relationship \\
$\bm{X} \propto^{\bm{U},\bm{U}}_{[c_1, c_2, \ldots]} \bm{Y}$, then, these two  infinite-dimensional matrices $\bm{X}$ and $\bm{Y}$ are commute under the conventional matrix multiplication. 
\end{proposition}
\textbf{Proof:}
Becuase two analogous infinite-dimensional matrices $\bm{X}$ and $\bm{Y}$ are given by the relationship \\
$\bm{X} \propto^{\bm{U},\bm{U}}_{[c_1, c_2, \ldots, c_K]} \bm{Y}$, we have
\begin{eqnarray}\label{eq1: prop2: Analogous Operators count-inf}
\bm{X}&=& \bm{U}\left(\bigoplus\limits_{k=1}^{\infty}\bigoplus\limits_{i=1}^{\alpha_k^{(\mathrm{G})}}\bm{J}_{m_{k,i}}(\lambda_k)\right)\bm{U}^{-1}\nonumber \\
&=&\sum\limits_{k=1}^\infty\sum\limits_{i=1}^{\alpha_k^{\mathrm{G}}} \lambda_k \bm{P}_{k,i}+
\sum\limits_{k=1}^\infty\sum\limits_{i=1}^{\alpha_k^{\mathrm{G}}} \bm{N}_{k,i}, 
\end{eqnarray}   
and
\begin{eqnarray}\label{eq2: prop2: Analogous Operators count-inf}
\bm{Y}&=& \bm{U}\left(\bigoplus\limits_{k=1}^{\infty}\bigoplus\limits_{i=1}^{\alpha_k^{(\mathrm{G})}}\bm{J}_{m_{k,i}}(c_k \times \lambda_k)\right)\bm{U}^{-1}\nonumber \\
&=&\sum\limits_{k=1}^\infty\sum\limits_{i=1}^{\alpha_k^{\mathrm{G}}}c_k \lambda_k \bm{P}_{k,i}+
\sum\limits_{k=1}^\infty\sum\limits_{i=1}^{\alpha_k^{\mathrm{G}}} \bm{N}_{k,i}.
\end{eqnarray}   

Besides, we also have the following relations among $\bm{P}_{k,i}$ and $\bm{N}_{k,i}$ from Eq.~\eqref{eq5: lma: projectors and nilpotents of a matrix count-inf} and Eq.~\eqref{eq6: lma: projectors and nilpotents of a matrix count-inf}:
\begin{eqnarray}\label{eq3: prop2: Analogous Operators count-inf}
\bm{P}_{k,i}\bm{P}_{k',i'}&=&\bm{P}_{k,i}\delta(k,k')\delta(i,i'), \nonumber \\
\bm{P}_{k',i'}\bm{N}_{k,i}&=&\bm{N}_{k,i}\bm{P}_{k',i'}=\bm{N}_{k,i}\delta(k,k')\delta(i,i'), \nonumber \\
\bm{N}_{k,i}\bm{N}_{k',i'}&=&\bm{N}^2_{k,i}\delta(k,k')\delta(i,i').
\end{eqnarray}   
Then, we have 
\begin{eqnarray}\label{eq4: prop2: Analogous Operators count-inf}
\bm{X}\bm{Y}&=& \sum\limits_{k=1}^\infty\sum\limits_{i=1}^{\alpha_k^{\mathrm{G}}} c_k \lambda^2_k \bm{P}_{k,i}+ \sum\limits_{k=1}^\infty\sum\limits_{i=1}^{\alpha_k^{\mathrm{G}}} c_k \lambda_k \bm{N}_{k,i} \nonumber \\
&&+ \sum\limits_{k=1}^\infty\sum\limits_{i=1}^{\alpha_k^{\mathrm{G}}} \lambda_k \bm{N}_{k,i}+
\sum\limits_{k=1}^\infty\sum\limits_{i=1}^{\alpha_k^{\mathrm{G}}} \bm{N}^2_{k,i} \nonumber \\
&=&  \bm{Y}\bm{X}
\end{eqnarray}   
$\hfill\Box$

The next proposition is to show that the proposed analogous relationship between $\bm{X}$ and $\bm{Y}$ provided by $\bm{X} \propto^{\bm{U},\bm{V}}_{[c_1, c_2, \ldots]} \bm{Y}$ is an equivalence relation.

\begin{proposition}\label{prop3: Analogous Operators}
The relationship $\bm{X} \propto^{\bm{U},\bm{V}}_{[c_1, c_2, \ldots, c_K]} \bm{Y}$ is an equivalence relation.
\end{proposition}
\textbf{Proof:}
Given the infinite-dimensional matrices $\bm{X}, \bm{Y}$ and $\bm{Z}$ with the following formats:
\begin{eqnarray}\label{eq1: prop3: Analogous Operators count-inf}
\bm{X}&=& \bm{U}\left(\bigoplus\limits_{k=1}^{\infty}\bigoplus\limits_{i=1}^{\alpha_k^{(\mathrm{G})}}\bm{J}_{m_{k,i}}(\lambda_k)\right)\bm{U}^{-1}, \nonumber \\
\bm{Y}&=& \bm{V}\left(\bigoplus\limits_{k=1}^{\infty}\bigoplus\limits_{i=1}^{\alpha_k^{(\mathrm{G})}}\bm{J}_{m_{k,i}}(c_k \lambda_k)\right)\bm{V}^{-1},\nonumber \\
\bm{Z}&=& \bm{W}\left(\bigoplus\limits_{k=1}^{\infty}\bigoplus\limits_{i=1}^{\alpha_k^{(\mathrm{G})}}\bm{J}_{m_{k,i}}(d_k \lambda_k)\right)\bm{W}^{-1},
\end{eqnarray} 
we have 
\begin{eqnarray}\label{eq2s: prop3: Analogous Operators count-inf}
\bm{X}&\propto^{\bm{U},\bm{U}}_{[1, 1, \ldots]}&\bm{X}, ~(\mbox{reflexive})\nonumber \\
\bm{X}&\propto^{\bm{U},\bm{V}}_{[c_1, c_2, \ldots]}&\bm{Y}, \mbox{~and $\bm{Y}\propto^{\bm{V},\bm{U}}_{[1/c_1, 1/c_2, \ldots]}\bm{X}$}~(\mbox{symmetric})\nonumber \\
\bm{X}&\propto^{\bm{U},\bm{V}}_{[c_1, c_2, \ldots]}&\bm{Y}, \mbox{~and $\bm{Y}\propto^{\bm{V},\bm{W}}_{[d_1/c_1, d_2/c_2, \ldots]}\bm{Z}$}\nonumber \\
&\Longrightarrow& \bm{X}\propto^{\bm{U},\bm{W}}_{[d_1, d_2, \ldots]}\bm{Z}~(\mbox{transitive})
\end{eqnarray} 
$\hfill\Box$

The next proposition is to show the matrix determinant relationship between two analogous infinite-dimensional matrices. 
\begin{proposition}\label{prop4: Analogous Operators count-inf}
Given two analogous matrices $\bm{X}$ and $\bm{Y}$ provided by the relationship $\bm{X} \propto^{\bm{U},\bm{V}}_{[c_1, c_2, \ldots]} \bm{Y}$ with
\begin{eqnarray}\label{eq1: prop4: Analogous Operators count-inf}
\bm{X}&=& \bm{U}\left(\bigoplus\limits_{k=1}^{\infty}\bigoplus\limits_{i=1}^{\alpha_k^{(\mathrm{G})}}\bm{J}_{m_{k,i}}(\lambda_k)\right)\bm{U}^{-1}, \nonumber \\
\bm{Y}&=& \bm{V}\left(\bigoplus\limits_{k=1}^{\infty}\bigoplus\limits_{i=1}^{\alpha_k^{(\mathrm{G})}}\bm{J}_{m_{k,i}}(c_k \lambda_k)\right)\bm{V}^{-1};
\end{eqnarray}   
we have 
\begin{eqnarray}\label{eq2: prop4: Analogous Operators count-inf}
\det(\bm{X})&=&\frac{\det(\bm{Y})}{\left(\prod\limits_{k=1}^\infty c_k^{\alpha_k^{(\mathrm{A})}}\right)}. 
\end{eqnarray}
\end{proposition}
\textbf{Proof:}
Because we have
\begin{eqnarray}\label{eq3: prop4: Analogous Operators count-inf}
\det(\bm{X})&=&\prod\limits_{k=1}^\infty \lambda_k^{\alpha_k^{(\mathrm{A})}},
\end{eqnarray}
and 
\begin{eqnarray}\label{eq4: prop4: Analogous Operators count-inf}
\det(\bm{Y})&=&\prod\limits_{k=1}^\infty\left(c_k\lambda_k\right)^{\alpha_k^{(\mathrm{A})}},
\end{eqnarray}
this proposition is proved by comparing Eq.~\eqref{eq3: prop4: Analogous Operators count-inf} and Eq.~\eqref{eq4: prop4: Analogous Operators count-inf}.
$\hfill\Box$

The next proposition is to show the matrix trace relationship for two analogous infinite dimensional matrices.

\begin{proposition}\label{prop5: Analogous Operators count-inf}
Given two analogous infinite dimensional matrices $\bm{X}$ and $\bm{Y}$ provided by the relationship $\bm{X} \propto^{\bm{U},\bm{V}}_{[c, c, \ldots]} \bm{Y}$ with
\begin{eqnarray}\label{eq1: prop5: Analogous Operators count-inf}
\bm{X}&=& \bm{U}\left(\bigoplus\limits_{k=1}^{\infty}\bigoplus\limits_{i=1}^{\alpha_k^{(\mathrm{G})}}\bm{J}_{m_{k,i}}(\lambda_k)\right)\bm{U}^{-1}, \nonumber \\
\bm{Y}&=& \bm{V}\left(\bigoplus\limits_{k=1}^{\infty}\bigoplus\limits_{i=1}^{\alpha_k^{(\mathrm{G})}}\bm{J}_{m_{k,i}}(c \lambda_k)\right)\bm{V}^{-1};
\end{eqnarray}   
we have 
\begin{eqnarray}\label{eq2: prop5: Analogous Operators count-inf}
\mathrm{Trace}(\bm{X})&=&c \mathrm{Trace}(\bm{Y})
\end{eqnarray}
\end{proposition}
\textbf{Proof:}
Because we have
\begin{eqnarray}\label{eq3: prop5: Analogous Operators count-inf}
\mathrm{Trace}(\bm{X})&=&\sum\limits_{k=1}^\infty \alpha_k^{(\mathrm{A})}\lambda_k,
\end{eqnarray}
and 
\begin{eqnarray}\label{eq4: prop5: Analogous Operators count-inf}
\mathrm{Trace}(\bm{Y})&=&\sum\limits_{k=1}^\infty \alpha_k^{(\mathrm{A})}c\lambda_k,
\end{eqnarray}
this proposition is proved by comparing Eq.~\eqref{eq3: prop5: Analogous Operators count-inf} and Eq.~\eqref{eq4: prop5: Analogous Operators count-inf}.
$\hfill\Box$

\subsection{Graph Representation of Projectors and Nilpotents in Infinite-Dimensional Operators}\label{sec: Graph Representation of Projectors and Nilpotents in Infinite-Dimensional Operators}

The graph representations in Figures~\ref{fig:Projector} and~\ref{fig:Nilpotent} can be extended to incorporate infinite projectors and nilpotent elements, providing a structured visualization that facilitates both visual and analytical exploration of their algebraic interactions.

\subsection{Asymptotic Size of Analogous Operator Types}\label{sec: Asymptotic Size of Analogous Operator Types}

In this section, we consider the asymptotic size of analogous operator types when the dimension approaches to infinity. We have the followin theorem that applies the results from Section~\ref{sec: Enumeration of Finite-Dimensional Operator Types}

\begin{theorem}\label{thm: Asymptotic Size of Analogous Operator Types}
Consider finite dimensional matrices with the dimensions $m \times m$ and $K$ (fixed) distinct eigenvalues such that $K \ll m$. If we assume that $\alpha_k^{(A)}=\frac{m}{K}$, asymptotically, the total number of analogous families for such matrices is
\begin{eqnarray}\label{eq0: thm: Asymptotic Size of Analogous Operator Types}
\left[\frac{\exp(\pi \sqrt{2\frac{m}{K}/3})}{4\sqrt{3}\frac{m}{K}}\right]^K,
\end{eqnarray}
where $m \rightarrow \infty$.
\end{theorem}
\textbf{Proof:} 
From Eq.~\eqref{eq8: enum}, the total number of analogous families for such matrices becomes
\begin{eqnarray}\label{eq1: thm: Asymptotic Size of Analogous Operator Types}
\prod\limits_{k=1}^{K}\mathrm{P}\left(\alpha_k^{(A)}\right)&=&\prod\limits_{k=1}^{K}\mathrm{P}\left(\frac{m}{K}\right)\nonumber \\
&\approx_1&\left[\frac{\exp(\pi \sqrt{\frac{2m}{3K}})}{\frac{4\sqrt{3}m}{K}}\right]^K,
\end{eqnarray}
where we apply approximation formula in $\approx_1$ for the partition function $\mathrm{P}\left(\frac{m}{K}\right)$~\cite{brigham1950general}.  
$\hfill\Box$

\subsection{Spectral Mapping Theorem for Multivariate Countable Infinite-Dimensional Operators}\label{sec: Spectral Mapping Theorem for Multivariate Countable Infinite-Dimensional Operators}

In this section, we will establish spectral mapping theorem for multivariate countable infinite-dimensional operators. 

\subsubsection{Single Variable}\label{sec: Single Variable}

\begin{theorem}\label{thm: Spectral Mapping Theorem for Single Variable count-inf}
Given an analytic function $f(z)$ within the domain for $|z| < R$, a square infinite-dimensional matrix $\bm{X}$ with infinite distinct eigenvalues $\lambda_k$ for $k=1,2,\ldots$ such that
\begin{eqnarray}\label{eq1: thm: Spectral Mapping Theorem for Single Variable count-inf}
\bm{X}&=&\sum\limits_{k=1}^\infty\sum\limits_{i=1}^{\alpha_k^{\mathrm{G}}} \lambda_k \bm{P}_{k,i}+
\sum\limits_{k=1}^\infty\sum\limits_{i=1}^{\alpha_k^{\mathrm{G}}} \bm{N}_{k,i},
\end{eqnarray}
where $\left\vert\lambda_k\right\vert<R$, then, we have
\begin{eqnarray}\label{eq2: thm: Spectral Mapping Theorem for Single Variable count-inf}
f(\bm{X})&=&\sum\limits_{k=1}^\infty \left[\sum\limits_{i=1}^{\alpha_k^{(\mathrm{G})}}f(\lambda_k)\bm{P}_{k,i}+\sum\limits_{i=1}^{\alpha_k^{(\mathrm{G})}}\sum\limits_{q=1}^{m_{k,i}-1}\frac{f^{(q)}(\lambda_k)}{q!}\bm{N}_{k,i}^q\right].
\end{eqnarray}
\end{theorem}
\textbf{Proof:}
Recall we have
\begin{eqnarray}\label{eq2-1: thm: Spectral Mapping Theorem for Single Variable count-inf}
\bm{P}_{k,i}\bm{P}_{k',i'}&=&\bm{P}_{k,i}\delta(k,k')\delta(i,i'), \nonumber \\
\bm{P}_{k',i'}\bm{N}_{k,i}&=&\bm{N}_{k,i}\bm{P}_{k',i'}=\bm{N}_{k,i}\delta(k,k')\delta(i,i'), \nonumber \\
\bm{N}_{k,i}\bm{N}_{k',i'}&=&\bm{N}^2_{k,i}\delta(k,k')\delta(i,i').
\end{eqnarray}   

Because $f(z)$ is an analytic function, we have
\begin{eqnarray}\label{eq3: thm: Spectral Mapping Theorem for Single Variable count-inf}
f(\bm{X})&=&\sum\limits_{\ell=0}^{\infty}a_\ell \bm{X}^\ell \nonumber \\
&=&\sum\limits_{\ell=0}^{\infty}a_\ell\left(\sum\limits_{k=1}^\infty\sum\limits_{i=1}^{\alpha_k^{\mathrm{G}}} \lambda_k \bm{P}_{k,i}+
\sum\limits_{k=1}^\infty\sum\limits_{i=1}^{\alpha_k^{\mathrm{G}}} \bm{N}_{k,i}\right)^\ell \nonumber \\
&=_1&\sum\limits_{\ell=0}^{\infty}\sum\limits_{k=1}^\infty\sum\limits_{i=1}^{\alpha_k^{\mathrm{G}}} a_\ell\left(\lambda_k \bm{P}_{k,i}+\bm{N}_{k,i}\right)^\ell  \nonumber \\
&=& \sum\limits_{k=1}^\infty\sum\limits_{i=1}^{\alpha_k^{\mathrm{G}}}\sum\limits_{\ell=0}^{\infty}\sum\limits_{q=0}^{\ell} \frac{a_\ell \ell!}{q! (\ell-q)!}\lambda^{\ell-q}_k\bm{P}^{\ell-q}_{k,i}\bm{N}^q_{k,i}\nonumber \\
&=& \sum\limits_{k=1}^\infty\sum\limits_{i=1}^{\alpha_k^{\mathrm{G}}}\sum\limits_{\ell=0}^{\infty}\left[a_{\ell}\lambda_k^{\ell}\bm{P}_{k,i}+\sum\limits_{q=1}^{\ell} \frac{a_\ell \ell!}{q! (\ell-q)!}\lambda^{\ell-q}_k\bm{P}^{\ell-q}_{k,i}\bm{N}^q_{k,i}\right]\nonumber \\
&=_2&\sum\limits_{k=1}^\infty\left[\sum\limits_{i=1}^{\alpha_k^{(\mathrm{G})}}f(\lambda_k)\bm{P}_{k,i}+\sum\limits_{i=1}^{\alpha_k^{(\mathrm{G})}}\sum\limits_{q=1}^{m_{k,i}-1}\frac{f^{(q)}(\lambda_k)}{q!}\bm{N}_{k,i}^q\right]. 
\end{eqnarray}
where we apply Eq.~\eqref{eq2-1: thm: Spectral Mapping Theorem for Single Variable count-inf} in $=_1$ and $=_2$, and $f^{(q)}(z) = \sum\limits_{\ell=q}^{\infty}\frac{a_\ell \ell!}{(\ell-q)!}z^{\ell-q}$ in $=_2$.
$\hfill \Box$

\subsubsection{Multiple Variables}\label{sec: Multiple Variables count-inf}

We first consider spectral mapping theorem for two input matrices in Theorem~\ref{thm: Spectral Mapping Theorem for Two Variables count-inf}.

\begin{theorem}\label{thm: Spectral Mapping Theorem for Two Variables count-inf}
Given an analytic function $f(z_1,z_2)$ within the domain for $|z_1| < R_1$ and $|z_2| < R_2$, the first matrix $\bm{X}_1$ with distinct eigenvalues $\lambda_{k_1}$ for $k_1=1,2,\ldots$ such that
\begin{eqnarray}\label{eq1-1: thm: Spectral Mapping Theorem for Two Variables count-inf}
\bm{X}_1&=&\sum\limits_{k_1=1}^{\infty}\sum\limits_{i_1=1}^{\alpha_{k_1}^{\mathrm{G}}} \lambda_{k_1} \bm{P}_{k_1,i_1}+
\sum\limits_{k_1=1}^{\infty}\sum\limits_{i_1=1}^{\alpha_{k_1}^{\mathrm{G}}} \bm{N}_{k_1,i_1},
\end{eqnarray}
where $\left\vert\lambda_{k_1}\right\vert<R_1$, and second matrix $\bm{X}_2$ with distinct eigenvalues $\lambda_{k_2}$ for $k_2=1,2,\ldots$ such that
\begin{eqnarray}\label{eq1-2: thm: Spectral Mapping Theorem for Two Variables count-inf}
\bm{X}_2&=&\sum\limits_{k_2=1}^{\infty}\sum\limits_{i_2=1}^{\alpha_{k_2}^{\mathrm{G}}} \lambda_{k_2} \bm{P}_{k_2,i_2}+
\sum\limits_{k_2=1}^{\infty}\sum\limits_{i_2=1}^{\alpha_{k_2}^{\mathrm{G}}} \bm{N}_{k_2,i_2},
\end{eqnarray}
where $\left\vert\lambda_{k_2}\right\vert<R_2$.

Then, we have
\begin{eqnarray}\label{eq2: thm: Spectral Mapping Theorem for Two Variables count-inf}
f(\bm{X}_1, \bm{X}_2)&=&\sum\limits_{k_1=1}^{\infty}\sum\limits_{k_2=1}^{\infty}\sum\limits_{i_1=1}^{\alpha_{k_1}^{(\mathrm{G})}}\sum\limits_{i_2=1}^{\alpha_{k_2}^{(\mathrm{G})}}f(\lambda_{k_1}, \lambda_{k_2})\bm{P}_{k_1,i_1}\bm{P}_{k_2,i_2} \nonumber \\
&&+\sum\limits_{k_1=1}^{\infty}\sum\limits_{k_2=1}^{\infty}\sum\limits_{i_1=1}^{\alpha_{k_1}^{(\mathrm{G})}}\sum\limits_{i_2=1}^{\alpha_{k_2}^{(\mathrm{G})}}\sum_{q_2=1}^{m_{k_2,i_2}-1}\frac{f^{(-,q_2)}(\lambda_{k_1},\lambda_{k_2})}{q_2!}\bm{P}_{k_1,i_1}\bm{N}_{k_2,i_2}^{q_2} \nonumber \\
&&+\sum\limits_{k_1=1}^{\infty}\sum\limits_{k_2=1}^{\infty}\sum\limits_{i_1=1}^{\alpha_{k_1}^{(\mathrm{G})}}\sum\limits_{i_2=1}^{\alpha_{k_2}^{(\mathrm{G})}}\sum_{q_1=1}^{m_{k_1,i_1}-1}\frac{f^{(q_1,-)}(\lambda_{k_1},\lambda_{k_2})}{q_1!}\bm{N}_{k_1,i_1}^{q_1}\bm{P}_{k_2,i_2}  \nonumber \\
&&+\sum\limits_{k_1=1}^{\infty}\sum\limits_{k_2=1}^{\infty}\sum\limits_{i_1=1}^{\alpha_{k_1}^{(\mathrm{G})}}\sum\limits_{i_2=1}^{\alpha_{k_2}^{(\mathrm{G})}}\sum_{q_1=1}^{m_{k_1,i_1}-1}\sum_{q_2=1}^{m_{k_2,i_2}-1}\frac{f^{(q_1,q_2)}(\lambda_{k_1},\lambda_{k_2})}{q_1!q_2!}\bm{N}_{k_1,i_1}^{q_1}\bm{N}_{k_2,i_2}^{q_2}
\end{eqnarray}
\end{theorem}
\textbf{Proof:}
Recall we still have
\begin{eqnarray}\label{eq2-1: thm: Spectral Mapping Theorem for Two Variables count-inf}
\bm{P}_{k_1,i_1}\bm{P}_{k'_1,i'_1}&=&\bm{P}_{k_1,i_1}\delta(k_1,k'_1)\delta(i_1,i'_1), \nonumber \\
\bm{P}_{k'_1,i'_1}\bm{N}_{k_1,i_1}&=&\bm{N}_{k_1,i_1}\bm{P}_{k'_1,i'_1}=\bm{N}_{k_1,i_1}\delta(k_1,k'_1)\delta(i_1,i'_1), \nonumber \\
\bm{N}_{k_1,i_1}\bm{N}_{k'_1,i'_1}&=&\bm{N}^2_{k_1,i_1}\delta(k_1,k'_1)\delta(i_1,i'_1).
\end{eqnarray}   
\begin{eqnarray}\label{eq2-2: thm: Spectral Mapping Theorem for Two Variables count-inf}
\bm{P}_{k_2,i_2}\bm{P}_{k'_2,i'_2}&=&\bm{P}_{k_2,i_2}\delta(k_2,k'_2)\delta(i_2,i'_2), \nonumber \\
\bm{P}_{k'_2,i'_2}\bm{N}_{k_2,i_2}&=&\bm{N}_{k_2,i_2}\bm{P}_{k'_2,i'_2}=\bm{N}_{k_2,i_2}\delta(k_2,k'_2)\delta(i_2,i'_2), \nonumber \\
\bm{N}_{k_2,i_2}\bm{N}_{k'_2,i'_2}&=&\bm{N}^2_{k_2,i_2}\delta(k_2,k'_2)\delta(i_2,i'_2).
\end{eqnarray}   

Because $f(z)$ is an analytic function, we have
\begin{eqnarray}\label{eq3: thm: Spectral Mapping Theorem for Two Variables count-inf}
f(\bm{X}_1,\bm{X}_2)&=&\sum\limits_{\ell_1=0,\ell_2=0}^{\infty}a_{\ell_1,\ell_2} \bm{X}_1^{\ell_1}\bm{X}_2^{\ell_2}\nonumber \\
&=&\sum\limits_{\ell_1=0,\ell_2=0}^{\infty}a_{\ell_1,\ell_2}\left(\sum\limits_{k_1=1}^{\infty}\sum\limits_{i_1=1}^{\alpha_{k_1}^{\mathrm{G}}} \lambda_{k_1} \bm{P}_{k_1,i_1}+
\sum\limits_{k_1=1}^{\infty}\sum\limits_{i_1=1}^{\alpha_{k_1}^{\mathrm{G}}} \bm{N}_{k_1,i_1}\right)^{\ell_1}\nonumber \\
&&\times \left(\sum\limits_{k_2=1}^{\infty}\sum\limits_{i_2=1}^{\alpha_{k_2}^{\mathrm{G}}} \lambda_{k_2} \bm{P}_{k_2,i_2}+
\sum\limits_{k_2=1}^{\infty}\sum\limits_{i_2=1}^{\alpha_{k_2}^{\mathrm{G}}} \bm{N}_{k_2,i_2}\right)^{\ell_2}\nonumber \\
&=_1&\sum\limits_{\ell_1=0,\ell_2=0}^{\infty}a_{\ell_1,\ell_2}\left[\sum\limits_{k_1=1}^\infty\sum\limits_{i_1=1}^{\alpha_{k_1}^{\mathrm{G}}}\left(\lambda_{k_1} \bm{P}_{k_1,i_1}+\bm{N}_{k_1,i_1}\right)^{\ell_1}\right]\left[\sum\limits_{k_2=1}^\infty\sum\limits_{i_2=1}^{\alpha_{k_2}^{\mathrm{G}}}\left(\lambda_{k_2} \bm{P}_{k_2,i_2}+\bm{N}_{k_2,i_2}\right)^{\ell_2}\right]\nonumber \\
&=& \sum\limits_{k_1=1}^{\infty}\sum\limits_{k_2=1}^{\infty}\sum\limits_{i_1=1}^{\alpha_{k_1}^{(\mathrm{G})}}\sum\limits_{i_2=1}^{\alpha_{k_2}^{(\mathrm{G})}}\sum\limits_{\ell_1=0,\ell_2=0}^{\infty}a_{\ell_1,\ell_2}\left[\sum\limits_{q_1=0}^{m_{k_1,i_1}-1}\frac{\ell_1!}{q_1! (\ell_1-q_1)!}\lambda^{\ell_1-q_1}_{k_1}\bm{P}^{\ell_1-q_1}_{k_1,i_1}\bm{N}^{q_1}_{k_1,i_1}\right]\nonumber \\
&&\times \left[\sum\limits_{q_2=0}^{m_{k_2,i_2}-1}\frac{\ell_2!}{q_2! (\ell_2-q_2)!}\lambda^{\ell_2-q_2}_{k_2}\bm{P}^{\ell_2-q_2}_{k_2,i_2}\bm{N}^{q_2}_{k_2,i_2}\right]\nonumber \\
&=& \sum\limits_{k_1=1}^{\infty}\sum\limits_{k_2=1}^{\infty}\sum\limits_{i_1=1}^{\alpha_{k_1}^{(\mathrm{G})}}\sum\limits_{i_2=1}^{\alpha_{k_2}^{(\mathrm{G})}}\sum\limits_{\ell_1=0,\ell_2=0}^{\infty}a_{\ell_1,\ell_2}\left[\lambda_{k_1}^{\ell_1}\bm{P}^{\ell_1}_{k_1,i_1}+\sum\limits_{q_1=1}^{m_{k_1,i_1}-1}\frac{\ell_1!}{q_1! (\ell_1-q_1)!}\lambda^{\ell_1-q_1}_{k_1}\bm{P}^{\ell_1-q_1}_{k_1,i_1}\bm{N}^{q_1}_{k_1,i_1}\right]\nonumber \\
&&\times \left[\lambda_{k_2}^{\ell_2}\bm{P}^{\ell_2}_{k_2,i_2}+\sum\limits_{q_2=1}^{m_{k_2,i_2}-1}\frac{\ell_2!}{q_2! (\ell_2-q_2)!}\lambda^{\ell_2-q_2}_{k_2}\bm{P}^{\ell_2-q_2}_{k_2,i_2}\bm{N}^{q_2}_{k_2,i_2}\right]\nonumber \\
&=_2&\sum\limits_{k_1=1}^{\infty}\sum\limits_{k_2=1}^{\infty}\sum\limits_{i_1=1}^{\alpha_{k_1}^{(\mathrm{G})}}\sum\limits_{i_2=1}^{\alpha_{k_2}^{(\mathrm{G})}}f(\lambda_{k_1}, \lambda_{k_2})\bm{P}_{k_1,i_1}\bm{P}_{k_2,i_2} \nonumber \\
&&+\sum\limits_{k_1=1}^{\infty}\sum\limits_{k_2=1}^{\infty}\sum\limits_{i_1=1}^{\alpha_{k_1}^{(\mathrm{G})}}\sum\limits_{i_2=1}^{\alpha_{k_2}^{(\mathrm{G})}}\sum_{q_2=1}^{m_{k_2,i_2}-1}\frac{f^{(-,q_2)}(\lambda_{k_1},\lambda_{k_2})}{q_2!}\bm{P}_{k_1,i_1}\bm{N}_{k_2,i_2}^{q_2} \nonumber \\
&&+\sum\limits_{k_1=1}^{\infty}\sum\limits_{k_2=1}^{\infty}\sum\limits_{i_1=1}^{\alpha_{k_1}^{(\mathrm{G})}}\sum\limits_{i_2=1}^{\alpha_{k_2}^{(\mathrm{G})}}\sum_{q_1=1}^{m_{k_1,i_1}-1}\frac{f^{(q_1,-)}(\lambda_{k_1},\lambda_{k_2})}{q_1!}\bm{N}_{k_1,i_1}^{q_1}\bm{P}_{k_2,i_2}  \nonumber \\
&&+\sum\limits_{k_1=1}^{\infty}\sum\limits_{k_2=1}^{\infty}\sum\limits_{i_1=1}^{\alpha_{k_1}^{(\mathrm{G})}}\sum\limits_{i_2=1}^{\alpha_{k_2}^{(\mathrm{G})}}\sum_{q_1=1}^{m_{k_1,i_1}-1}\sum_{q_2=1}^{m_{k_2,i_2}-1}\frac{f^{(q_1,q_2)}(\lambda_{k_1},\lambda_{k_2})}{q_1!q_2!}\bm{N}_{k_1,i_1}^{q_1}\bm{N}_{k_2,i_2}^{q_2}. 
\end{eqnarray}
where we apply Eq.~\eqref{eq2-1: thm: Spectral Mapping Theorem for Two Variables count-inf} and Eq.~\eqref{eq2-2: thm: Spectral Mapping Theorem for Two Variables count-inf} in $=_1$ and $=_2$, and following relations about the partial derivatives of $f(z_1,z_2)$ in $=_2$:
\begin{eqnarray}\label{eq4: thm: Spectral Mapping Theorem for Two Variables count-inf}
f^{(-,q_2)}(z_1,z_2)&=&\sum_{\ell_1=0}^{\infty}\sum_{\ell_2=q_2}^{\infty}\frac{a_{\ell_1,\ell_2}\ell_2!}{(\ell_2-q_2)!}z_1^{\ell_1}z_2^{\ell_2-q_2},\nonumber\\
f^{(q_1,-)}(z_1,z_2)&=&\sum_{\ell_1=q_1}^{\infty}\sum_{\ell_2=0}^{\infty}\frac{a_{\ell_1,\ell_2}\ell_1!}{(\ell_1-q_1)!}z_1^{\ell_1-q_1}z_2^{\ell_2},\nonumber\\
f^{(q_1,q_2)}(z_1,z_2)&=&\sum_{\ell_1=q_1}^{\infty}\sum_{\ell_2=q_2}^{\infty}\frac{a_{\ell_1,\ell_2}\ell_1!\ell_2!}{(\ell_1-q_1)!(\ell_2-q_2)!}z_1^{\ell_1-q_1}z_2^{\ell_2-q_2}.
\end{eqnarray}
$\hfill \Box$

We are ready to present Theorem~\ref{thm: Spectral Mapping Theorem for r Variables count-inf} about spectral mapping theorem for $r$ infinite dimensional matrices. 
\begin{theorem}\label{thm: Spectral Mapping Theorem for r Variables count-inf}
Given an analytic function $f(z_1,z_2,\ldots,z_r)$ within the domain for $|z_l| < R_l$, and the matrix $\bm{X}_l$ with distinct eigenvalues $\lambda_{k_l}$ for $k_l=1,2,\ldots$ such that
\begin{eqnarray}\label{eq1-1: thm: Spectral Mapping Theorem for r Variables count-inf}
\bm{X}_l&=&\sum\limits_{k_l=1}^{\infty}\sum\limits_{i_l=1}^{\alpha_{k_l}^{\mathrm{G}}} \lambda_{k_l} \bm{P}_{k_l,i_l}+
\sum\limits_{k_l=1}^{\infty}\sum\limits_{i_l=1}^{\alpha_{k_l}^{\mathrm{G}}} \bm{N}_{k_l,i_l},
\end{eqnarray}
where $\left\vert\lambda_{k_l}\right\vert<R_l$ for $l=1,2,\ldots,r$.

Then, we have
\begin{eqnarray}\label{eq2: thm: Spectral Mapping Theorem for kappa Variables count-inf}
\lefteqn{f(\bm{X}_1,\ldots,\bm{X}_r)=}\nonumber \\
&& \sum\limits_{k_1=\ldots=k_r=1}^{\infty,\ldots,\infty} \sum\limits_{i_1=\ldots=i_r=1}^{\alpha_{k_1}^{(\mathrm{G})},\ldots,\alpha_{k_r}^{(\mathrm{G})}}
f(\lambda_{k_1},\ldots,\lambda_{k_r})\bm{P}_{k_1,i_1}\ldots\bm{P}_{k_r,i_r}\nonumber \\
&&+\sum\limits_{k_1=\ldots=k_r=1}^{\infty,\ldots,\infty} \sum\limits_{i_1=\ldots=i_r=1}^{\alpha_{k_1}^{(\mathrm{G})},\ldots,\alpha_{k_r}^{(\mathrm{G})}}\sum\limits_{\kappa=1}^{r-1}\sum\limits_{\alpha_\kappa(q_1,\ldots,q_r)}\Bigg(\sum\limits_{\alpha_{\kappa}(q_1,\ldots,q_r)=1}^{m_{k_{\mbox{Ind}(\alpha_{\kappa}(q_1,\ldots,q_r))},i_{\mbox{Ind}(\alpha_{\kappa}(q_1,\ldots,q_r))}}-1}\nonumber \\
&&~~~~~ \frac{f^{\alpha_{\kappa}(q_1,\ldots,q_r)}(\lambda_{k_1},\ldots,\lambda_{k_r})}{q_{\iota_1}!q_{\iota_2}!\ldots q_{\iota_\kappa}!}\times \prod\limits_{\substack{\beta =\mbox{Ind}(\alpha_{\kappa}(q_1,\ldots,q_r)), \bm{Y}=\bm{N}^{q_\beta}_{k_\beta,i_\beta} \\ \beta \neq \mbox{Ind}(\alpha_{\kappa}(q_1,\ldots,q_r)), \bm{Y}=\bm{P}_{k_\beta,i_\beta} }
}^{r} \bm{Y}\Bigg) 
\nonumber \\
&&+\sum\limits_{k_1=\ldots=k_r=1}^{\infty,\ldots,\infty} \sum\limits_{i_1=\ldots=i_r=1}^{\alpha_{k_1}^{(\mathrm{G})},\ldots,\alpha_{k_r}^{(\mathrm{G})}} \sum\limits_{q_1=\ldots=q_r=1}^{m_{k_1,i_1}-1,\ldots,m_{k_r,i_r}-1}
\frac{f^{(q_1,\ldots,q_r)}(\lambda_{k_1},\ldots,\lambda_{k_r})}{q_1!\cdots q_r!}\bm{N}^{q_1}_{k_1,i_1}\ldots\bm{N}^{q_r}_{k_r,i_r}.
\end{eqnarray}
\end{theorem}
\textbf{Proof:}
The proof follows directly from Theorem~\ref{thm: Spectral Mapping Theorem for Two Variables count-inf}, with the remaining steps consisting primarily of routine notational manipulations.
$\hfill\Box$

\section{Continous Spectrum Operators}\label{sec: Continous Spectrum Operators}

\subsection{Spectral Mapping Theorem for Continous Spectrum Operators}\label{sec: Spectral Mapping Theorem for Continous Spectrum Operators}

In this section, spectral mapping theorem for continous spectrum operators and its functional calculus will be presented. We will focus in discussing spectral operators~\cite{dunford1958survey}. We first introduce the following Theorem~\ref{thm: spectral theorem of inf-dim} about the spectral decomposition.

\begin{theorem}\label{thm: spectral theorem of inf-dim}
Given a spectral operator $\bm{X}$, we can express it by
\begin{eqnarray}
\bm{X}&=&\int\limits_{\lambda \in \sigma(\bm{X})}\lambda d\bm{E}_{\bm{X}}(\lambda)+
\int\limits_{\lambda \in \sigma(\bm{X})}\left(\bm{X}-\lambda\bm{I}\right)d\bm{E}_{\bm{X}}(\lambda),
\end{eqnarray}
where $\sigma(\bm{X})$ is the spectrum of the operator $\bm{X}$, i.e., the eigenvalue set of the operaor $\bm{X}$.
\end{theorem}
\textbf{Proof:}
From~\cite{dunford1958survey}, we have the spectrum of $\bm{X}$ identical to the spectrum of $\int\limits_{\lambda \in \sigma(\bm{X})}\lambda d\bm{E}_{\bm{X}}(\lambda)$, then, we have the following decomposition
\begin{eqnarray}\label{eq1: thm: spectral theorem of inf-dim}
\bm{X}&=&\int\limits_{\lambda \in \sigma(\bm{X})}\lambda d\bm{E}_{\bm{X}}(\lambda)+\bm{N}_{\bm{X}},
\end{eqnarray}
where $\bm{E}_{\bm{X}}(\lambda)$ is the spectral measure with respect to the eigenvalue $\lambda$, and $\bm{N}_{\bm{X}}$ is the nilpotent part of the operator $\bm{X}$. 

Therefore, we have
\begin{eqnarray}\label{eq2: thm: spectral theorem of inf-dim}
\bm{N}_{\bm{X}}&=&\bm{X}-\int\limits_{\lambda \in \sigma(\bm{X})}\lambda d\bm{E}_{\bm{X}}(\lambda)\nonumber \\
&=_1&\bm{X}\int\limits_{\lambda \in \sigma(\bm{X})}d\bm{E}_{\bm{X}}(\lambda)-\int\limits_{\lambda \in \sigma(\bm{X})}\lambda\bm{I}d\bm{E}_{\bm{X}}(\lambda)\nonumber \\
&=&\int\limits_{\lambda \in \sigma(\bm{X})}\left(\bm{X}-\lambda\bm{I}\right)d\bm{E}_{\bm{X}}(\lambda),
\end{eqnarray}
where we apply the property $\bm{I}=\int\limits_{\lambda \in \sigma(\bm{X})}\bm{E}_{\bm{X}}(\lambda)$ in $=_1$. This thoerem is proved. 
$\hfill\Box$

Lemma~\ref{lma: orthogonal of spectral basis} below will give relationships among $d\bm{E}_{\bm{X}}$ and $(\bm{X}-\lambda\bm{I})d\bm{E}_{\bm{X}}(\lambda)$ which will be used to establish the functional calculus for the operator $\bm{X}$.

\begin{lemma}\label{lma: orthogonal of spectral basis}
We have the following relationships for a spectral operator $\bm{X}$:
\begin{eqnarray}\label{eq1: lma: orthogonal of spectral basis}
d\bm{E}_{\bm{X}}(\lambda)d\bm{E}_{\bm{X}}(\lambda')&=&d\bm{E}_{\bm{X}}(\lambda)\delta(\lambda,\lambda'),
\end{eqnarray}
\begin{eqnarray}\label{eq1: lma: orthogonal of spectral basis}
(\bm{X}-\lambda\bm{I})d\bm{E}_{\bm{X}}(\lambda)(\bm{X}-\lambda'\bm{I})d\bm{E}_{\bm{X}}(\lambda')=(\bm{X}-\lambda\bm{I})^2d\bm{E}_{\bm{X}}(\lambda)\delta(\lambda,\lambda'),
\end{eqnarray} 
\begin{eqnarray}\label{eq3: lma: orthogonal of spectral basis}
d\bm{E}_{\bm{X}}(\lambda')(\bm{X}-\lambda\bm{I})d\bm{E}_{\bm{X}}(\lambda)&=&(\bm{X}-\lambda\bm{I})d\bm{E}_{\bm{X}}(\lambda)d\bm{E}_{\bm{X}}(\lambda') \nonumber \\
&=&(\bm{X}-\lambda\bm{I})d\bm{E}_{\bm{X}}(\lambda)\delta(\lambda,\lambda').
\end{eqnarray} 
\end{lemma}
\textbf{Proof:}
Let the spectrum of the operator $\bm{X}$ be decomposed by the following disjoint Borel sets:
\begin{eqnarray}\label{eq4: lma: orthogonal of spectral basis}
\sigma(\bm{X})&=&\bigcup\limits_{k=1}^K \sigma_k
\end{eqnarray} 
We define $\bm{N}_{\bm{X}}(\sigma_k)$ as
\begin{eqnarray}\label{eq5: lma: orthogonal of spectral basis}
\bm{N}_{\bm{X}}(\sigma_k)&\define&\int\limits_{\lambda \in \sigma_k}(\bm{X}-\lambda\bm{I})d\bm{E}_{\bm{X}}(\lambda).
\end{eqnarray} 
Also, we also define $\bm{E}_{\bm{X}}(\sigma_k)$ as
\begin{eqnarray}\label{eq5.5: lma: orthogonal of spectral basis}
\bm{E}_{\bm{X}}(\sigma_k)&\define&\int\limits_{\lambda \in \sigma_k}d\bm{E}_{\bm{X}}(\lambda).
\end{eqnarray} 
From the spectral measure properties with respect to $\sigma_k$, we have
\begin{eqnarray}\label{eq6: lma: orthogonal of spectral basis}
\bm{E}_{\bm{X}}(\sigma_k)\bm{E}_{\bm{X}}(\sigma_{k'})&=&\bm{E}_{\bm{X}}(\sigma_k)\delta(k,k'),\nonumber \\
\bm{E}_{\bm{X}}(\sigma_k)\bm{N}_{\bm{X}}(\sigma_{k'})&=&\bm{N}_{\bm{X}}(\sigma_{k'})\bm{E}_{\bm{X}}(\sigma_k)\nonumber\\
&=&\bm{N}_{\bm{X}}(\sigma_{k'})\delta(k,k'),\nonumber \\
\bm{N}_{\bm{X}}(\sigma_k)\bm{N}_{\bm{X}}(\sigma_{k'})&=&\bm{N}^2_{\bm{X}}(\sigma_k)\delta(k,k'),
\end{eqnarray} 
This lemma is proved by applying a limiting process as $K \to \infty$.
$\hfill\Box$

Theorem~\ref{thm: functional cal of single inf-dim} below is the main result of this section. Our functional calculus formula for spectral operators provides a more detailed elaboration compared to the functional calculus formulas given by~\cite{dunford1958survey,kantorovitz1965jordan}, as we offer a more precise characterization of the degree of the nilpotent part corresponding to different eigenvalues.

\begin{theorem}\label{thm: functional cal of single inf-dim}
Given an analytic function $f(z)$ within the domain for $|z| < R$, a spectral operator $\bm{X}$ with the following spectral decomposition
\begin{eqnarray}\label{eq1: thm: functional cal of single inf-dim}
\bm{X}&=&\int\limits_{\lambda \in \sigma(\bm{X})}\lambda d\bm{E}_{\bm{X}}(\lambda)+
\int\limits_{\lambda \in \sigma(\bm{X})}\left(\bm{X}-\lambda\bm{I}\right)d\bm{E}_{\bm{X}}(\lambda),
\end{eqnarray}
where $|\lambda|<R$ for $\lambda \in \sigma(\bm{X})$, then, we have
\begin{eqnarray}\label{eq2: thm: functional cal of single inf-dim}
f(\bm{X})&=&\int\limits_{\lambda \in \sigma(\bm{X})}f(\lambda)d\bm{E}_{\bm{X}}(\lambda)+\int\limits_{\lambda \in \sigma(\bm{X})}\sum\limits_{q=1}^{m_{\lambda}-1}\frac{f^{(q)}(\lambda)}{q!}(\bm{X}-\lambda\bm{I})^q d\bm{E}_{\bm{X}}(\lambda),
\end{eqnarray}
where $m_{\lambda}$ is a positive integer to measure the degree of nilpotency with respect to the eigenvalue $\lambda$, i.e., 
\begin{eqnarray}\label{eq2.5: thm: functional cal of single inf-dim}
(\bm{X}-\lambda\bm{I})^q d\bm{E}_{\bm{X}}(\lambda)&& 
    \begin{cases}
      \neq \bm{0} & \text{if $q < m_{\lambda}$.}\\
      =\bm{0} & \text{if $q \geq m_{\lambda}$.}
    \end{cases}  
\end{eqnarray}
\end{theorem}
\textbf{Proof:}

Because $f(z)$ is an analytic function, we have
\begin{eqnarray}\label{eq3: thm: Spectral Mapping Theorem for Single Variable}
f(\bm{X})&=&\sum\limits_{\ell=0}^{\infty}a_\ell \bm{X}^\ell \nonumber \\
&=&\sum\limits_{\ell=0}^{\infty}a_\ell\left(\int\limits_{\lambda \in \sigma(\bm{X})}\lambda d\bm{E}_{\bm{X}}(\lambda)+
\int\limits_{\lambda \in \sigma(\bm{X})}\left(\bm{X}-\lambda\bm{I}\right)d\bm{E}_{\bm{X}}(\lambda)\right)^\ell \nonumber \\
&=_1&\sum\limits_{\ell=0}^{\infty}\int\limits_{\lambda \in \sigma(\bm{X})}a_\ell\left(\lambda d\bm{E}_{\bm{X}}(\lambda)+\left(\bm{X}-\lambda\bm{I}\right)d\bm{E}_{\bm{X}}(\lambda)\right)^\ell  \nonumber \\
&=&\int\limits_{\lambda \in \sigma(\bm{X})}\sum\limits_{\ell=0}^{\infty}\sum\limits_{q=0}^{\ell} \frac{a_\ell \ell!}{q! (\ell-q)!}\lambda^{\ell-q}\left[d\bm{E}_{\bm{X}}(\lambda)\right]^{\ell-q}\left[(\bm{X}-\lambda\bm{I})d\bm{E}_{\bm{X}}(\lambda)\right]^q\nonumber \\
&=&\int\limits_{\lambda \in \sigma(\bm{X})}\sum\limits_{\ell=0}^{\infty}\left[a_{\ell}\lambda^{\ell}\left[d\bm{E}_{\bm{X}}(\lambda)\right]^\ell+\sum\limits_{q=1}^{\ell}\frac{a_\ell \ell!}{q! (\ell-q)!}\lambda^{\ell-q}\left[d\bm{E}_{\bm{X}}(\lambda)\right]^{\ell-q}\left[(\bm{X}-\lambda\bm{I})d\bm{E}_{\bm{X}}(\lambda)\right]^q\right]\nonumber \\
&=_2&\int\limits_{\lambda \in \sigma(\bm{X})}\left[f(\lambda)d\bm{E}_{\bm{X}}(\lambda)+\sum\limits_{q=1}^{m_{\lambda}-1}\frac{f^{(q)}(\lambda)}{q!}(\bm{X}-\lambda\bm{I})^q d\bm{E}_{\bm{X}}(\lambda)\right]\nonumber \\
&=&\int\limits_{\lambda \in \sigma(\bm{X})}f(\lambda)d\bm{E}_{\bm{X}}(\lambda)+\int\limits_{\lambda \in \sigma(\bm{X})}\sum\limits_{q=1}^{m_{\lambda}-1}\frac{f^{(q)}(\lambda)}{q!}(\bm{X}-\lambda\bm{I})^q d\bm{E}_{\bm{X}}(\lambda),
\end{eqnarray}
where we apply Lemma~\ref{lma: orthogonal of spectral basis} in $=_1$ and $=_2$, and $f^{(q)}(z) = \sum\limits_{\ell=q}^{\infty}\frac{a_\ell \ell!}{(\ell-q)!}z^{\ell-q}$ in $=_2$.
$\hfill \Box$

\subsection{Analogous Continous Spectrum Operators and Their Properties}\label{sec: Analogous Continous Spectrum Operators and Their Properties}

In Section~\ref{sec: Analogous Operators}, we explore concepts related to analogous matrices and their properties. In parallel, we introduce \emph{analogous operators} based on Theorem~\ref{thm: spectral theorem of inf-dim}.

Given two spectral operators $\bm{X}$ and $\bm{Y}$, we say that $\bm{X}$ is analogous to $\bm{Y}$ with respect to the invertible operator $\bm{U}$ and a ratio function $c(\lambda)$, where $c(\lambda) \neq 0$ for $\lambda \in \sigma(\bm{X})$. This relationship, denoted by $\bm{X} \propto^{\bm{U}}_{c(\lambda)}\bm{Y}$, holds if $\bm{X}$ and $\bm{Y}$ can be expressed as following:
\begin{eqnarray}\label{eq1: Analogous Operators inf}
\bm{X}&=&\left(\int\limits_{\lambda \in \sigma(\bm{X})}\lambda d\bm{E}_{\bm{X}}(\lambda)+
\int\limits_{\lambda \in \sigma(\bm{X})}\left(\bm{X}-\lambda\bm{I}\right)d\bm{E}_{\bm{X}}(\lambda)\right), \nonumber \\
\bm{Y}&=& \bm{U}\left(\int\limits_{\lambda \in \sigma(\bm{X})}c(\lambda)\lambda d\bm{E}_{\bm{X}}(\lambda)+
\int\limits_{\lambda \in \sigma(\bm{X})}\left(\bm{X}-c(\lambda)\lambda\bm{I}\right)d\bm{E}_{\bm{X}}(\lambda)\right)\bm{U}^{-1}.
\end{eqnarray}

We will have the following properties about two analogous operators $\bm{X}$ and $\bm{Y}$ given by the relationship $\bm{X} \propto^{\bm{U}}_{c(\lambda)}\bm{Y}$. 

\iffalse
\begin{proposition}\label{prop0: Analogous Operators inf}
% 
Two analogous operators $\bm{X}$ and $\bm{Y}$ given by the relationship $\bm{X} \propto^{\bm{U}}_{c(\lambda)=1} \bm{Y}$, then, these two operators $\bm{X}$ and $\bm{Y}$ have the same number of $m_{\lambda}$.
%
\end{proposition}
%
\textbf{Proof:}
%
From Eq.~\eqref{eq1: Analogous Operators inf}, we know that the number of zero-valued eigenvalues of the matrix $\bm{X}$ should be equal to the number of zero-valued eigenvalues of the matrix $\bm{Y}$.
%
$\hfill\Box$
\fi

Proposition~\ref{prop1: Analogous Operators inf} is to show that the similarity relationship between two operators is the special case of the proposed analogous relationship. 

\begin{proposition}\label{prop1: Analogous Operators inf}
Two analogous operators $\bm{X}$ and $\bm{Y}$ given by the relationship $\bm{X} \propto^{\bm{U}}_{c(\lambda)=1}\bm{Y}$, then, these two operators $\bm{X}$ and $\bm{Y}$ are similar.
\end{proposition}
\textbf{Proof:}
Since we have 
\begin{eqnarray}\label{eq1: prop1: Analogous Operators inf}
\bm{X}&=&\left(\int\limits_{\lambda \in \sigma(\bm{X})}\lambda d\bm{E}_{\bm{X}}(\lambda)+
\int\limits_{\lambda \in \sigma(\bm{X})}\left(\bm{X}-\lambda\bm{I}\right)d\bm{E}_{\bm{X}}(\lambda)\right), \nonumber \\
\bm{Y}&=& \bm{U}\left(\int\limits_{\lambda \in \sigma(\bm{X})}1 \times \lambda d\bm{E}_{\bm{X}}(\lambda)+
\int\limits_{\lambda \in \sigma(\bm{X})}\left(\bm{X}-1 \times \lambda\bm{I}\right)d\bm{E}_{\bm{X}}(\lambda)\right)\bm{U}^{-1}.
\end{eqnarray}   
therefore, we can obtain the operator $\bm{Y}$ from the operator $\bm{X}$ by 
\begin{eqnarray}\label{eq2: prop1: Analogous Operators inf}
\bm{U}\bm{X}\bm{U}^{-1}&=& \bm{Y}.
\end{eqnarray}   
$\hfill\Box$

The next proposition is about the commutative condition for two analogous operators.

\begin{proposition}\label{prop2: Analogous Operators}
Two analogous operators $\bm{X}$ and $\bm{Y}$ are given by the relationship $\bm{X} \propto^{\bm{I}}_{c(\lambda)} \bm{Y}$, then, these two operators $\bm{X}$ and $\bm{Y}$ are commute under the multiplication. 
\end{proposition}
\textbf{Proof:}
Becuase two analogous operators $\bm{X}$ and $\bm{Y}$ are given by the relationship $\bm{X} \propto^{\bm{I}}_{c(\lambda)} \bm{Y}$, we have
\begin{eqnarray}\label{eq1: prop2: Analogous Operators}
\bm{X}&=&\left(\int\limits_{\lambda \in \sigma(\bm{X})}\lambda d\bm{E}_{\bm{X}}(\lambda)+
\int\limits_{\lambda \in \sigma(\bm{X})}\left(\bm{X}-\lambda\bm{I}\right)d\bm{E}_{\bm{X}}(\lambda)\right), 
\end{eqnarray}   
and
\begin{eqnarray}\label{eq2: prop2: Analogous Operators}
\bm{Y}&=&\bm{I}\left(\int\limits_{\lambda \in \sigma(\bm{X})}c(\lambda)\lambda d\bm{E}_{\bm{X}}(\lambda)+
\int\limits_{\lambda \in \sigma(\bm{X})}\left(\bm{X}-c(\lambda)\lambda\bm{I}\right)d\bm{E}_{\bm{X}}(\lambda)\right)\bm{I}^{-1}.
\end{eqnarray}   

Moreover, from Lemma~\ref{lma: orthogonal of spectral basis}, we have 
\begin{eqnarray}\label{eq4: prop2: Analogous Operators}
\bm{X}\bm{Y}&=&\int\limits_{\lambda \in \sigma(\bm{X})}c(\lambda)\lambda^2 d\bm{E}_{\bm{X}}(\lambda)+\int\limits_{\lambda \in \sigma(\bm{X})}\lambda\left(\bm{X}-c(\lambda)\lambda\bm{I}\right)d\bm{E}_{\bm{X}}(\lambda)\nonumber \\
&&+\int\limits_{\lambda \in \sigma(\bm{X})}c(\lambda)\lambda\left(\bm{X}-\lambda\bm{I}\right)d\bm{E}_{\bm{X}}(\lambda)+\int\limits_{\lambda \in \sigma(\bm{X})}\left(\bm{X}-\lambda\bm{I}\right)\left(\bm{X}-c(\lambda)\lambda\bm{I}\right)d\bm{E}_{\bm{X}}(\lambda) \nonumber \\
&=&  \bm{Y}\bm{X}
\end{eqnarray}   
$\hfill\Box$

The next proposition is to show that the proposed analogous relationship between $\bm{X}$ and $\bm{Y}$ provided by $\bm{X} \propto^{\bm{U}}_{c(\lambda)} \bm{Y}$ is an equivalence relation.

\begin{proposition}\label{prop3: Analogous Operators}
The relationship $\bm{X} \propto^{\bm{U}}_{c(\lambda)} \bm{Y}$ is an equivalence relation.
\end{proposition}
\textbf{Proof:}
Given the operators $\bm{X}, \bm{Y}$ and $\bm{Z}$ with the following formats:
\begin{eqnarray}\label{eq1: prop3: Analogous Operators inf}
\bm{X}&=&\left(\int\limits_{\lambda \in \sigma(\bm{X})}\lambda d\bm{E}_{\bm{X}}(\lambda)+
\int\limits_{\lambda \in \sigma(\bm{X})}\left(\bm{X}-\lambda\bm{I}\right)d\bm{E}_{\bm{X}}(\lambda)\right), \nonumber \\
\bm{Y}&=&\bm{U}\left(\int\limits_{\lambda \in \sigma(\bm{X})}c_1(\lambda)\lambda d\bm{E}_{\bm{X}}(\lambda)+
\int\limits_{\lambda \in \sigma(\bm{X})}\left(\bm{X}-c_1(\lambda)\lambda\bm{I}\right)d\bm{E}_{\bm{X}}(\lambda)\right)\bm{U}^{-1},\nonumber \\
\bm{Z}&=&\bm{V}\left(\int\limits_{\lambda \in \sigma(\bm{X})}c_2(\lambda)\lambda d\bm{E}_{\bm{X}}(\lambda)+
\int\limits_{\lambda \in \sigma(\bm{X})}\left(\bm{X}-c_2(\lambda)\lambda\bm{I}\right)d\bm{E}_{\bm{X}}(\lambda)\right)\bm{V}^{-1}
\end{eqnarray} 
we have 
\begin{eqnarray}\label{eq2s: prop3: Analogous Operators inf}
\bm{X}&\propto^{\bm{I}}_{c_1(\lambda)=1}&\bm{X}, ~(\mbox{reflexive})\nonumber \\
\bm{X}&\propto^{\bm{U}}_{c_1(\lambda)}&\bm{Y}, \mbox{~and $\bm{Y}\propto^{\bm{U}^{-1}}_{c_1(\lambda)^{-1}}\bm{X}$}~(\mbox{symmetric})\nonumber \\
\bm{X}&\propto^{\bm{U}}_{c_1(\lambda)}&\bm{Y}, \mbox{~and $\bm{Y}\propto^{\bm{U}^{-1}\bm{V}}_{c_2(\lambda)/c_1(\lambda)}\bm{Z}$}\nonumber \\
&\Longrightarrow& \bm{X}\propto^{\bm{V}}_{c_2(\lambda)}\bm{Z}~(\mbox{transitive})
\end{eqnarray} 
$\hfill\Box$

\subsection{Spectral Mapping Theorem for Multivariate Continous Spectrum Operators}\label{sec: Spectral Mapping Theorem for Multivariate Continous Spectrum Operators}

We first consider spectral mapping theorem for two input operators in Theorem~\ref{thm: Spectral Mapping Theorem for Two Variables inf}.

\begin{theorem}\label{thm: Spectral Mapping Theorem for Two Variables inf}
Given an analytic function $f(z_1,z_2)$ within the domain for $|z_1| < R_1$ and $|z_2| < R_2$, the first operator $\bm{X}_1$ decomposed by:
\begin{eqnarray}\label{eq1-1: thm: Spectral Mapping Theorem for Two Variables inf}
\bm{X}_1&=&\int\limits_{\lambda_1 \in \sigma(\bm{X}_1)}\lambda_1 d\bm{E}_{\bm{X}_1}(\lambda_1)+
\int\limits_{\lambda_1 \in \sigma(\bm{X}_1)}\left(\bm{X}_1-\lambda_1\bm{I}\right)d\bm{E}_{\bm{X}_1}(\lambda_1),
\end{eqnarray}
where $\left\vert\lambda_{1}\right\vert<R_1$, and second operator $\bm{X}_2$ decomposed by:
\begin{eqnarray}\label{eq1-2: thm: Spectral Mapping Theorem for Two Variables inf}
\bm{X}_2&=&\int\limits_{\lambda_2 \in \sigma(\bm{X}_2)}\lambda_2 d\bm{E}_{\bm{X}_2}(\lambda_2)+
\int\limits_{\lambda_2 \in \sigma(\bm{X}_2)}\left(\bm{X}_2-\lambda_2\bm{I}\right)d\bm{E}_{\bm{X}_2}(\lambda_2),
\end{eqnarray}
where $\left\vert\lambda_{2}\right\vert<R_2$.

Then, we have
\begin{eqnarray}\label{eq2: thm: Spectral Mapping Theorem for Two Variables inf}
f(\bm{X}_1, \bm{X}_2)&=&\int\limits_{\lambda_1 \in \sigma(\bm{X}_1)}\int\limits_{\lambda_2 \in \sigma(\bm{X}_2)}f(\lambda_{1}, \lambda_{2})d\bm{E}_{\bm{X}_1}(\lambda_1)d\bm{E}_{\bm{X}_2}(\lambda_2) \nonumber \\
&&+\int\limits_{\lambda_1 \in \sigma(\bm{X}_1)}\int\limits_{\lambda_2 \in \sigma(\bm{X}_2)}\sum_{q_2=1}^{m_{\lambda_2}-1}\frac{f^{(-,q_2)}(\lambda_{1},\lambda_{2})}{q_2!}d\bm{E}_{\bm{X}_1}(\lambda_1)\left(\bm{X}_2-\lambda_2\bm{I}\right)^{q_2}d\bm{E}_{\bm{X}_2}(\lambda_2) \nonumber \\
&&+\int\limits_{\lambda_1 \in \sigma(\bm{X}_1)}\int\limits_{\lambda_2 \in \sigma(\bm{X}_2)}\sum_{q_1=1}^{m_{\lambda_1}-1}\frac{f^{(q_1,-)}(\lambda_{1},\lambda_{2})}{q_1!}\left(\bm{X}_1-\lambda_1\bm{I}\right)^{q_1}d\bm{E}_{\bm{X}_1}(\lambda_1)d\bm{E}_{\bm{X}_2}(\lambda_2)\nonumber \\
&&+\int\limits_{\lambda_1 \in \sigma(\bm{X}_1)}\int\limits_{\lambda_2 \in \sigma(\bm{X}_2)}\sum_{q_1=1}^{m_{\lambda_1}-1}\sum_{q_2=1}^{m_{\lambda_2}-1}\nonumber \\
&&\frac{f^{(q_1,q_2)}(\lambda_{1},\lambda_{2})}{q_1!q_2!}\left(\bm{X}_1-\lambda_1\bm{I}\right)^{q_1}d\bm{E}_{\bm{X}_1}(\lambda_1)\left(\bm{X}_2-\lambda_2\bm{I}\right)^{q_2}d\bm{E}_{\bm{X}_2}(\lambda_2).
\end{eqnarray}
\end{theorem}
\textbf{Proof:}
From Lemma~\ref{lma: orthogonal of spectral basis}, we have the following relationships for a spectral operator $\bm{X}_1$:
\begin{eqnarray}
d\bm{E}_{\bm{X}_1}(\lambda_1)d\bm{E}_{\bm{X}_1}(\lambda_1')&=&d\bm{E}_{\bm{X}_1}(\lambda_1)\delta(\lambda_1,\lambda_1'), \nonumber 
\end{eqnarray}
\begin{eqnarray}
(\bm{X}_1-\lambda_1\bm{I})d\bm{E}_{\bm{X}_1}(\lambda_1)(\bm{X}_1-\lambda_1'\bm{I})d\bm{E}_{\bm{X}_1}(\lambda_1')=(\bm{X}_1-\lambda_1\bm{I})^2d\bm{E}_{\bm{X}_1}(\lambda_1)\delta(\lambda_1,\lambda_1'),\nonumber 
\end{eqnarray} 
\begin{eqnarray}\label{eq3-1: thm: Spectral Mapping Theorem for Two Variables inf}
d\bm{E}_{\bm{X}_1}(\lambda_1')(\bm{X}_1-\lambda_1\bm{I})d\bm{E}_{\bm{X}_1}(\lambda_1)&=&(\bm{X}_1-\lambda_1\bm{I})d\bm{E}_{\bm{X}_1}(\lambda)d\bm{E}_{\bm{X}_1}(\lambda_1') \nonumber \\
&=&(\bm{X}_1-\lambda_1\bm{I})d\bm{E}_{\bm{X}_1}(\lambda_1)\delta(\lambda_1,\lambda_1').
\end{eqnarray} 
Similarly, we also have the following relationships for a spectral operator $\bm{X}_2$:
\begin{eqnarray}
d\bm{E}_{\bm{X}_2}(\lambda_2)d\bm{E}_{\bm{X}_2}(\lambda_2')&=&d\bm{E}_{\bm{X}_2}(\lambda_2)\delta(\lambda_2,\lambda_2'), \nonumber 
\end{eqnarray}
\begin{eqnarray}
(\bm{X}_2-\lambda_2\bm{I})d\bm{E}_{\bm{X}_2}(\lambda_2)(\bm{X}_2-\lambda_2'\bm{I})d\bm{E}_{\bm{X}_2}(\lambda_2')=(\bm{X}_2-\lambda_2\bm{I})^2d\bm{E}_{\bm{X}_2}(\lambda_2)\delta(\lambda_2,\lambda_2'),\nonumber 
\end{eqnarray} 
\begin{eqnarray}\label{eq3-2: thm: Spectral Mapping Theorem for Two Variables inf}
d\bm{E}_{\bm{X}_2}(\lambda_2')(\bm{X}_2-\lambda_2\bm{I})d\bm{E}_{\bm{X}_2}(\lambda_2)&=&(\bm{X}_2-\lambda_2\bm{I})d\bm{E}_{\bm{X}_2}(\lambda)d\bm{E}_{\bm{X}_2}(\lambda_2') \nonumber \\
&=&(\bm{X}_2-\lambda_2\bm{I})d\bm{E}_{\bm{X}_2}(\lambda_2)\delta(\lambda_2,\lambda_2').
\end{eqnarray}

Because $f(z_1, z_2)$ is an analytic function, we have
\begin{eqnarray}\label{eq3: thm: Spectral Mapping Theorem for Two Variables inf}
f(\bm{X}_1,\bm{X}_2)&=&\sum\limits_{\ell_1=0,\ell_2=0}^{\infty}a_{\ell_1,\ell_2} \bm{X}_1^{\ell_1}\bm{X}_2^{\ell_2}\nonumber \\
&=&\sum\limits_{\ell_1=0,\ell_2=0}^{\infty}a_{\ell_1,\ell_2}\left(\int\limits_{\lambda_1 \in \sigma(\bm{X}_1)}\lambda_1 d\bm{E}_{\bm{X}_1}(\lambda_1)+
\int\limits_{\lambda_1 \in \sigma(\bm{X}_1)}\left(\bm{X}_1-\lambda_1\bm{I}\right)d\bm{E}_{\bm{X}_1}(\lambda_1)\right)^{\ell_1}\nonumber \\
&&\times \left(\int\limits_{\lambda_2 \in \sigma(\bm{X}_2)}\lambda_2 d\bm{E}_{\bm{X}_2}(\lambda_2)+
\int\limits_{\lambda_2 \in \sigma(\bm{X}_2)}\left(\bm{X}_2-\lambda_2\bm{I}\right)d\bm{E}_{\bm{X}_2}(\lambda_2)\right)^{\ell_2}\nonumber \\
&=_1&\sum\limits_{\ell_1=0,\ell_2=0}^{\infty}a_{\ell_1,\ell_2}\left[\int\limits_{\lambda_1 \in \sigma(\bm{X}_1)}\left(\lambda_1 d\bm{E}_{\bm{X}_1}(\lambda_1)+\left(\bm{X}_1-\lambda_1\bm{I}\right)d\bm{E}_{\bm{X}_1}(\lambda_1)\right)^{\ell_1}\right]\nonumber \\
&&\times\left[\int\limits_{\lambda_2 \in \sigma(\bm{X}_2)}\left(\lambda_2 d\bm{E}_{\bm{X}_2}(\lambda_2)+\left(\bm{X}_2-\lambda_2\bm{I}\right)d\bm{E}_{\bm{X}_2}(\lambda_2)\right)^{\ell_2}\right]\nonumber \\
&=&\int\limits_{\lambda_1 \in \sigma(\bm{X}_1)}\int\limits_{\lambda_2 \in \sigma(\bm{X}_2)}\sum\limits_{\ell_1=0,\ell_2=0}^{\infty}\Bigg\{a_{\ell_1,\ell_2}\nonumber \\
&&\times\left[\sum\limits_{q_1=0}^{m_{\lambda_1}-1}\frac{\ell_1!}{q_1! (\ell_1-q_1)!}\lambda^{\ell_1-q_1}_{1}\left(d\bm{E}_{\bm{X}_1}(\lambda_1)\right)^{\ell_1-q_1}\left((\bm{X}_1 - \lambda_1\bm{I})d\bm{E}_{\bm{X}_1}(\lambda_1)\right)^{q_1}\right]\nonumber \\
&&\times \left[\sum\limits_{q_2=0}^{m_{\lambda_2}-1}\frac{\ell_2!}{q_2! (\ell_2-q_2)!}\lambda^{\ell_2-q_2}_{2}\left(d\bm{E}_{\bm{X}_2}(\lambda_2)\right)^{\ell_2-q_2}\left((\bm{X}_2 - \lambda_2\bm{I})d\bm{E}_{\bm{X}_2}(\lambda_2)\right)^{q_2}\right]\Bigg\}\nonumber \\
&=&\int\limits_{\lambda_1 \in \sigma(\bm{X}_1)}\int\limits_{\lambda_2 \in \sigma(\bm{X}_2)}\sum\limits_{\ell_1=0,\ell_2=0}^{\infty}\Bigg\{a_{\ell_1,\ell_2}\nonumber \\
&&\times\left[\lambda_{1}^{\ell_1}\left(d\bm{E}_{\bm{X}_1}(\lambda_1)\right)^{\ell_1}\right.\nonumber \\
&&\left.~~~~+\sum\limits_{q_1=1}^{m_{\lambda_1}-1}\frac{\ell_1!}{q_1! (\ell_1-q_1)!}\lambda^{\ell_1-q_1}_{1}\left(d\bm{E}_{\bm{X}_1}(\lambda_1)\right)^{\ell_1-q_1}\left((\bm{X}_1 - \lambda_1\bm{I})d\bm{E}_{\bm{X}_1}(\lambda_1)\right)^{q_1}\right]\nonumber \\
&&\times \left[\lambda_{2}^{\ell_2}\left(d\bm{E}_{\bm{X}_1}(\lambda_1)\right)^{\ell_2}\right.\nonumber \\
&&\left.~~~~+\sum\limits_{q_2=1}^{m_{\lambda_2}-1}\frac{\ell_2!}{q_2! (\ell_2-q_2)!}\lambda^{\ell_2-q_2}_{2}\left(d\bm{E}_{\bm{X}_2}(\lambda_2)\right)^{\ell_2-q_2}\left((\bm{X}_2 - \lambda_2\bm{I})d\bm{E}_{\bm{X}_2}(\lambda_2)\right)^{q_2}\right]\Bigg\}\nonumber 
\end{eqnarray}
\begin{eqnarray}
&=_2&\int\limits_{\lambda_1 \in \sigma(\bm{X}_1)}\int\limits_{\lambda_2 \in \sigma(\bm{X}_2)}f(\lambda_{1}, \lambda_{2})d\bm{E}_{\bm{X}_1}(\lambda_1)d\bm{E}_{\bm{X}_2}(\lambda_2) \nonumber \\
&&+\int\limits_{\lambda_1 \in \sigma(\bm{X}_1)}\int\limits_{\lambda_2 \in \sigma(\bm{X}_2)}\sum_{q_2=1}^{m_{\lambda_2}-1}\frac{f^{(-,q_2)}(\lambda_{1},\lambda_{2})}{q_2!}d\bm{E}_{\bm{X}_1}(\lambda_1)\left(\bm{X}_2-\lambda_2\bm{I}\right)^{q_2}d\bm{E}_{\bm{X}_2}(\lambda_2) \nonumber \\
&&+\int\limits_{\lambda_1 \in \sigma(\bm{X}_1)}\int\limits_{\lambda_2 \in \sigma(\bm{X}_2)}\sum_{q_1=1}^{m_{\lambda_1}-1}\frac{f^{(q_1,-)}(\lambda_{1},\lambda_{2})}{q_1!}\left(\bm{X}_1-\lambda_1\bm{I}\right)^{q_1}d\bm{E}_{\bm{X}_1}(\lambda_1)d\bm{E}_{\bm{X}_2}(\lambda_2)\nonumber \\
&&+\int\limits_{\lambda_1 \in \sigma(\bm{X}_1)}\int\limits_{\lambda_2 \in \sigma(\bm{X}_2)}\sum_{q_1=1}^{m_{\lambda_1}-1}\sum_{q_2=1}^{m_{\lambda_2}-1}\nonumber \\
&&\frac{f^{(q_1,q_2)}(\lambda_{1},\lambda_{2})}{q_1!q_2!}\left(\bm{X}_1-\lambda_1\bm{I}\right)^{q_1}d\bm{E}_{\bm{X}_1}(\lambda_1)\left(\bm{X}_2-\lambda_2\bm{I}\right)^{q_2}d\bm{E}_{\bm{X}_2}(\lambda_2),
\end{eqnarray}
where we apply Eq.~\eqref{eq3-1: thm: Spectral Mapping Theorem for Two Variables inf} and Eq.~\eqref{eq3-2: thm: Spectral Mapping Theorem for Two Variables inf} in $=_1$, and $=_2$, and following relations about the partial derivatives of $f(z_1,z_2)$ in $=_2$:
\begin{eqnarray}\label{eq4: thm: Spectral Mapping Theorem for Two Variables inf}
f^{(-,q_2)}(z_1,z_2)&=&\sum_{\ell_1=0}^{\infty}\sum_{\ell_2=q_2}^{\infty}\frac{a_{\ell_1,\ell_2}\ell_2!}{(\ell_2-q_2)!}z_1^{\ell_1}z_2^{\ell_2-q_2},\nonumber\\
f^{(q_1,-)}(z_1,z_2)&=&\sum_{\ell_1=q_1}^{\infty}\sum_{\ell_2=0}^{\infty}\frac{a_{\ell_1,\ell_2}\ell_1!}{(\ell_1-q_1)!}z_1^{\ell_1-q_1}z_2^{\ell_2},\nonumber\\
f^{(q_1,q_2)}(z_1,z_2)&=&\sum_{\ell_1=q_1}^{\infty}\sum_{\ell_2=q_2}^{\infty}\frac{a_{\ell_1,\ell_2}\ell_1!\ell_2!}{(\ell_1-q_1)!(\ell_2-q_2)!}z_1^{\ell_1-q_1}z_2^{\ell_2-q_2}.
\end{eqnarray}
$\hfill \Box$

We are ready to present Theorem~\ref{thm: Spectral Mapping Theorem for r Variables inf} about spectral mapping theorem for $r$ spectral operators. 
\begin{theorem}\label{thm: Spectral Mapping Theorem for r Variables inf}
Given an analytic function $f(z_1,z_2,\ldots,z_r)$ within the domain for $|z_l| < R_l$, and the operator $\bm{X}_l$ decomposed by:
\begin{eqnarray}\label{eq1-1: thm: Spectral Mapping Theorem for r Variables inf}
\bm{X}_l&=&\int\limits_{\lambda_l \in \sigma(\bm{X}_l)}\lambda_l d\bm{E}_{\bm{X}_l}(\lambda_l)+
\int\limits_{\lambda_l \in \sigma(\bm{X}_l)}\left(\bm{X}_l-\lambda_l\bm{I}\right)d\bm{E}_{\bm{X}_l}(\lambda_l),
\end{eqnarray}
where $\left\vert\lambda_{l}\right\vert<R_l$ for $l=1,2,\ldots,r$.

Then, we have
\begin{eqnarray}\label{eq2: thm: Spectral Mapping Theorem for kappa Variables inf}
\lefteqn{f(\bm{X}_1,\ldots,\bm{X}_r)=}\nonumber \\
&&\int\limits_{\lambda_1 \in \sigma(\bm{X}_1)}\cdots\int\limits_{\lambda_r \in \sigma(\bm{X}_r)}
f(\lambda_1,\ldots,\lambda_r)d\bm{E}_{\bm{X}_1}(\lambda_1)\cdots d\bm{E}_{\bm{X}_r}(\lambda_r) \nonumber \\
&&+\int\limits_{\lambda_1 \in \sigma(\bm{X}_1)}\cdots\int\limits_{\lambda_r \in \sigma(\bm{X}_r)}\sum\limits_{\kappa=1}^{r-1}\sum\limits_{\alpha_\kappa(q_1,\ldots,q_r)}\Bigg(\sum\limits_{\alpha_{\kappa}(q_1,\ldots,q_r)=1}^{m_{\lambda_{\mbox{Ind}(\alpha_{\kappa}(q_1,\ldots,q_r))}}-1}\nonumber \\
&&~~~~~ \frac{f^{\alpha_{\kappa}(q_1,\ldots,q_r)}(\lambda_1,\ldots,\lambda_r)}{q_{\iota_1}!q_{\iota_2}!\ldots q_{\iota_\kappa}!}\times \prod\limits_{\substack{\beta =\mbox{Ind}(\alpha_{\kappa}(q_1,\ldots,q_r)), \bm{Y}=\left(\bm{X}_\beta - \lambda_\beta \bm{I}\right)^{q_\beta}d\bm{E}_{\bm{X}_\beta}(\lambda_\beta) \\ \beta \neq \mbox{Ind}(\alpha_{\kappa}(q_1,\ldots,q_r)), \bm{Y}=d\bm{E}_{\bm{X}_\beta}(\lambda_\beta)}
}^{r} \bm{Y}\Bigg) 
\nonumber \\
&&+\int\limits_{\lambda_1 \in \sigma(\bm{X}_1)}\cdots\int\limits_{\lambda_r \in \sigma(\bm{X}_r)}\sum\limits_{q_1=\ldots=q_r=1}^{m_{\lambda_1}-1,\ldots,m_{\lambda_r}-1}
\frac{f^{(q_1,\ldots,q_r)}(\lambda_1,\ldots,\lambda_r)}{q_1!\cdots q_r!}\nonumber \\
&&\times \left(\bm{X}_1 - \lambda_1\bm{I}\right)^{q_1}d\bm{E}_{\bm{X}_1}(\lambda_1) \left(\bm{X}_2 - \lambda_2\bm{I}\right)^{q_2}d\bm{E}_{\bm{X}_2}(\lambda_2)\cdots \left(\bm{X}_r - \lambda_r\bm{I}\right)^{q_r}d\bm{E}_{\bm{X}_r}(\lambda_r),
\end{eqnarray}
where we have
\begin{itemize}
\item $\sum\limits_{\alpha_\kappa(q_1,\ldots,q_r)}$ runs over all $\kappa$ selections of $q_1,\ldots,q_r$ by $\alpha_\kappa(q_1,\ldots,q_r)$;
\item $m_{\lambda_{\mbox{Ind}(\alpha_{\kappa}(q_1,\ldots,q_r))}}-1 =$ $m_{\lambda_{\iota_1}}-1$,$\ldots,m_{\lambda_{\iota_\kappa}}-1$;
\item $f^{\alpha_{\kappa}(q_1,\ldots,q_r)}(\lambda_1,\ldots,\lambda_r)$ represents the partial derivatives with respect to variables with indices $\iota_1,\iota_2,\ldots,\iota_\kappa$ and the orders of derivatives given by $q_{\iota_1},q_{\iota_2},\ldots,q_{\iota_\kappa}$.
\end{itemize}
\end{theorem}
\textbf{Proof:}
The proof follows directly from Theorem~\ref{thm: Spectral Mapping Theorem for Two Variables inf}, with the remaining steps consisting primarily of routine notational manipulations.
$\hfill\Box$

\section{Hybrid Spectrum Operators}\label{sec: Hybrid Spectrum Operators}

Operators that have both discrete and continuous eigenvalues arise in several physical and mathematical contexts, particularly in quantum mechanics, spectral theory, and differential operators on non-compact domains. Operators with mixed spectra often appear in systems where there is a combination of confined (bounded) and unbounded regions, or where localized perturbations create bound states in an otherwise continuous spectrum. The resulting mixed spectra reflect both discrete states (e.g., bound states or modes) and a continuous range of states (e.g., free or scattering states), common in many physical and mathematical settings. In this section, we examine spectral properties of operators that feature both discrete eigenvalues and continuous spectra. 

\subsection{Spectral Mapping Theorem for Hybrid Spectrum Operators}\label{sec: Spectral Mapping Theorem for Hybrid Spectrum Operators}

In this section, spectral mapping theorem for hybrid operator and its functional calculus will be presented. We still focus in discussing spectral operators. We first introduce the following Theorem~\ref{thm: spectral theorem of inf-dim h} about the spectral decomposition.

\begin{theorem}\label{thm: spectral theorem of inf-dim h}
Given a spectral hybrid operator $\bm{X}$, we can express it by
\begin{eqnarray}\label{eq0: thm: spectral theorem of inf-dim h}
\bm{X}&=&\sum_{\lambda_i \in \sigma_d(\bm{X})}\lambda_i \bm{F}_{\bm{X}}(\lambda_i)+\int\limits_{\lambda \in \sigma_c(\bm{X})}\lambda d\bm{E}_{\bm{X}}(\lambda)\nonumber \\
&&+
\sum_{\lambda_i \in \sigma_d(\bm{X})}\left(\bm{X}-\lambda_i\bm{I}\right)\bm{F}_{\bm{X}}(\lambda_i)+\int\limits_{\lambda \in \sigma_c(\bm{X})}\left(\bm{X}-\lambda\bm{I}\right)d\bm{E}_{\bm{X}}(\lambda),
\end{eqnarray}
where $\bm{F}_{\bm{X}}(\lambda_i)$ are the projectors corresponding to the eigenvalue $\lambda_i$, $\sigma_d(\bm{X})$ and $\sigma_c(\bm{X})$ are sets of the discrete eigenvalues and continous spectra of the operator $\bm{X}$, respectively. 
\end{theorem}
\textbf{Proof:}
From~\cite{dunford1958survey}, we have the spectrum of $\bm{X}$ identical to the spectrum of 
\begin{eqnarray}\label{eq0.1: thm: spectral theorem of inf-dim h}
\sum_{\lambda_i \in \sigma_d(\bm{X})}\lambda_i \bm{F}_{\bm{X}}(\lambda_i) + \int\limits_{\lambda \in \sigma_c(\bm{X})}\lambda d\bm{E}_{\bm{X}}(\lambda),
\end{eqnarray}
then, we have the following decomposition
\begin{eqnarray}\label{eq1: thm: spectral theorem of inf-dim h}
\bm{X}&=&\sum_{\lambda_i \in \sigma_d(\bm{X})}\lambda_i \bm{F}_{\bm{X}}(\lambda_i) + \int\limits_{\lambda \in \sigma_c(\bm{X})}\lambda d\bm{E}_{\bm{X}}(\lambda)+\bm{N}_{\bm{X}},
\end{eqnarray}
where $\bm{N}_{\bm{X}}$ is the nilpotent part of the operator $\bm{X}$. 

Therefore, we have
\begin{eqnarray}\label{eq2: thm: spectral theorem of inf-dim h}
\bm{N}_{\bm{X}}&=&\bm{X}-\left[\sum_{\lambda_i \in \sigma_d(\bm{X})}\lambda_i \bm{F}_{\bm{X}}(\lambda_i) + \int\limits_{\lambda \in \sigma_c(\bm{X})}\lambda d\bm{E}_{\bm{X}}(\lambda)\right]\nonumber \\
&=_1&\bm{X}\left[\sum_{\lambda_i \in \sigma_d(\bm{X})}\bm{F}_{\bm{X}}(\lambda_i)+\int\limits_{\lambda \in \sigma_c(\bm{X})}d\bm{E}_{\bm{X}}(\lambda)\right]-\left[\sum_{\lambda_i \in \sigma_d(\bm{X})}\lambda_i \bm{F}_{\bm{X}}(\lambda_i) + \int\limits_{\lambda \in \sigma_c(\bm{X})}\lambda d\bm{E}_{\bm{X}}(\lambda)\right]\nonumber \\
&=&\sum_{\lambda_i \in \sigma_d(\bm{X})}\left(\bm{X}-\lambda_i\bm{I}\right)\bm{F}_{\bm{X}}(\lambda_i)+\int\limits_{\lambda \in \sigma_c(\bm{X})}\left(\bm{X}-\lambda\bm{I}\right)d\bm{E}_{\bm{X}}(\lambda),
\end{eqnarray}
where we apply the property 
\begin{eqnarray}\label{eq3: thm: spectral theorem of inf-dim h}
\bm{I}=\left[\sum_{\lambda_i \in \sigma_d(\bm{X})}\bm{F}_{\bm{X}}(\lambda_i)+\int\limits_{\lambda \in \sigma_c(\bm{X})}d\bm{E}_{\bm{X}}(\lambda)\right]
\end{eqnarray}
in $=_1$. This thoerem is proved. 
$\hfill\Box$

Lemma~\ref{lma: orthogonal of spectral basis h} below will give relationships among $\bm{F}(\lambda_i)$, $\left(\bm{X}-\lambda_i\bm{I}\right)\bm{F}_{\bm{X}}(\lambda_i)$,$d\bm{E}_{\bm{X}}(\lambda)$ and $(\bm{X}-\lambda\bm{I})d\bm{E}_{\bm{X}}(\lambda)$ which will be used to establish the functional calculus for the operator $\bm{X}$.

\begin{lemma}\label{lma: orthogonal of spectral basis h}
We have the following relationships for a spectral operator $\bm{X}$:
\begin{eqnarray}\label{eq1: lma: orthogonal of spectral basis h}
\bm{F}_{\bm{X}}(\lambda_i)\bm{F}_{\bm{X}}(\lambda_{i'})&=&\bm{F}_{\bm{X}}(\lambda_i)\delta(\lambda_i,\lambda_{i'}),
\end{eqnarray}
\begin{eqnarray}\label{eq2: lma: orthogonal of spectral basis h}
(\bm{X}-\lambda_i\bm{I})\bm{F}_{\bm{X}}(\lambda_i)(\bm{X}-\lambda_{i'}\bm{I})\bm{F}_{\bm{X}}(\lambda_{i'})=(\bm{X}-\lambda_i\bm{I})^2\bm{F}_{\bm{X}}(\lambda_i)\delta(\lambda_i,\lambda_{i'}),
\end{eqnarray} 
\begin{eqnarray}\label{eq3: lma: orthogonal of spectral basis h}
\bm{F}_{\bm{X}}(\lambda_{i'})(\bm{X}-\lambda_i\bm{I})\bm{F}_{\bm{X}}(\lambda_i)&=&(\bm{X}-\lambda_i\bm{I})\bm{F}_{\bm{X}}(\lambda_i)\bm{F}_{\bm{X}}(\lambda_{i'}) \nonumber \\
&=&(\bm{X}-\lambda_i\bm{I})\bm{F}_{\bm{X}}(\lambda_i)\delta(\lambda_i,\lambda_{i'}).
\end{eqnarray} 
\begin{eqnarray}\label{eq4: lma: orthogonal of spectral basis h}
\bm{F}_{\bm{X}}(\lambda_i)d\bm{E}_{\bm{X}}(\lambda)&=&d\bm{E}_{\bm{X}}(\lambda)\bm{F}_{\bm{X}}(\lambda_i)=\bm{0}\nonumber \\
\bm{F}_{\bm{X}}(\lambda_i)(\bm{X}-\lambda\bm{I})d\bm{E}_{\bm{X}}(\lambda)&=&(\bm{X}-\lambda\bm{I})d\bm{E}_{\bm{X}}(\lambda)\bm{F}_{\bm{X}}(\lambda_i)=\bm{0}\nonumber \\
\left(\bm{X}-\lambda_i\bm{I}\right)\bm{F}_{\bm{X}}(\lambda_i)d\bm{E}_{\bm{X}}(\lambda)&=&d\bm{E}_{\bm{X}}(\lambda)\left(\bm{X}-\lambda_i\bm{I}\right)\bm{F}_{\bm{X}}(\lambda_i)=\bm{0}\nonumber \\
\left(\bm{X}-\lambda_i\bm{I}\right)\bm{F}_{\bm{X}}(\lambda_i)(\bm{X}-\lambda\bm{I})d\bm{E}_{\bm{X}}(\lambda)&=&(\bm{X}-\lambda\bm{I})d\bm{E}_{\bm{X}}(\lambda)\left(\bm{X}-\lambda_i\bm{I}\right)\bm{F}_{\bm{X}}(\lambda_i)=\bm{0}.
\end{eqnarray} 
\end{lemma}
\textbf{Proof:}
Since $\lambda_i$, $\lambda_{i'}$ and the neighborhood region around $\lambda$ are exclusive, this Lemma follows from the proof in Lemma~\ref{lma: orthogonal of spectral basis}.
$\hfill\Box$

Theorem~\ref{thm: functional cal of single inf-dim h} below is the main result of this section.

\begin{theorem}\label{thm: functional cal of single inf-dim h}
Given an analytic function $f(z)$ within the domain for $|z| < R$, a spectral operator $\bm{X}$ with the following spectral decomposition
\begin{eqnarray}\label{eq1: thm: functional cal of single inf-dim h}
\bm{X}&=&\sum_{\lambda_i \in \sigma_d(\bm{X})}\lambda_i \bm{F}_{\bm{X}}(\lambda_i)+\int\limits_{\lambda \in \sigma_c(\bm{X})}\lambda d\bm{E}_{\bm{X}}(\lambda)\nonumber \\
&&+
\sum_{\lambda_i \in \sigma_d(\bm{X})}\left(\bm{X}-\lambda_i\bm{I}\right)\bm{F}_{\bm{X}}(\lambda_i)+\int\limits_{\lambda \in \sigma_c(\bm{X})}\left(\bm{X}-\lambda\bm{I}\right)d\bm{E}_{\bm{X}}(\lambda),
\end{eqnarray}
where $|\lambda_i|<R$ for $\lambda_i \in \sigma_d(\bm{X})$ and $|\lambda|<R$ for $\lambda \in \sigma_c(\bm{X})$, then, we have
\begin{eqnarray}\label{eq2: thm: functional cal of single inf-dim h}
f(\bm{X})&=&\sum\limits_{\lambda_i \in \sigma_d(\bm{X})}f(\lambda_i)\bm{F}_{\bm{X}}(\lambda_i)+\int\limits_{\lambda \in \sigma_c(\bm{X})}f(\lambda)d\bm{E}_{\bm{X}}(\lambda)\nonumber \\
&& + \sum\limits_{\lambda_i \in \sigma_d(\bm{X})}\sum\limits_{q=1}^{m_{\lambda_i}-1}\frac{f^{(q)}(\lambda_i)}{q!}(\bm{X}-\lambda_i\bm{I})^q \bm{F}_{\bm{X}}(\lambda_i)\nonumber \\
&& + \int\limits_{\lambda \in \sigma_c(\bm{X})}\sum\limits_{q=1}^{m_{\lambda}-1}\frac{f^{(q)}(\lambda)}{q!}(\bm{X}-\lambda\bm{I})^q d\bm{E}_{\bm{X}}(\lambda),
\end{eqnarray}
where $m_{\lambda_i}$ is a positive integer to measure the degree of nilpotency with respect to the eigenvalue $\lambda_i$, i.e., 
\begin{eqnarray}\label{eq2.5: thm: functional cal of single inf-dim h}
(\bm{X}-\lambda_i\bm{I})^q \bm{F}_{\bm{X}}(\lambda_i)&& 
    \begin{cases}
      \neq \bm{0} & \text{if $q < m_{\lambda_i}$,}\\
      =\bm{0} & \text{if $q \geq m_{\lambda_i}$;}
    \end{cases}  
\end{eqnarray}
and
$m_{\lambda}$ is a positive integer to measure the degree of nilpotency with respect to the eigenvalue $\lambda$, i.e., 
\begin{eqnarray}\label{eq2.6: thm: functional cal of single inf-dim h}
(\bm{X}-\lambda\bm{I})^q d\bm{E}_{\bm{X}}(\lambda)&& 
    \begin{cases}
      \neq \bm{0} & \text{if $q < m_{\lambda}$,}\\
      =\bm{0} & \text{if $q \geq m_{\lambda}$.}
    \end{cases}  
\end{eqnarray}

\end{theorem}
\textbf{Proof:}

Because $f(z)$ is an analytic function, we have
\begin{eqnarray}\label{eq3: thm: Spectral Mapping Theorem for Single Variable h}
f(\bm{X})&=&\sum\limits_{\ell=0}^{\infty}a_\ell \bm{X}^\ell \nonumber \\
&=&\sum\limits_{\ell=0}^{\infty}a_\ell\left(\sum_{\lambda_i \in \sigma_d(\bm{X})}\lambda_i \bm{F}_{\bm{X}}(\lambda_i)+\int\limits_{\lambda \in \sigma_c(\bm{X})}\lambda d\bm{E}_{\bm{X}}(\lambda)\right.\nonumber \\
&&\left.+
\sum_{\lambda_i \in \sigma_d(\bm{X})}\left(\bm{X}-\lambda_i\bm{I}\right)\bm{F}_{\bm{X}}(\lambda_i)+\int\limits_{\lambda \in \sigma_c(\bm{X})}\left(\bm{X}-\lambda\bm{I}\right)d\bm{E}_{\bm{X}}(\lambda)\right)^\ell \nonumber \\
&=_1&\sum\limits_{\ell=0}^{\infty}\left[\sum\limits_{\lambda_i \in \sigma_d(\bm{X})}a_\ell\left(\lambda_i \bm{F}_{\bm{X}}(\lambda_i)+\left(\bm{X}-\lambda_i\bm{I}\right)\bm{F}_{\bm{X}}(\lambda_i)\right)^\ell \right. \nonumber \\
&&\left.  + \int\limits_{\lambda \in \sigma_c(\bm{X})}a_\ell\left(\lambda d\bm{E}_{\bm{X}}(\lambda)+\left(\bm{X}-\lambda\bm{I}\right)d\bm{E}_{\bm{X}}(\lambda)\right)^\ell\right] \nonumber \\
&=&\int\limits_{\lambda \in \sigma_c(\bm{X})}\sum\limits_{\ell=0}^{\infty}\sum\limits_{q=0}^{\ell} \frac{a_\ell \ell!}{q! (\ell-q)!}\lambda^{\ell-q}\left[d\bm{E}_{\bm{X}}(\lambda)\right]^{\ell-q}\left[(\bm{X}-\lambda\bm{I})d\bm{E}_{\bm{X}}(\lambda)\right]^q\nonumber \\
& &+ \sum\limits_{\lambda_i \in \sigma_d(\bm{X})}\sum\limits_{\ell=0}^{\infty}\sum\limits_{q=0}^{\ell} \frac{a_\ell \ell!}{q! (\ell-q)!}\lambda_i^{\ell-q}\left[\bm{F}_{\bm{X}}(\lambda_i)\right]^{\ell-q}\left[(\bm{X}-\lambda_i\bm{I})\bm{F}_{\bm{X}}(\lambda_i)\right]^q\nonumber \\
&=&\int\limits_{\lambda \in \sigma_c(\bm{X})}\sum\limits_{\ell=0}^{\infty}\left[a_{\ell}\lambda^{\ell}\left[d\bm{E}_{\bm{X}}(\lambda)\right]^\ell+\sum\limits_{q=1}^{\ell}\frac{a_\ell \ell!}{q! (\ell-q)!}\lambda^{\ell-q}\left[d\bm{E}_{\bm{X}}(\lambda)\right]^{\ell-q}\left[(\bm{X}-\lambda\bm{I})d\bm{E}_{\bm{X}}(\lambda)\right]^q\right]\nonumber \\
&&+ \sum\limits_{\lambda_i \in \sigma_d(\bm{X})}\sum\limits_{\ell=0}^{\infty}\left[a_{\ell}\lambda_i^{\ell}\left[\bm{F}_{\bm{X}}(\lambda_i)\right]^\ell+\sum\limits_{q=1}^{\ell}\frac{a_\ell \ell!}{q! (\ell-q)!}\lambda_i^{\ell-q}\left[\bm{F}_{\bm{X}}(\lambda_i)\right]^{\ell-q}\left[(\bm{X}-\lambda_i\bm{I})\bm{F}_{\bm{X}}(\lambda_i)\right]^q\right]\nonumber \\
&=_2&\int\limits_{\lambda \in \sigma_c(\bm{X})}\left[f(\lambda)d\bm{E}_{\bm{X}}(\lambda)+\sum\limits_{q=1}^{m_{\lambda}-1}\frac{f^{(q)}(\lambda)}{q!}(\bm{X}-\lambda\bm{I})^q d\bm{E}_{\bm{X}}(\lambda)\right]\nonumber \\
&&+ \sum\limits_{\lambda_i \in \sigma_d(\bm{X})}\left[f(\lambda_i)\bm{F}_{\bm{X}}(\lambda_i)+\sum\limits_{q=1}^{m_{\lambda_i}-1}\frac{f^{(q)}(\lambda_i)}{q!}(\bm{X}-\lambda_i\bm{I})^q \bm{F}_{\bm{X}}(\lambda_i)\right] \nonumber \\
&=&\int\limits_{\lambda \in \sigma_c(\bm{X})}f(\lambda)d\bm{E}_{\bm{X}}(\lambda)+\int\limits_{\lambda \in \sigma_c(\bm{X})}\sum\limits_{q=1}^{m_{\lambda}-1}\frac{f^{(q)}(\lambda)}{q!}(\bm{X}-\lambda\bm{I})^q d\bm{E}_{\bm{X}}(\lambda) \nonumber \\
&&+\sum\limits_{\lambda_i \in \sigma_d(\bm{X})}f(\lambda_i)\bm{F}_{\bm{X}}(\lambda_i)+\sum\limits_{\lambda_i \in \sigma_d(\bm{X})}\sum\limits_{q=1}^{m_{\lambda_i}-1}\frac{f^{(q)}(\lambda_i)}{q!}(\bm{X}-\lambda_i\bm{I})^q\bm{F}_{\bm{X}}(\lambda_i),
\end{eqnarray}
where we apply Lemma~\ref{lma: orthogonal of spectral basis} and Lemma~\ref{lma: orthogonal of spectral basis h} in $=_1$ and $=_2$, and $f^{(q)}(z) = \sum\limits_{\ell=q}^{\infty}\frac{a_\ell \ell!}{(\ell-q)!}z^{\ell-q}$ in $=_2$.
$\hfill \Box$

\subsection{Analogous Hybrid Spectrum Operators and Their Properties}\label{sec: Analogous Hybrid  Spectrum Operators and Their Properties}

Given two hybrid spectrum operators $\bm{X}$ and $\bm{Y}$, we say that $\bm{X}$ is analogous to $\bm{Y}$ with respect to the invertible operator $\bm{U}$, a sequence of complex numbers $[c_i] \neq 0$ and a ratio function $c(\lambda)$, where $c(\lambda) \neq 0$ for $\lambda \in \sigma_c(\bm{X})$. This relationship, denoted by $\bm{X} \propto^{\bm{U}}_{[c_i], c(\lambda)}\bm{Y}$, holds if $\bm{X}$ and $\bm{Y}$ can be expressed as following:
\begin{eqnarray}\label{eq1: Analogous Operators inf h}
\bm{X}&=&\left(\sum_{\lambda_i \in \sigma_d(\bm{X})}\lambda_i \bm{F}_{\bm{X}}(\lambda_i)+\int\limits_{\lambda \in \sigma_c(\bm{X})}\lambda d\bm{E}_{\bm{X}}(\lambda)\right.\nonumber \\
&&\left.+
\sum_{\lambda_i \in \sigma_d(\bm{X})}\left(\bm{X}-\lambda_i\bm{I}\right)\bm{F}_{\bm{X}}(\lambda_i)+\int\limits_{\lambda \in \sigma_c(\bm{X})}\left(\bm{X}-\lambda\bm{I}\right)d\bm{E}_{\bm{X}}(\lambda),\right), \nonumber \\
\bm{Y}&=& \bm{U}\left(\sum_{\lambda_i \in \sigma_d(\bm{X})}c_i\lambda_i \bm{F}_{\bm{X}}(\lambda_i)+\int\limits_{\lambda \in \sigma_c(\bm{X})}c(\lambda)\lambda d\bm{E}_{\bm{X}}(\lambda)\right.\nonumber \\
&&\left.+
\sum_{\lambda_i \in \sigma_d(\bm{X})}\left(\bm{X}-c_i\lambda_i\bm{I}\right)\bm{F}_{\bm{X}}(\lambda_i)+\int\limits_{\lambda \in \sigma_c(\bm{X})}\left(\bm{X}-c(\lambda)\lambda\bm{I}\right)d\bm{E}_{\bm{X}}(\lambda),\right)\bm{U}^{-1}.
\end{eqnarray}

We will have the following properties about two analogous hybrid spectrum operators $\bm{X}$ and $\bm{Y}$ given by the relationship $\bm{X} \propto^{\bm{U}}_{[c_i],c(\lambda)}\bm{Y}$. Proposition~\ref{prop1: Analogous Operators inf h} is to show that the similarity relationship between two hybrid spectrum operators is the special case of the proposed analogous relationship. 

\begin{proposition}\label{prop1: Analogous Operators inf h}
Two analogous hybrid spectrum operators $\bm{X}$ and $\bm{Y}$ given by the relationship $\bm{X} \propto^{\bm{U}}_{[1],c(\lambda)=1}\bm{Y}$, then, these two hybrid spectrum operators $\bm{X}$ and $\bm{Y}$ are similar.
\end{proposition}
\textbf{Proof:}
Since we have 
\begin{eqnarray}\label{eq1: prop1: Analogous Operators inf h}
\bm{X}&=&\left(\sum_{\lambda_i \in \sigma_d(\bm{X})}\lambda_i \bm{F}_{\bm{X}}(\lambda_i)+\int\limits_{\lambda \in \sigma_c(\bm{X})}\lambda d\bm{E}_{\bm{X}}(\lambda)\right.\nonumber \\
&&\left.+
\sum_{\lambda_i \in \sigma_d(\bm{X})}\left(\bm{X}-\lambda_i\bm{I}\right)\bm{F}_{\bm{X}}(\lambda_i)+\int\limits_{\lambda \in \sigma_c(\bm{X})}\left(\bm{X}-\lambda\bm{I}\right)d\bm{E}_{\bm{X}}(\lambda),\right), \nonumber \\
\bm{Y}&=& \bm{U}\left(\sum_{\lambda_i \in \sigma_d(\bm{X})}1 \times \lambda_i \bm{F}_{\bm{X}}(\lambda_i)+\int\limits_{\lambda \in \sigma_c(\bm{X})}1 \times \lambda d\bm{E}_{\bm{X}}(\lambda)\right.\nonumber \\
&&\left.+
\sum_{\lambda_i \in \sigma_d(\bm{X})}\left(\bm{X}-1 \times \lambda_i\bm{I}\right)\bm{F}_{\bm{X}}(\lambda_i)+\int\limits_{\lambda \in \sigma_c(\bm{X})}\left(\bm{X}-1 \times \lambda\bm{I}\right)d\bm{E}_{\bm{X}}(\lambda),\right)\bm{U}^{-1}.
\end{eqnarray}   
therefore, we can obtain the operator $\bm{Y}$ from the operator $\bm{X}$ by 
\begin{eqnarray}\label{eq2: prop1: Analogous Operators inf}
\bm{U}\bm{X}\bm{U}^{-1}&=& \bm{Y}.
\end{eqnarray}   
$\hfill\Box$

The next proposition is about the commutative condition for two analogous hybrid spectrum operators.

\begin{proposition}\label{prop2: Analogous Operators h}
Two analogous hybrid spectrum operators $\bm{X}$ and $\bm{Y}$ are given by the relationship $\bm{X} \propto^{\bm{I}}_{[c_i], c(\lambda)} \bm{Y}$, then, these two hybrid spectrum operators $\bm{X}$ and $\bm{Y}$ are commute under the multiplication. 
\end{proposition}
\textbf{Proof:}
Becuase two analogous hybrid spectrum operators $\bm{X}$ and $\bm{Y}$ are given by the relationship $\bm{X} \propto^{\bm{I}}_{[c_i], c(\lambda)} \bm{Y}$, we have
\begin{eqnarray}\label{eq1: prop2: Analogous Operators h}
\bm{X}&=&\left(\sum_{\lambda_i \in \sigma_d(\bm{X})}\lambda_i \bm{F}_{\bm{X}}(\lambda_i)+\int\limits_{\lambda \in \sigma_c(\bm{X})}\lambda d\bm{E}_{\bm{X}}(\lambda)\right.\nonumber \\
&&\left.+
\sum_{\lambda_i \in \sigma_d(\bm{X})}\left(\bm{X}-\lambda_i\bm{I}\right)\bm{F}_{\bm{X}}(\lambda_i)+\int\limits_{\lambda \in \sigma_c(\bm{X})}\left(\bm{X}-\lambda\bm{I}\right)d\bm{E}_{\bm{X}}(\lambda),\right), 
\end{eqnarray}   
and
\begin{eqnarray}\label{eq2: prop2: Analogous Operators h}
\bm{Y}&=&\bm{I}\left(\sum_{\lambda_i \in \sigma_d(\bm{X})}c_i \lambda_i \bm{F}_{\bm{X}}(\lambda_i)+\int\limits_{\lambda \in \sigma_c(\bm{X})}c(\lambda) \lambda d\bm{E}_{\bm{X}}(\lambda)\right.\nonumber \\
&&\left.+
\sum_{\lambda_i \in \sigma_d(\bm{X})}\left(\bm{X}-c_i \lambda_i\bm{I}\right)\bm{F}_{\bm{X}}(\lambda_i)+\int\limits_{\lambda \in \sigma_c(\bm{X})}\left(\bm{X}-c(\lambda)\lambda\bm{I}\right)d\bm{E}_{\bm{X}}(\lambda),\right)\bm{I}^{-1}.
\end{eqnarray}   

Moreover, from Lemma~\ref{lma: orthogonal of spectral basis} and Lemma~\ref{lma: orthogonal of spectral basis h}, we have 
\begin{eqnarray}\label{eq4: prop2: Analogous Operators}
\bm{X}\bm{Y}&=&\int\limits_{\lambda \in \sigma_c(\bm{X})}c(\lambda)\lambda^2 d\bm{E}_{\bm{X}}(\lambda)+\int\limits_{\lambda \in \sigma_c(\bm{X})}\lambda\left(\bm{X}-c(\lambda)\lambda\bm{I}\right)d\bm{E}_{\bm{X}}(\lambda)\nonumber \\
&&+\int\limits_{\lambda \in \sigma_c(\bm{X})}c(\lambda)\lambda\left(\bm{X}-\lambda\bm{I}\right)d\bm{E}_{\bm{X}}(\lambda)+\int\limits_{\lambda \in \sigma_c(\bm{X})}\left(\bm{X}-\lambda\bm{I}\right)\left(\bm{X}-c(\lambda)\lambda\bm{I}\right)d\bm{E}_{\bm{X}}(\lambda) \nonumber \\
&&+\sum\limits_{\lambda_i \in \sigma_d(\bm{X})}c_i\lambda_i^2 \bm{F}_{\bm{X}}(\lambda_i)+\sum\limits_{\lambda_i \in \sigma_d(\bm{X})}\lambda_i\left(\bm{X}-c_i\lambda_i\bm{I}\right)\bm{F}_{\bm{X}}(\lambda_i)\nonumber \\
&&+\sum\limits_{\lambda_i \in \sigma_d(\bm{X})}c_i\lambda_i\left(\bm{X}-\lambda_i\bm{I}\right)\bm{F}_{\bm{X}}(\lambda_i)+\sum\limits_{\lambda_i \in \sigma_d(\bm{X})}\left(\bm{X}-\lambda_i\bm{I}\right)\left(\bm{X}-c_i\lambda_i\bm{I}\right)\bm{F}_{\bm{X}}(\lambda_i) \nonumber \\
&=&  \bm{Y}\bm{X} 
\end{eqnarray}   
$\hfill\Box$

The next proposition is to show that the proposed analogous relationship between $\bm{X}$ and $\bm{Y}$ provided by $\bm{X} \propto^{\bm{U}}_{[c_i], c(\lambda)} \bm{Y}$ is an equivalence relation.

\begin{proposition}\label{prop3: Analogous Operators}
The relationship $\bm{X} \propto^{\bm{U}}_{[c_i], c(\lambda)} \bm{Y}$ is an equivalence relation.
\end{proposition}
\textbf{Proof:}
Given the operators $\bm{X}, \bm{Y}$ and $\bm{Z}$ with the following formats:
\begin{eqnarray}\label{eq1: prop3: Analogous Operators inf}
\bm{X}&=&\left(\sum_{\lambda_i \in \sigma_d(\bm{X})}\lambda_i \bm{F}_{\bm{X}}(\lambda_i)+\int\limits_{\lambda \in \sigma_c(\bm{X})}\lambda d\bm{E}_{\bm{X}}(\lambda)\right.\nonumber \\
&&\left.+
\sum_{\lambda_i \in \sigma_d(\bm{X})}\left(\bm{X}-\lambda_i\bm{I}\right)\bm{F}_{\bm{X}}(\lambda_i)+\int\limits_{\lambda \in \sigma_c(\bm{X})}\left(\bm{X}-\lambda\bm{I}\right)d\bm{E}_{\bm{X}}(\lambda),\right), \nonumber \\
\bm{Y}&=&\bm{U}\left(\sum_{\lambda_i \in \sigma_d(\bm{X})}c_{1,i}\lambda_i \bm{F}_{\bm{X}}(\lambda_i)+\int\limits_{\lambda \in \sigma_c(\bm{X})}c_1(\lambda)\lambda d\bm{E}_{\bm{X}}(\lambda)\right.\nonumber \\
&&\left.+
\sum_{\lambda_i \in \sigma_d(\bm{X})}\left(\bm{X}-c_{1,i}\lambda_i\bm{I}\right)\bm{F}_{\bm{X}}(\lambda_i)+\int\limits_{\lambda \in \sigma_c(\bm{X})}\left(\bm{X}-c_1(\lambda)\lambda\bm{I}\right)d\bm{E}_{\bm{X}}(\lambda),\right)\bm{U}^{-1},\nonumber \\
\bm{Z}&=&\bm{V}\left(\sum_{\lambda_i \in \sigma_d(\bm{X})}c_{2,i}\lambda_i \bm{F}_{\bm{X}}(\lambda_i)+\int\limits_{\lambda \in \sigma_c(\bm{X})}c_2(\lambda)\lambda d\bm{E}_{\bm{X}}(\lambda)\right.\nonumber \\
&&\left.+
\sum_{\lambda_i \in \sigma_d(\bm{X})}\left(\bm{X}-c_{2,i}\lambda_i\bm{I}\right)\bm{F}_{\bm{X}}(\lambda_i)+\int\limits_{\lambda \in \sigma_c(\bm{X})}\left(\bm{X}-c_2(\lambda)\lambda\bm{I}\right)d\bm{E}_{\bm{X}}(\lambda),\right)\bm{V}^{-1}
\end{eqnarray} 
we have 
\begin{eqnarray}\label{eq2s: prop3: Analogous Operators inf}
\bm{X}&\propto^{\bm{I}}_{c_{1,i}=1, c_1(\lambda)=1}&\bm{X}, ~(\mbox{reflexive})\nonumber \\
\bm{X}&\propto^{\bm{U}}_{c_{1,i},c_1(\lambda)}&\bm{Y}, \mbox{~and $\bm{Y}\propto^{\bm{U}^{-1}}_{c_{1,i}^{-1}, c_1(\lambda)^{-1}}\bm{X}$}~(\mbox{symmetric})\nonumber \\
\bm{X}&\propto^{\bm{U}}_{c_{1,i},c_1(\lambda)}&\bm{Y}, \mbox{~and $\bm{Y}\propto^{\bm{U}^{-1}\bm{V}}_{c_{2,i}/c_{1,i}, c_2(\lambda)/c_1(\lambda)}\bm{Z}$}\nonumber \\
&\Longrightarrow& \bm{X}\propto^{\bm{V}}_{c_{2,i}, c_2(\lambda)}\bm{Z}~(\mbox{transitive})
\end{eqnarray} 
$\hfill\Box$

The next proposition is to show the Fredholm determinant of a operator $\bm{X}$, denoted by 
$\mbox{det}_{\mathrm{F}}(\bm{X})$, of the trace-class operators, a special case of spectral operators, for two analogous operators. 

\begin{proposition}\label{prop4: Analogous Operators}
Given two analogous  trace-class operators $\bm{X}$ and $\bm{Y}$ provided by the relationship \\ $\bm{X} \propto^{\bm{U}}_{[c_i], i=1,2,\ldots} \bm{Y}$ with
\begin{eqnarray}\label{eq1: prop4: Analogous Operators inf}
\bm{X}&=&\left(\sum\limits_{i=1}^{\infty}\lambda_i \bm{E}_{\bm{X}}(\lambda_i)+\sum\limits_{i=1}^{\infty}\left(\bm{X}-\lambda_i\bm{I}\right)\bm{E}_{\bm{X}}(\lambda_i)\right), \nonumber \\
\bm{Y}&=& \bm{U}\left(\sum\limits_{i=1}^{\infty}c_i\lambda_i \bm{E}_{\bm{X}}(\lambda_i)+\sum\limits_{i=1}^{\infty}\left(\bm{X}-c_i\lambda_i\bm{I}\right)\bm{E}_{\bm{X}}(\lambda_i)\right)\bm{U}^{-1};
\end{eqnarray}   
we have 
\begin{eqnarray}\label{eq2: prop4: Analogous Operators inf}
\mbox{det}_{\mathrm{F}}(\bm{Y})&=&\prod\limits_{i=1}^{\infty}(1+c_i \lambda_i) \nonumber \\
&\define& \mbox{det}_{\mathrm{F}}(\bm{X}) \otimes [c_i]
\end{eqnarray}
\end{proposition}
\textbf{Proof:}
Because we have
\begin{eqnarray}\label{eq3: prop4: Analogous Operators inf}
\mbox{det}_{\mathrm{F}}(\bm{Y})&=&\prod\limits_{i=1}^{\infty}(1+c_i \lambda_i),
\end{eqnarray}
and 
\begin{eqnarray}\label{eq4: prop4: Analogous Operators inf}
\mbox{det}_{\mathrm{F}}(\bm{Y})&=&\prod\limits_{i=1}^{\infty}(1+\lambda_i),
\end{eqnarray}
this proposition is proved by comparing Eq.~\eqref{eq3: prop4: Analogous Operators inf} and Eq.~\eqref{eq4: prop4: Analogous Operators inf}.
$\hfill\Box$

The next proposition is to show the trace relationship for two trace-class operators.

\begin{proposition}\label{prop5: Analogous Operators inf}
Given two analogous trace-class operators $\bm{X}$ and $\bm{Y}$ provided by the relationship $\bm{X} \propto^{\bm{U}}_{[c, c, \ldots,]} \bm{Y}$ with
\begin{eqnarray}\label{eq1: prop5: Analogous Operators inf}
\bm{X}&=&\left(\sum\limits_{i=1}^{\infty}\lambda_i \bm{E}_{\bm{X}}(\lambda_i)+\sum\limits_{i=1}^{\infty}\left(\bm{X}-\lambda_i\bm{I}\right)\bm{E}_{\bm{X}}(\lambda_i)\right), \nonumber \\
\bm{Y}&=& \bm{U}\left(\sum\limits_{i=1}^{\infty}c_i\lambda_i \bm{E}_{\bm{X}}(\lambda_i)+\sum\limits_{i=1}^{\infty}\left(\bm{X}-c_i\lambda_i\bm{I}\right)\bm{E}_{\bm{X}}(\lambda_i)\right)\bm{U}^{-1};
\end{eqnarray}   
we have 
\begin{eqnarray}\label{eq2: prop5: Analogous Operators inf}
\mathrm{Trace}(\bm{X})&=&c \mathrm{Trace}(\bm{Y})
\end{eqnarray}
\end{proposition}
\textbf{Proof:}
Because we have
\begin{eqnarray}\label{eq3: prop5: Analogous Operators inf}
\mathrm{Trace}(\bm{X})&=&\sum\limits_{i=1}^\infty \lambda_i,
\end{eqnarray}
and 
\begin{eqnarray}\label{eq4: prop5: Analogous Operators inf}
\mathrm{Trace}(\bm{Y})&=&\sum\limits_{i=1}^\infty c\lambda_i,
\end{eqnarray}
this proposition is proved by comparing Eq.~\eqref{eq3: prop5: Analogous Operators inf} and Eq.~\eqref{eq4: prop5: Analogous Operators inf}.
$\hfill\Box$

\subsection{Spectral Mapping Theorem for Multivariate Hybrid Operators}\label{sec: Spectral Mapping Theorem for Multivariate Hybrid Operators}

We first consider spectral mapping theorem for two input  hybrid operators in Theorem~\ref{thm: Spectral Mapping Theorem for Two Variables inf h}.

\begin{theorem}\label{thm: Spectral Mapping Theorem for Two Variables inf h}
Given an analytic function $f(z_1,z_2)$ within the domain for $|z_1| < R_1$ and $|z_2| < R_2$, the first hybrid operator $\bm{X}_1$ decomposed by:
\begin{eqnarray}\label{eq1-1: thm: Spectral Mapping Theorem for Two Variables inf h}
\bm{X}_1&=&\sum_{\lambda_{1,i} \in \sigma_d(\bm{X})}\lambda_{1,i} \bm{F}_{\bm{X}}(\lambda_{1,i})+\int\limits_{\lambda_1 \in \sigma_c(\bm{X})}\lambda_1 d\bm{E}_{\bm{X}}(\lambda_1)\nonumber \\
&&+
\sum_{\lambda_{1,i} \in \sigma_d(\bm{X})}\left(\bm{X}-\lambda_{1,i}\bm{I}\right)\bm{F}_{\bm{X}}(\lambda_{1,i})+\int\limits_{\lambda_1 \in \sigma_c(\bm{X})}\left(\bm{X}-\lambda_1\bm{I}\right)d\bm{E}_{\bm{X}}(\lambda_1),
\end{eqnarray}
where $\left\vert\lambda_{1,i}\right\vert<R_1$ and $\left\vert\lambda_{1}\right\vert<R_1$, and the second hybrid operator $\bm{X}_2$ decomposed by:
\begin{eqnarray}\label{eq1-2: thm: Spectral Mapping Theorem for Two Variables inf h}
\bm{X}_2&=&\sum_{\lambda_{2,i} \in \sigma_d(\bm{X})}\lambda_{2,i} \bm{F}_{\bm{X}}(\lambda_{2,i})+\int\limits_{\lambda_2 \in \sigma_c(\bm{X})}\lambda_2 d\bm{E}_{\bm{X}}(\lambda_2)\nonumber \\
&&+
\sum_{\lambda_{2,i} \in \sigma_d(\bm{X})}\left(\bm{X}-\lambda_{2,i}\bm{I}\right)\bm{F}_{\bm{X}}(\lambda_{2,i})+\int\limits_{\lambda_2 \in \sigma_c(\bm{X})}\left(\bm{X}-\lambda_2\bm{I}\right)d\bm{E}_{\bm{X}}(\lambda_2),
\end{eqnarray}
where $\left\vert\lambda_{2,i}\right\vert<R_2$ and $\left\vert\lambda_{2}\right\vert<R_2$.

Then, we have
\begin{eqnarray}
f(\bm{X}_1, \bm{X}_2)&=&\int\limits_{\lambda_1 \in \sigma_c(\bm{X}_1)}\int\limits_{\lambda_2 \in \sigma_c(\bm{X}_2)}f(\lambda_{1}, \lambda_{2})d\bm{E}_{\bm{X}_1}(\lambda_1)d\bm{E}_{\bm{X}_2}(\lambda_2) \nonumber \\
&&+\int\limits_{\lambda_1 \in \sigma_c(\bm{X}_1)}\int\limits_{\lambda_2 \in \sigma_c(\bm{X}_2)}\sum_{q_2=1}^{m_{\lambda_2}-1}\frac{f^{(-,q_2)}(\lambda_{1},\lambda_{2})}{q_2!}d\bm{E}_{\bm{X}_1}(\lambda_1)\left(\bm{X}_2-\lambda_2\bm{I}\right)^{q_2}d\bm{E}_{\bm{X}_2}(\lambda_2) \nonumber \\
&&+\int\limits_{\lambda_1 \in \sigma_c(\bm{X}_1)}\int\limits_{\lambda_2 \in \sigma_c(\bm{X}_2)}\sum_{q_1=1}^{m_{\lambda_1}-1}\frac{f^{(q_1,-)}(\lambda_{1},\lambda_{2})}{q_1!}\left(\bm{X}_1-\lambda_1\bm{I}\right)^{q_1}d\bm{E}_{\bm{X}_1}(\lambda_1)d\bm{E}_{\bm{X}_2}(\lambda_2)\nonumber \\
&&+\int\limits_{\lambda_1 \in \sigma_c(\bm{X}_1)}\int\limits_{\lambda_2 \in \sigma_c(\bm{X}_2)}\sum_{q_1=1}^{m_{\lambda_1}-1}\sum_{q_2=1}^{m_{\lambda_2}-1}\nonumber \\
&&\frac{f^{(q_1,q_2)}(\lambda_{1},\lambda_{2})}{q_1!q_2!}\left(\bm{X}_1-\lambda_1\bm{I}\right)^{q_1}d\bm{E}_{\bm{X}_1}(\lambda_1)\left(\bm{X}_2-\lambda_2\bm{I}\right)^{q_2}d\bm{E}_{\bm{X}_2}(\lambda_2) \nonumber \\
&&{\color{red}+}\int\limits_{\lambda_1 \in \sigma_c(\bm{X}_1)}\sum\limits_{\lambda_{2,i} \in \sigma_d(\bm{X}_2)}f(\lambda_{1}, \lambda_{2,i})d\bm{E}_{\bm{X}_1}(\lambda_1)\bm{F}_{\bm{X}_2}(\lambda_{2,i}) \nonumber \\
&&+\int\limits_{\lambda_1 \in \sigma_c(\bm{X}_1)}\sum\limits_{\lambda_{2,i} \in \sigma_d(\bm{X}_2)}\sum_{q_2=1}^{m_{\lambda_{2,i}}-1}\frac{f^{(-,q_2)}(\lambda_{1},\lambda_{2,i})}{q_2!}d\bm{E}_{\bm{X}_1}(\lambda_1)\left(\bm{X}_2-\lambda_{2,i}\bm{I}\right)^{q_2}\bm{F}_{\bm{X}_2}(\lambda_{2,i}) \nonumber \\
&&+\int\limits_{\lambda_1 \in \sigma_c(\bm{X}_1)}\sum\limits_{\lambda_2 \in \sigma_d(\bm{X}_2)}\sum_{q_1=1}^{m_{\lambda_1}-1}\frac{f^{(q_1,-)}(\lambda_{1},\lambda_{2,i})}{q_1!}\left(\bm{X}_1-\lambda_1\bm{I}\right)^{q_1}d\bm{E}_{\bm{X}_1}(\lambda_1)\bm{F}_{\bm{X}_2}(\lambda_{2,i})\nonumber \\
&&+\int\limits_{\lambda_1 \in \sigma_c(\bm{X}_1)}\sum\limits_{\lambda_{2,i} \in \sigma_d(\bm{X}_2)}\sum_{q_1=1}^{m_{\lambda_1}-1}\sum_{q_2=1}^{m_{\lambda_{2,i}}-1}\nonumber \\
&&\frac{f^{(q_1,q_2)}(\lambda_{1},\lambda_{2,i})}{q_1!q_2!}\left(\bm{X}_1-\lambda_1\bm{I}\right)^{q_1}d\bm{E}_{\bm{X}_1}(\lambda_1)\left(\bm{X}_2-\lambda_{2,i}\bm{I}\right)^{q_2}\bm{F}_{\bm{X}_2}(\lambda_{2,i}) \nonumber \\ 
\end{eqnarray}
\begin{eqnarray}\label{eq2: thm: Spectral Mapping Theorem for Two Variables inf h}
&&+\sum\limits_{\lambda_{1,i} \in \sigma_d(\bm{X}_1)}\int\limits_{\lambda_2 \in \sigma_c(\bm{X}_2)}f(\lambda_{1,i}, \lambda_{2})\bm{F}_{\bm{X}_1}(\lambda_{1,i})d\bm{E}_{\bm{X}_2}(\lambda_2) \nonumber \\
&&+\sum\limits_{\lambda_{1,i} \in \sigma_d(\bm{X}_1)}\int\limits_{\lambda_2 \in \sigma_c(\bm{X}_2)}\sum_{q_2=1}^{m_{\lambda_2}-1}\frac{f^{(-,q_2)}(\lambda_{1,i},\lambda_{2})}{q_2!}\bm{F}_{\bm{X}_1}(\lambda_{1,i})\left(\bm{X}_2-\lambda_2\bm{I}\right)^{q_2}d\bm{E}_{\bm{X}_2}(\lambda_2) \nonumber \\
&&+\sum\limits_{\lambda_{1,i} \in \sigma_d(\bm{X}_1)}\int\limits_{\lambda_2 \in \sigma_c(\bm{X}_2)}\sum_{q_1=1}^{m_{\lambda_{1,i}}-1}\frac{f^{(q_1,-)}(\lambda_{1,i},\lambda_{2})}{q_1!}\left(\bm{X}_1-\lambda_{1,i}\bm{I}\right)^{q_1}\bm{F}_{\bm{X}_1}(\lambda_1)d\bm{E}_{\bm{X}_2}(\lambda_2)\nonumber \\
&&+\sum\limits_{\lambda_{1,i} \in \sigma_d(\bm{X}_1)}\int\limits_{\lambda_2 \in \sigma_c(\bm{X}_2)}\sum_{q_1=1}^{m_{\lambda_{1,i}}-1}\sum_{q_2=1}^{m_{\lambda_2}-1}\nonumber \\
&&\frac{f^{(q_1,q_2)}(\lambda_{1,i},\lambda_{2})}{q_1!q_2!}\left(\bm{X}_1-\lambda_{1,i}\bm{I}\right)^{q_1}\bm{F}_{\bm{X}_1}(\lambda_{1,i})\left(\bm{X}_2-\lambda_2\bm{I}\right)^{q_2}d\bm{E}_{\bm{X}_2}(\lambda_2) \nonumber \\
&&{\color{red}+}\sum\limits_{\lambda_{1,i} \in \sigma_d(\bm{X}_1)}\sum\limits_{\lambda_{2,i} \in \sigma_d(\bm{X}_2)}f(\lambda_{1,i}, \lambda_{2,i})\bm{F}_{\bm{X}_{1}}(\lambda_{1,i})\bm{F}_{\bm{X}_2}(\lambda_{2,i}) \nonumber \\
&&+\sum\limits_{\lambda_{1,i} \in \sigma_d(\bm{X}_1)}\sum\limits_{\lambda_{2,i} \in \sigma_d(\bm{X}_2)}\sum_{q_2=1}^{m_{\lambda_{2,i}}-1}\frac{f^{(-,q_2)}(\lambda_{1,i},\lambda_{2,i})}{q_2!}\bm{F}_{\bm{X}_1}(\lambda_{1,i})\left(\bm{X}_2-\lambda_{2,i}\bm{I}\right)^{q_2}\bm{F}_{\bm{X}_2}(\lambda_{2,i}) \nonumber \\
&&+\sum\limits_{\lambda_{1,i} \in \sigma_d(\bm{X}_1)}\sum\limits_{\lambda_2 \in \sigma_d(\bm{X}_2)}\sum_{q_1=1}^{m_{\lambda_{1,i}}-1}\frac{f^{(q_1,-)}(\lambda_{1,i},\lambda_{2,i})}{q_1!}\left(\bm{X}_1-\lambda_{1,i}\bm{I}\right)^{q_1}\bm{F}_{\bm{X}_1}(\lambda_{1,i})\bm{F}_{\bm{X}_2}(\lambda_{2,i})\nonumber \\
&&+\sum\limits_{\lambda_{1,i} \in \sigma_d(\bm{X}_1)}\sum\limits_{\lambda_{2,i} \in \sigma_d(\bm{X}_2)}\sum_{q_1=1}^{m_{\lambda_{1,i}}-1}\sum_{q_2=1}^{m_{\lambda_{2,i}}-1}\nonumber \\
&&\frac{f^{(q_1,q_2)}(\lambda_{1,i},\lambda_{2,i})}{q_1!q_2!}\left(\bm{X}_1-\lambda_{1,i}\bm{I}\right)^{q_1}\bm{F}_{\bm{X}_1}(\lambda_{1,i})\left(\bm{X}_2-\lambda_{2,i}\bm{I}\right)^{q_2}\bm{F}_{\bm{X}_2}(\lambda_{2,i})
\end{eqnarray}
\end{theorem}
\textbf{Proof:}
From Lemma~\ref{lma: orthogonal of spectral basis}, we have the following relationships for the continous spectrum part of the operator $\bm{X}_1$:
\begin{eqnarray}
d\bm{E}_{\bm{X}_1}(\lambda_1)d\bm{E}_{\bm{X}_1}(\lambda_1')&=&d\bm{E}_{\bm{X}_1}(\lambda_1)\delta(\lambda_1,\lambda_1'), \nonumber 
\end{eqnarray}
\begin{eqnarray}
(\bm{X}_1-\lambda_1\bm{I})d\bm{E}_{\bm{X}_1}(\lambda_1)(\bm{X}_1-\lambda_1'\bm{I})d\bm{E}_{\bm{X}_1}(\lambda_1')=(\bm{X}_1-\lambda_1\bm{I})^2d\bm{E}_{\bm{X}_1}(\lambda_1)\delta(\lambda_1,\lambda_1'),\nonumber 
\end{eqnarray} 
\begin{eqnarray}\label{eq3-1: thm: Spectral Mapping Theorem for Two Variables inf h}
d\bm{E}_{\bm{X}_1}(\lambda_1')(\bm{X}_1-\lambda_1\bm{I})d\bm{E}_{\bm{X}_1}(\lambda_1)&=&(\bm{X}_1-\lambda_1\bm{I})d\bm{E}_{\bm{X}_1}(\lambda)d\bm{E}_{\bm{X}_1}(\lambda_1') \nonumber \\
&=&(\bm{X}_1-\lambda_1\bm{I})d\bm{E}_{\bm{X}_1}(\lambda_1)\delta(\lambda_1,\lambda_1').
\end{eqnarray} 
Similarly, we also have the following relationships for the continous spectrum part of the operator $\bm{X}_2$:
\begin{eqnarray}
d\bm{E}_{\bm{X}_2}(\lambda_2)d\bm{E}_{\bm{X}_2}(\lambda_2')&=&d\bm{E}_{\bm{X}_2}(\lambda_2)\delta(\lambda_2,\lambda_2'), \nonumber 
\end{eqnarray}
\begin{eqnarray}
(\bm{X}_2-\lambda_2\bm{I})d\bm{E}_{\bm{X}_2}(\lambda_2)(\bm{X}_2-\lambda_2'\bm{I})d\bm{E}_{\bm{X}_2}(\lambda_2')=(\bm{X}_2-\lambda_2\bm{I})^2d\bm{E}_{\bm{X}_2}(\lambda_2)\delta(\lambda_2,\lambda_2'),\nonumber 
\end{eqnarray} 
\begin{eqnarray}\label{eq3-2: thm: Spectral Mapping Theorem for Two Variables inf h}
d\bm{E}_{\bm{X}_2}(\lambda_2')(\bm{X}_2-\lambda_2\bm{I})d\bm{E}_{\bm{X}_2}(\lambda_2)&=&(\bm{X}_2-\lambda_2\bm{I})d\bm{E}_{\bm{X}_2}(\lambda)d\bm{E}_{\bm{X}_2}(\lambda_2') \nonumber \\
&=&(\bm{X}_2-\lambda_2\bm{I})d\bm{E}_{\bm{X}_2}(\lambda_2)\delta(\lambda_2,\lambda_2').
\end{eqnarray} 

From Lemma~\ref{lma: orthogonal of spectral basis h}, we have the following relationships for the discrete spectrum part of the operator $\bm{X}_1$:
\begin{eqnarray}
\bm{F}_{\bm{X}_1}(\lambda_1)\bm{F}_{\bm{X}_1}(\lambda_1')&=&\bm{F}_{\bm{X}_1}(\lambda_1)\delta(\lambda_1,\lambda_1'), \nonumber 
\end{eqnarray}
\begin{eqnarray}
(\bm{X}_1-\lambda_1\bm{I})\bm{F}_{\bm{X}_1}(\lambda_1)(\bm{X}_1-\lambda_1'\bm{I})\bm{F}_{\bm{X}_1}(\lambda_1')=(\bm{X}_1-\lambda_1\bm{I})^2\bm{F}_{\bm{X}_1}(\lambda_1)\delta(\lambda_1,\lambda_1'),\nonumber 
\end{eqnarray} 
\begin{eqnarray}\label{eq3-3: thm: Spectral Mapping Theorem for Two Variables inf h}
\bm{F}_{\bm{X}_1}(\lambda_1')(\bm{X}_1-\lambda_1\bm{I})\bm{F}_{\bm{X}_1}(\lambda_1)&=&(\bm{X}_1-\lambda_1\bm{I})\bm{F}_{\bm{X}_1}(\lambda)\bm{F}_{\bm{X}_1}(\lambda_1') \nonumber \\
&=&(\bm{X}_1-\lambda_1\bm{I})\bm{F}_{\bm{X}_1}(\lambda_1)\delta(\lambda_1,\lambda_1').
\end{eqnarray} 
\begin{eqnarray}\label{eq3-4: thm: Spectral Mapping Theorem for Two Variables inf h}
\bm{F}_{\bm{X}_1}(\lambda_i)d\bm{E}_{\bm{X}_1}(\lambda)&=&d\bm{E}_{\bm{X}_1}(\lambda)\bm{F}_{\bm{X}_1}(\lambda_i)=\bm{0}\nonumber \\
\bm{F}_{\bm{X}_1}(\lambda_i)(\bm{X}_1-\lambda\bm{I})d\bm{E}_{\bm{X}_1}(\lambda)&=&(\bm{X}_1-\lambda\bm{I})d\bm{E}_{\bm{X}_1}(\lambda)\bm{F}_{\bm{X}_1}(\lambda_i)=\bm{0}\nonumber \\
\left(\bm{X}_1-\lambda_i\bm{I}\right)\bm{F}_{\bm{X}_1}(\lambda_i)d\bm{E}_{\bm{X}_1}(\lambda)&=&d\bm{E}_{\bm{X}_1}(\lambda)\left(\bm{X}_1-\lambda_i\bm{I}\right)\bm{F}_{\bm{X}_1}(\lambda_i)=\bm{0}\nonumber \\
\left(\bm{X}_1-\lambda_i\bm{I}\right)\bm{F}_{\bm{X}_1}(\lambda_i)(\bm{X}_1-\lambda\bm{I})d\bm{E}_{\bm{X}_1}(\lambda)&=&(\bm{X}_1-\lambda\bm{I})d\bm{E}_{\bm{X}_1}(\lambda)\left(\bm{X}_1-\lambda_i\bm{I}\right)\bm{F}_{\bm{X}_1}(\lambda_i)=\bm{0}.
\end{eqnarray}

Similarly, we also have the following relationships for the discrete spectrum part of the operator $\bm{X}_2$:
\begin{eqnarray}
\bm{F}_{\bm{X}_2}(\lambda_2)\bm{F}_{\bm{X}_2}(\lambda_2')&=&d\bm{F}_{\bm{X}_2}(\lambda_2)\delta(\lambda_2,\lambda_2'), \nonumber 
\end{eqnarray}
\begin{eqnarray}
(\bm{X}_2-\lambda_2\bm{I})\bm{F}_{\bm{X}_2}(\lambda_2)(\bm{X}_2-\lambda_2'\bm{I})\bm{F}_{\bm{X}_2}(\lambda_2')=(\bm{X}_2-\lambda_2\bm{I})^2\bm{F}_{\bm{X}_2}(\lambda_2)\delta(\lambda_2,\lambda_2'),\nonumber 
\end{eqnarray} 
\begin{eqnarray}\label{eq3-5: thm: Spectral Mapping Theorem for Two Variables inf h}
\bm{F}_{\bm{X}_2}(\lambda_2')(\bm{X}_2-\lambda_2\bm{I})\bm{F}_{\bm{X}_2}(\lambda_2)&=&(\bm{X}_2-\lambda_2\bm{I})\bm{F}_{\bm{X}_2}(\lambda)\bm{F}_{\bm{X}_2}(\lambda_2') \nonumber \\
&=&(\bm{X}_2-\lambda_2\bm{I})\bm{F}_{\bm{X}_2}(\lambda_2)\delta(\lambda_2,\lambda_2').
\end{eqnarray} 
\begin{eqnarray}\label{eq3-6: thm: Spectral Mapping Theorem for Two Variables inf h}
\bm{F}_{\bm{X}_2}(\lambda_i)d\bm{E}_{\bm{X}_2}(\lambda)&=&d\bm{E}_{\bm{X}_2}(\lambda)\bm{F}_{\bm{X}_2}(\lambda_i)=\bm{0}\nonumber \\
\bm{F}_{\bm{X}_2}(\lambda_i)(\bm{X}_2-\lambda\bm{I})d\bm{E}_{\bm{X}_2}(\lambda)&=&(\bm{X}_2-\lambda\bm{I})d\bm{E}_{\bm{X}_2}(\lambda)\bm{F}_{\bm{X}_2}(\lambda_i)=\bm{0}\nonumber \\
\left(\bm{X}_2-\lambda_i\bm{I}\right)\bm{F}_{\bm{X}_2}(\lambda_i)d\bm{E}_{\bm{X}_2}(\lambda)&=&d\bm{E}_{\bm{X}_2}(\lambda)\left(\bm{X}_2-\lambda_i\bm{I}\right)\bm{F}_{\bm{X}_2}(\lambda_i)=\bm{0}\nonumber \\
\left(\bm{X}_2-\lambda_i\bm{I}\right)\bm{F}_{\bm{X}_2}(\lambda_i)(\bm{X}_2-\lambda\bm{I})d\bm{E}_{\bm{X}_2}(\lambda)&=&(\bm{X}_2-\lambda\bm{I})d\bm{E}_{\bm{X}_2}(\lambda)\left(\bm{X}_2-\lambda_i\bm{I}\right)\bm{F}_{\bm{X}_2}(\lambda_i)=\bm{0}.
\end{eqnarray} 

Because $f(z_1, z_2)$ is an analytic function, we have
\begin{eqnarray}\label{eq3: thm: Spectral Mapping Theorem for Two Variables inf}
f(\bm{X}_1,\bm{X}_2)&=&\sum\limits_{\ell_1=0,\ell_2=0}^{\infty}a_{\ell_1,\ell_2} \bm{X}_1^{\ell_1}\bm{X}_2^{\ell_2}\nonumber \\
&=&\sum\limits_{\ell_1=0,\ell_2=0}^{\infty}a_{\ell_1,\ell_2}\left(\sum_{\lambda_{1,i} \in \sigma_d(\bm{X})}\lambda_{1,i} \bm{F}_{\bm{X}}(\lambda_{1,i})+\int\limits_{\lambda_1 \in \sigma_c(\bm{X})}\lambda_1 d\bm{E}_{\bm{X}}(\lambda_1)\right.\nonumber \\
&&\left.+
\sum_{\lambda_{1,i} \in \sigma_d(\bm{X})}\left(\bm{X}-\lambda_{1,i}\bm{I}\right)\bm{F}_{\bm{X}}(\lambda_{1,i})+\int\limits_{\lambda_1 \in \sigma_c(\bm{X})}\left(\bm{X}-\lambda_1\bm{I}\right)d\bm{E}_{\bm{X}}(\lambda_1)\right)^{\ell_1}\nonumber \\
&&\times \left(\sum_{\lambda_{2,i} \in \sigma_d(\bm{X})}\lambda_{2,i} \bm{F}_{\bm{X}}(\lambda_{2,i})+\int\limits_{\lambda_2 \in \sigma_c(\bm{X})}\lambda_2 d\bm{E}_{\bm{X}}(\lambda_2)\right.\nonumber \\
&&\left.+
\sum_{\lambda_{2,i} \in \sigma_d(\bm{X})}\left(\bm{X}-\lambda_{2,i}\bm{I}\right)\bm{F}_{\bm{X}}(\lambda_{2,i})+\int\limits_{\lambda_2 \in \sigma_c(\bm{X})}\left(\bm{X}-\lambda_2\bm{I}\right)d\bm{E}_{\bm{X}}(\lambda_2)\right)^{\ell_2}\nonumber \\
&=_1&\sum\limits_{\ell_1=0,\ell_2=0}^{\infty}a_{\ell_1,\ell_2}\left[
\sum\limits_{\lambda_{1,i} \in \sigma_d(\bm{X}_1)}\left(\lambda_{1,i}\bm{F}_{\bm{X}_1}(\lambda_{1,i})+\left(\bm{X}_1-\lambda_{1,i}\bm{I}\right)\bm{F}_{\bm{X}_1}(\lambda_{1,i})\right)^{\ell_1}\right.\nonumber \\
&&\left.+
\int\limits_{\lambda_1 \in \sigma_c(\bm{X}_1)}\left(\lambda_1 d\bm{E}_{\bm{X}_1}(\lambda_1)+\left(\bm{X}_1-\lambda_1\bm{I}\right)d\bm{E}_{\bm{X}_1}(\lambda_1)\right)^{\ell_1}\right]\nonumber \\
&&\times\left[
\sum\limits_{\lambda_{2,i} \in \sigma_d(\bm{X}_2)}\left(\lambda_{2,i}\bm{F}_{\bm{X}_2}(\lambda_{2,i})+\left(\bm{X}_2-\lambda_{2,i}\bm{I}\right)\bm{F}_{\bm{X}_2}(\lambda_{2,i})\right)^{\ell_2}\right.\nonumber \\
&&\left.+
\int\limits_{\lambda_2 \in \sigma_c(\bm{X}_2)}\left(\lambda_2 d\bm{E}_{\bm{X}_2}(\lambda_2)+\left(\bm{X}_2-\lambda_2\bm{I}\right)d\bm{E}_{\bm{X}_2}(\lambda_2)\right)^{\ell_2}\right]\nonumber \\
&=_2&\int\limits_{\lambda_1 \in \sigma_c(\bm{X}_1)}\int\limits_{\lambda_2 \in \sigma_c(\bm{X}_2)}\sum\limits_{\ell_1=0,\ell_2=0}^{\infty}\Bigg\{a_{\ell_1,\ell_2}\nonumber \\
&&\times\left[\sum\limits_{q_1=0}^{m_{\lambda_1}-1}\frac{\ell_1!}{q_1! (\ell_1-q_1)!}\lambda^{\ell_1-q_1}_{1}\left(d\bm{E}_{\bm{X}_1}(\lambda_1)\right)^{\ell_1-q_1}\left((\bm{X}_1 - \lambda_1\bm{I})d\bm{E}_{\bm{X}_1}(\lambda_1)\right)^{q_1}\right]\nonumber \\
&&\times \left[\sum\limits_{q_2=0}^{m_{\lambda_2}-1}\frac{\ell_2!}{q_2! (\ell_2-q_2)!}\lambda^{\ell_2-q_2}_{2}\left(d\bm{E}_{\bm{X}_2}(\lambda_2)\right)^{\ell_2-q_2}\left((\bm{X}_2 - \lambda_2\bm{I})d\bm{E}_{\bm{X}_2}(\lambda_2)\right)^{q_2}\right]\Bigg\}\nonumber \\
&&{\color{red}+}
\sum\limits_{\lambda_{1,i} \in \sigma_d(\bm{X}_1)}\int\limits_{\lambda_2 \in \sigma_c(\bm{X}_2)}\sum\limits_{\ell_1=0,\ell_2=0}^{\infty}\Bigg\{a_{\ell_1,\ell_2}\nonumber \\
&&\times\left[\sum\limits_{q_1=0}^{m_{\lambda_{1,i}}-1}\frac{\ell_1!}{q_1! (\ell_1-q_1)!}\lambda^{\ell_1-q_1}_{1,i}\left(\bm{F}_{\bm{X}_1}(\lambda_{1,i})\right)^{\ell_1-q_1}\left((\bm{X}_1 - \lambda_{1,i}\bm{I})\bm{F}_{\bm{X}_1}(\lambda_{1,i})\right)^{q_1}\right]\nonumber \\
&&\times \left[\sum\limits_{q_2=0}^{m_{\lambda_2}-1}\frac{\ell_2!}{q_2! (\ell_2-q_2)!}\lambda^{\ell_2-q_2}_{2}\left(d\bm{E}_{\bm{X}_2}(\lambda_2)\right)^{\ell_2-q_2}\left((\bm{X}_2 - \lambda_2\bm{I})d\bm{E}_{\bm{X}_2}(\lambda_2)\right)^{q_2}\right]\Bigg\}\nonumber
\end{eqnarray}
\begin{eqnarray}
&&{\color{red}+}\int\limits_{\lambda_1 \in \sigma_c(\bm{X}_1)}\sum\limits_{\lambda_{2,i} \in \sigma_d(\bm{X}_2)}\sum\limits_{\ell_1=0,\ell_2=0}^{\infty}\Bigg\{a_{\ell_1,\ell_2}\nonumber \\
&&\times\left[\sum\limits_{q_1=0}^{m_{\lambda_1}-1}\frac{\ell_1!}{q_1! (\ell_1-q_1)!}\lambda^{\ell_1-q_1}_{1}\left(d\bm{E}_{\bm{X}_1}(\lambda_1)\right)^{\ell_1-q_1}\left((\bm{X}_1 - \lambda_1\bm{I})d\bm{E}_{\bm{X}_1}(\lambda_1)\right)^{q_1}\right]\nonumber \\
&&\times \left[\sum\limits_{q_2=0}^{m_{\lambda_{2,i}}-1}\frac{\ell_2!}{q_2! (\ell_2-q_2)!}\lambda^{\ell_2-q_2}_{2,i}\left(\bm{F}_{\bm{X}_2}(\lambda_{2,i})\right)^{\ell_2-q_2}\left((\bm{X}_2 - \lambda_{2,i}\bm{I})\bm{F}_{\bm{X}_2}(\lambda_{2,i})\right)^{q_2}\right]\Bigg\}\nonumber \\
&&{\color{red}+}\sum\limits_{\lambda_{1,i} \in \sigma_d(\bm{X}_1)}\sum\limits_{\lambda_{2,i} \in \sigma_d(\bm{X}_2)}\sum\limits_{\ell_1=0,\ell_2=0}^{\infty}\Bigg\{a_{\ell_1,\ell_2}\nonumber \\
&&\times\left[\sum\limits_{q_1=0}^{m_{\lambda_{1,i}}-1}\frac{\ell_1!}{q_1! (\ell_1-q_1)!}\lambda^{\ell_1-q_1}_{1,i}\left(\bm{F}_{\bm{X}_1}(\lambda_{1,i})\right)^{\ell_1-q_1}\left((\bm{X}_1 - \lambda_{1,i}\bm{I})\bm{F}_{\bm{X}_1}(\lambda_{1,i})\right)^{q_1}\right]\nonumber \\
&&\times \left[\sum\limits_{q_2=0}^{m_{\lambda_{2,i}}-1}\frac{\ell_2!}{q_2! (\ell_2-q_2)!}\lambda^{\ell_2-q_2}_{2,i}\left(\bm{F}_{\bm{X}_2}(\lambda_{2,i})\right)^{\ell_2-q_2}\left((\bm{X}_2 - \lambda_{2,i}\bm{I})\bm{F}_{\bm{X}_2}(\lambda_{2,i})\right)^{q_2}\right]\Bigg\},
\end{eqnarray}
where we apply relationships given by Eqs~\eqref{eq3-1: thm: Spectral Mapping Theorem for Two Variables inf h}-\eqref{eq3-6: thm: Spectral Mapping Theorem for Two Variables inf h} in $=_1$ and $=_2$. The remaining steps in this theorem follow the same method as Theorem~\ref{thm: Spectral Mapping Theorem for Two Variables inf}.
$\hfill \Box$

We are ready to present Theorem~\ref{thm: Spectral Mapping Theorem for r Variables inf h} about spectral mapping theorem for $r$ spectral hybrid operators. 
\begin{theorem}\label{thm: Spectral Mapping Theorem for r Variables inf h}
Given an analytic function $f(z_1,z_2,\ldots,z_r)$ within the domain for $|z_l| < R_l~$ for $l=1,2,\ldots,r$, and the operator $\bm{X}_l$ decomposed by:
\begin{eqnarray}\label{eq1-1: thm: Spectral Mapping Theorem for r Variables inf h}
\bm{X}_l&=&\sum_{\lambda_{l,i} \in \sigma_d(\bm{X})}\lambda_{l,i} \bm{F}_{\bm{X}}(\lambda_{l,i})+\int\limits_{\lambda_l \in \sigma_c(\bm{X})}\lambda_l d\bm{E}_{\bm{X}}(\lambda_l)\nonumber \\
&&+
\sum_{\lambda_{l,i} \in \sigma_d(\bm{X})}\left(\bm{X}-\lambda_{l,i}\bm{I}\right)\bm{F}_{\bm{X}}(\lambda_{l,i})+\int\limits_{\lambda_l \in \sigma_c(\bm{X})}\left(\bm{X}-\lambda_l\bm{I}\right)d\bm{E}_{\bm{X}}(\lambda_l),
\end{eqnarray}
where $\left\vert\lambda_{l,i}\right\vert<R_l$ and $\left\vert\lambda_{l}\right\vert<R_l$ for $l=1,2,\ldots,r$.

Then, we have
\begin{eqnarray}\label{eq2: thm: Spectral Mapping Theorem for kappa Variables inf h}
\lefteqn{f(\bm{X}_1,\ldots,\bm{X}_r)=}\nonumber \\
&&\sum\limits_{\stirlingii{\lambda_{1,i}}{\lambda_1},\ldots,\stirlingii{\lambda_{r,i}}{\lambda_r}}\Bigg[\stirlingii{\sum\limits_{\lambda_{1,i} \in \sigma_d(\bm{X}_1)}}{\int\limits_{\lambda_1 \in \sigma_c(\bm{X}_1)}}\cdots \stirlingii{\sum\limits_{\lambda_{r,i} \in \sigma_d(\bm{X}_r)}}{\int\limits_{\lambda_r \in \sigma_c(\bm{X}_r)}}
f\left(\stirlingii{\lambda_{1,i}}{\lambda_1},\ldots,\stirlingii{\lambda_{r,i}}{\lambda_r}\right)\nonumber \\
&&\times \stirlingii{ \bm{F}_{\bm{X}_1}(\lambda_{1,i})}{d\bm{E}_{\bm{X}_1}(\lambda_1)}\cdots \stirlingii{ \bm{F}_{\bm{X}_r}(\lambda_{r,i})}{d\bm{E}_{\bm{X}_r}(\lambda_r)}\Bigg] \nonumber \\
&&+\sum\limits_{\stirlingii{\lambda_{1,i}}{\lambda_1},\ldots,\stirlingii{\lambda_{r,i}}{\lambda_r}}\Bigg[\stirlingii{\sum\limits_{\lambda_{1,i} \in \sigma_d(\bm{X}_1)}}{\int\limits_{\lambda_1 \in \sigma_c(\bm{X}_1)}}\cdots \stirlingii{\sum\limits_{\lambda_{r,i} \in \sigma_d(\bm{X}_r)}}{\int\limits_{\lambda_r \in \sigma_c(\bm{X}_r)}}\sum\limits_{\kappa=1}^{r-1} \sum\limits_{\alpha_\kappa(q_1,\ldots,q_r)}\Bigg(\sum\limits_{\alpha_{\kappa}(q_1,\ldots,q_r)=1}^{m_{\stirlingii{\lambda_{\mbox{Ind}(\alpha_{\kappa}(q_1,\ldots,q_r))},i}{\lambda_{\mbox{Ind}(\alpha_{\kappa}(q_1,\ldots,q_r))}}}-1}\nonumber \\
&&~~~~~ \frac{f^{\alpha_{\kappa}(q_1,\ldots,q_r)}\left(\stirlingii{\lambda_{1,i}}{\lambda_1},\ldots,\stirlingii{\lambda_{r,i}}{\lambda_r}\right)}{q_{\iota_1}!q_{\iota_2}!\ldots q_{\iota_\kappa}!}\nonumber \\&&\times \prod\limits_{\substack{\beta =\mbox{Ind}(\alpha_{\kappa}(q_1,\ldots,q_r)), \bm{Y}=\stirlingii{\left(\bm{X}_\beta - \lambda_{\beta,i} \bm{I}\right)^{q_\beta}\bm{F}_{\bm{X}_\beta}(\lambda_{\beta,i})}{\left(\bm{X}_\beta - \lambda_\beta \bm{I}\right)^{q_\beta}d\bm{E}_{\bm{X}_\beta}(\lambda_\beta)} \\ \beta \neq \mbox{Ind}(\alpha_{\kappa}(q_1,\ldots,q_r)), \bm{Y}=\stirlingii{\bm{F}_{\bm{X}_\beta}(\lambda_{\beta,i})}{d\bm{E}_{\bm{X}_\beta}(\lambda_\beta)}}
}^{r} \bm{Y}\Bigg)\Bigg] 
\nonumber \\
&&+\sum\limits_{\stirlingii{\lambda_{1,i}}{\lambda_1},\ldots,\stirlingii{\lambda_{r,i}}{\lambda_r}}\Bigg[\stirlingii{\sum\limits_{\lambda_{1,i} \in \sigma_d(\bm{X}_1)}}{\int\limits_{\lambda_1 \in \sigma_c(\bm{X}_1)}}\cdots \stirlingii{\sum\limits_{\lambda_{r,i} \in \sigma_d(\bm{X}_r)}}{\int\limits_{\lambda_r \in \sigma_c(\bm{X}_r)}}\sum\limits_{q_1=\ldots=q_r=1}^{m_{\stirlingii{\lambda_{1,i}}{\lambda_1}}-1,\ldots,m_{\stirlingii{\lambda_{r,i}}{\lambda_r}}-1}\nonumber \\
&&~~~~~ 
\frac{f^{(q_1,\ldots,q_r)}\left(\stirlingii{\lambda_{1,i}}{\lambda_1},\ldots,\stirlingii{\lambda_{r,i}}{\lambda_r}\right)}{q_1!\cdots q_r!}\nonumber \\
&&\times \stirlingii{\left(\bm{X}_1 - \lambda_{1,i}\bm{I}\right)^{q_1}\bm{F}_{\bm{X}_1}(\lambda_{1,i})}{\left(\bm{X}_1 - \lambda_1\bm{I}\right)^{q_1}d\bm{E}_{\bm{X}_1}(\lambda_1)} \cdots \stirlingii{\left(\bm{X}_r - \lambda_{r,i}\bm{I}\right)^{q_r}\bm{F}_{\bm{X}_r}(\lambda_{r,i})}{\left(\bm{X}_r - \lambda_r\bm{I}\right)^{q_r}d\bm{E}_{\bm{X}_r}(\lambda_r)} \Bigg],
\end{eqnarray}
where the summation $\sum\limits_{\stirlingii{\lambda_{1,i}}{\lambda_1},\ldots,\stirlingii{\lambda_{r,i}}{\lambda_r}}$ will run over all $2^r$ summands (combinations) of $\stirlingii{\lambda_{1,i}}{\lambda_1},\ldots,\stirlingii{\lambda_{r,i}}{\lambda_r}$. Each summand will be used to select corresponding term in $\stirlingii{\sum\limits_{\lambda_{l,i} \in \sigma_d(\bm{X}_1)}}{\int\limits_{\lambda_l \in \sigma_c(\bm{X}_1)}}$, $\stirlingii{\lambda_{l,i}}{\lambda_l}$, $\stirlingii{\lambda_{\mbox{Ind}(\alpha_{\kappa}(q_1,\ldots,q_r))},i}{\lambda_{\mbox{Ind}(\alpha_{\kappa}(q_1,\ldots,q_r))}}$, $\stirlingii{\left(\bm{X}_\beta - \lambda_{\beta,i} \bm{I}\right)^{q_\beta}\bm{F}_{\bm{X}_\beta}(\lambda_{\beta,i})}{\left(\bm{X}_\beta - \lambda_\beta \bm{I}\right)^{q_\beta}d\bm{E}_{\bm{X}_\beta}(\lambda_\beta)}$, $\stirlingii{\bm{F}_{\bm{X}_\beta}(\lambda_{\beta,i})}{d\bm{E}_{\bm{X}_\beta}(\lambda_\beta)}$,  $\stirlingii{ \bm{F}_{\bm{X}_l}(\lambda_{l,i})}{d\bm{E}_{\bm{X}_l}(\lambda_l)}$, and $\stirlingii{\left(\bm{X}_l - \lambda_{l,i}\bm{I}\right)^{q_l}\bm{F}_{\bm{X}_l}(\lambda_{l,i})}{\left(\bm{X}_l - \lambda_l\bm{I}\right)^{q_l}d\bm{E}_{\bm{X}_l}(\lambda_l)}$ where $l=1,2,\ldots,r$.
\end{theorem}
\textbf{Proof:}
The proof follows directly from Theorem~\ref{thm: Spectral Mapping Theorem for Two Variables inf h}, with the remaining steps consisting primarily of routine notational manipulations.
$\hfill\Box$

\bibliographystyle{IEEETran}
\bibliography{Operator_Characterization_Bib}

\end{document}